\documentclass[mnsc,nonblindrev]{informs3_hide}

\OneAndAHalfSpacedXII 

\usepackage{endnotes}
\let\footnote=\endnote

\usepackage{natbib}
 \bibpunct[, ]{(}{)}{,}{a}{}{,}%
 \def\bibfont{\small}%

\usepackage[hidelinks]{hyperref}

\def\EMAIL#1{\href{mailto:#1}{#1}}

\def\theARTICLETOP{}
\def\theARTICLETOPLEFT{}
\def\theARTICLETOPRIGHT{}
\def\theLRHFirstLine{}
\def\theLRHSecondLine{}
\def\theRRHFirstLine{}
\def\theRRHSecondLine{}

\def\theECLRHFirstLine{}
\def\theECLRHSecondLine{}
\def\theECRRHFirstLine{}
\def\theECRRHSecondLine{}
\usepackage{mathtools,amsfonts,bm,booktabs}
\usepackage{algorithm,algorithmicx,enumitem,setspace}
\usepackage[noend]{algpseudocode}
\algdef{SE}[SUBALG]{Indent}{EndIndent}{}{\algorithmicend\ }%
\algtext*{Indent}
\algtext*{EndIndent}
\usepackage{subfigure}

\newcommand{\I}[1]{I\{#1\}}

\newcommand{\bE}{\mathbb{E}}

\newcommand{\bP}{\mathbb{P}}

\newcommand{\cD}{\mathcal{D}}

\newcommand{\cF}{\mathcal{F}}
\newcommand{\cH}{\mathcal{H}}

\newcommand{\cN}{\mathcal{N}}

\newcommand{\cR}{\mathcal{R}}
\newcommand{\cT}{\mathcal{T}}
\newcommand{\cS}{\mathcal{S}}

\newcommand{\cY}{\mathcal{Y}}

\newcommand{\interior}{\operatorname{int}}
\newcommand{\cSd}{\cS_{\operatorname{br}}}
\newcommand{\cSg}{\cS_{\operatorname{part}}}
\newcommand{\gP}{g_{\operatorname{poly}}}

\TheoremsNumberedThrough     
\ECRepeatTheorems

\EquationsNumberedThrough    

\allowdisplaybreaks[4]
\usepackage{color}

\begin{document}
\RUNAUTHOR{Cao, Hu, and Hu}
\RUNTITLE{Estimation, Sensitivity Analysis, and Optimization of DRMs}
\TITLE{Estimation, Sensitivity Analysis, and Optimization of Distortion Risk Measures via Simulation}

\ARTICLEAUTHORS{%
\AUTHOR{Hao Cao}
\AFF{School of Management, Fudan University, Shanghai, 200433, China, \EMAIL{hcao21@m.fudan.edu.cn}}
\AUTHOR{Jian-Qiang Hu}
\AFF{School of Management, Fudan University, Shanghai, 200433, China, \EMAIL{hujq@fudan.edu.cn}}
\AUTHOR{Jiaqiao Hu}
\AFF{Department of Applied Mathematics and Statistics, \\
State University of New York at Stony Brook, Stony Brook, NY 11794, \EMAIL{jiaqiao.hu.1@stonybrook.edu}}
}

\ABSTRACT{We present a simulation-based method for the estimation, sensitivity analysis, and optimization of distortion risk measures (DRMs), a class of performance measures that capture a wide range of risk attitudes through an underlying distortion function. Existing simulation approaches often approximate DRMs using numerical partitions or empirical distributions. We show that, when the distortion function is piecewise polynomial, the nonlinear structure of a DRM can be handled exactly, without introducing distortion discretization error, using only simulated system outputs. This representation enables the development of a computational method that integrates finite-difference techniques for sensitivity analysis and for solving a broad class of differentiable risk-based optimization problems. Under suitable conditions, we establish strong consistency and convergence rates for DRM estimation, sensitivity analysis, and optimization. We further extend the approach to non-polynomial distortion functions through polynomial approximation. Numerical experiments on synthetic and calibrated financial examples demonstrate the effectiveness and robustness of the proposed approach across a range of applications.}

\maketitle

\section{Introduction}\label{intro}
Stochastic events such as hurricanes, earthquakes, floods, and wars can severely disrupt the operation of complex systems in areas including logistics, supply chain management, insurance, and finance. The potentially significant consequences of such events highlight the importance of effectively managing risk in stochastic systems. These systems are often characterized by nonlinear dynamics and multiple sources of uncertainty, making simulation an important tool for realistic performance evaluation and motivating growing interest in simulation-based and data-driven methods for risk analysis and performance improvement. For example, sensitivity analysis methods for value-at-risk (VaR) have been developed in a variety of simulation settings \citep{fu2009conditional, hong2009estimating, jiang2015estimating, hu2022stochastic}, with more recent work considering VaR optimization in black-box environments \citep{hu2025quantile}. Conditional value-at-risk (CVaR), another widely used risk measure, has also been extensively studied, with numerous methods developed for its estimation, sensitivity analysis, and optimization; see, e.g., \cite{rockafellar2000optimization, hong2009simulating, hu2025simulation, wang2025statistical} and references therein.

Both VaR and CVaR are special cases of distortion risk measures (DRMs), which allow a decision-maker to incorporate their own risk preferences by transforming the system's output distribution through a specified distortion function. In particular, by selecting an appropriate distortion function, DRMs can encompass a broad range of risk measures, including Glue value-at-risk (GlueVaR) \citep{belles2014beyond}, mean-Gini \citep{shalit1984mean}, Gini shortfall \citep{furman2017gini}, and DRMs with the Wang transform \citep{wang2000class}. Closely related are spectral risk measures, which can be represented as DRMs under absolutely continuous concave distortion functions \citep{acerbi2002spectral,adam2008spectral}. Stochastic optimization of empirical spectral risk measures with nondecreasing quantile weights has been studied by \cite{mehta2023stochastic}. These examples highlight the flexibility of DRMs in capturing different risk preferences. Therefore, to capture a broader range of risks, it is natural to move beyond individual risk measures and consider the more inclusive class of DRMs.


Unfortunately, the structural complexity of DRMs—particularly the discontinuity and/or nonlinearity of the distortion function—prevents the direct application or extension of the aforementioned measure-specific techniques.
\cite{glynn2021computing} proposed a simulation-based sensitivity analysis method, representing the DRM gradient using a Lebesgue-Stieltjes integral, which can be approximated through a discrete summation over a partition. For each partition point, the terms involved in the sum can be estimated using the generalized likelihood ratio technique \citep{peng2018new} applied to a batch of samples when the underlying stochastic models are known a priori.
Recently, \cite{jiang2024distortion, jiang2025stochastic} integrated this partitioning approach with a stochastic approximation (SA) scheme and developed an optimization algorithm for DRMs that effectively mitigates the computational overhead associated with batching.
However, both approaches by \cite{glynn2021computing} and \cite{jiang2025stochastic} rely on the explicit form of the underlying stochastic model and require increasingly dense partitions to reduce discretization errors, which could result in substantial computational costs.

In this paper, we consider the analysis of DRMs in a black-box simulation setting, where there is minimal or no knowledge of the underlying model generating the output. We focus on DRMs with piecewise polynomial distortion functions, a class that encompasses many widely used risk measures. A main contribution of this study is a reformulation of DRMs under this setting that reduces the nonlinear DRM functional to a simple linear combination of expectations and VaR terms, thereby eliminating the need to discretize the distortion function. This representation further allows the gradient of a DRM to be expressed as a linear combination of the gradients of these expectations and VaRs, which can be efficiently approximated using finite-difference (FD) estimators that rely solely on simulation outputs.

For differentiable DRM optimization, the gradient information can be easily incorporated into standard gradient descent methods, enabling a unified approach for the simultaneous estimation, sensitivity analysis, and optimization of DRMs. In contrast to the methods of \cite{glynn2021computing} and \cite{jiang2025stochastic}, our approach relies exclusively on output observations generated by independent runs of a black-box simulator at selected decisions. Moreover, instead of relying on progressively refined partitions to eliminate discretization errors, our approach approximates the exact DRM and its gradient using only a small number of output evaluations per iteration, with the number of evaluations proportional to the number of discontinuities and the degree of the distortion function. Also related is \cite{shao2024extreme}, which develops a robust DRM estimation approach to address distributional uncertainty in an essentially ``grey-box'' setting, whereas we develop sensitivity estimators and an optimization method for piecewise-polynomial DRMs using only simulation outputs.

Under appropriate technical conditions, we establish the strong consistency of the proposed estimators. We further derive bounds on the root mean squared errors (RMSEs) of the DRM and its gradient estimators at fixed input parameters, showing that their decay rates can be made arbitrarily close to $O(k^{-1/2})$ and $O(k^{-1/3})$, respectively, where $k$ denotes the number of algorithm iterations.
When the DRM is strongly convex with respect to the decision variables, we also show that an upper bound on the RMSE of the gradient descent iterates decays at a rate arbitrarily close to $O(k^{-1/4})$. For general non-polynomial distortion functions, we extend the proposed approach using polynomial approximation techniques, including Bernstein approximation and Chebyshev interpolation, and establish corresponding error bounds. Finally, we evaluate the proposed method on both synthetic and calibrated financial examples. The numerical results demonstrate its effectiveness, with improvements over existing partition-based and empirical DRM estimation approaches in both solution accuracy and computational efficiency.

The rest of this paper is organized as follows. In Section~\ref{sec_prob}, we introduce the problem setting. Section~\ref{sec_DRM_poly} motivates our focus on DRMs with piecewise polynomial distortion functions and develops the proposed approach for their estimation, sensitivity analysis, and optimization, along with examples of specific risk measures.
Section~\ref{sec_th} establishes the convergence and convergence rate properties of the proposed algorithm.
Section~\ref{sec_ext} extends the approach to more general DRMs using piecewise polynomial approximation techniques.
Section~\ref{sec_num} presents numerical experiments evaluating the performance of the proposed method. Finally, Section~\ref{sec_conclude} concludes the paper..

\section{Problem Setting}\label{sec_prob}
Let $Y(\theta)$ be an output performance sample from a simulation model. We assume that $Y(\theta)$ is a continuous random variable with a cumulative distribution function (CDF) $F(\cdot;\theta)$ and a continuous, strictly positive density function $f(\cdot;\theta)$ on its support $\cY(\theta)\subseteq\Re$,
where $\theta\in\Re^d$ represents the vector of model parameters (decision variables) of interest. We consider the general class of DRMs that takes the form:
\begin{align}\label{drm}
    \rho_g(\theta)
    =\int_0^\infty g\big(\bar{F}(y;\theta)\big)dy
    -\int_0^\infty\big[1-g\big(\bar{F}(-y;\theta)\big)\big]dy,
\end{align}
where $\bar{F}(y;\theta):=1-F(y;\theta)=\bP[Y(\theta)>y]$ is the complementary CDF, {$g:[0,1]\to[0,1]$ is a nondecreasing distortion function satisfying $g(0)=0$ and $g(1)=1$,}
and the integrals above are assumed to be well-defined (the detailed technical conditions are given in Section~\ref{sec_cvg}).
We also write \(\rho_g[Y(\theta)]=\rho_g(\theta)\), using square brackets for random variable inputs to emphasize that the DRM is a functional of their distributions.
For a given probability level \(\varphi\in(0,1)\), we define the \(\varphi\)-VaR by the generalized inverse \(q(\varphi;\theta):=\inf\{y:F(y;\theta)\ge\varphi\}\). Note that $q(\varphi;\theta)$ can be recovered as a special case of $\rho_g(\theta)$ by choosing $g(u)=I\{u>1-\varphi\}$, where $I\{\cdot\}$ is the indicator function.
{We consider risk estimation and sensitivity analysis, as well as the optimization problem}
\begin{equation}\label{eq_prob}
    {\min\nolimits_{\theta\in\Theta}} \rho_g(\theta),
\end{equation}
where $\Theta$ is a compact, convex, full-dimensional subset of $\Re^d$ representing the feasible region. We consider a smooth, differentiable scenario, so that the gradient information about $\rho_g(\theta)$ can be exploited to perform the optimization.

Central to our setting is that the distribution $F(\cdot;\theta)$ defining $\rho_g(\theta)$ is unknown and is instead viewed as the distribution of a simulation output random variable $Y(\theta)$. For any given $\theta$, independent replications of $Y(\theta)$ can be generated by repeatedly running the simulator, but direct evaluations of $F$, its density $f$, or their derivatives are unavailable. Thus, the underlying system is treated as a black box, and the algorithm must rely exclusively on simulation outputs generated at selected values of $\theta$. Such a setting is common in simulation-based optimization, where a complex system can be repeatedly evaluated at selected input parameters even when its input-output relationship is difficult to characterize analytically \citep{golovin2017google,cao2025blackbox,hu2025quantile}. The key question is therefore how to effectively estimate, analyze, and optimize $\rho_g(\theta)$ using only observations of $Y(\theta)$.

We begin by considering the class of distortion functions $g$ that are piecewise left-continuous with finitely many breakpoints. Specifically, we express $g$ as
\begin{equation}\label{model}
g(u)=\sum\nolimits_{\ell=1}^L g_\ell(u)I\{1-\varphi_\ell<u\le 1-\varphi_{\ell-1}\},
\qquad u\in[0,1],
\end{equation}
for some integer $L\ge1$ and constants $0=\varphi_0<\varphi_1<\cdots<\varphi_L=1$, where each $g_\ell:[1-\varphi_\ell,1-\varphi_{\ell-1}]\to[0,1]$ is continuously differentiable, with one-sided derivatives at the endpoints. Adjacent pieces need not agree at their common endpoints, thereby allowing $g$ to have jumps. We further require $g_L(0)=g(0)=0$ and $g_1(1)=g(1)=1$, with the value of $g$ at each interior breakpoint determined by the piece to its left.


Using the form of $g$ in (\ref{model}) and applying integration by parts, it can be derived that
\begin{align}
    \rho_g(\theta)
    &=\sum\nolimits_{\ell=1}^L
    \int_{q(\varphi_{\ell-1};\theta)}^{q(\varphi_\ell;\theta)}
    \big[g_\ell(\bar{F}(y;\theta))I\{y\ge0\}
    +(g_\ell(\bar{F}(y;\theta))-1)I\{y<0\}\big]dy\label{eq_rho_tmp}\\
    &\hspace{-0.8cm}=\sum\nolimits_{\ell=1}^L
    \Big[\int_{q(\varphi_{\ell-1};\theta)}^{q(\varphi_\ell;\theta)}
    yd\big[1-g_\ell\big(\bar{F}(y;\theta)\big)\big]\nonumber\\
    &+\big[g_\ell(1-\varphi_\ell)
    -I\{q(\varphi_\ell;\theta)<0\}\big]
    q(\varphi_\ell;\theta)\nonumber\\
    &-\big[g_\ell(1-\varphi_{\ell-1})
    -I\{q(\varphi_{\ell-1};\theta)<0\}\big]
    q(\varphi_{\ell-1};\theta)\Big]\nonumber\\
    &\hspace{-0.8cm}=\sum\nolimits_{\ell=1}^L
    \underbrace{\int_{q(\varphi_{\ell-1};\theta)}^{q(\varphi_\ell;\theta)}
    yd\big[1-g_\ell\big(\bar{F}(y;\theta)\big)\big]}_{\text{within intervals}}
    +\sum\nolimits_{\ell=1}^{L-1}
    \underbrace{
    \big[g_\ell(1-\varphi_\ell)-g_{\ell+1}(1-\varphi_\ell)\big]
    q(\varphi_\ell;\theta)
    }_{\text{at breakpoints}}.\label{eq_rho}
\end{align}
where we have defined $q(0;\theta):=\min\{0,\ \inf\cY(\theta)\}$, $q(1;\theta):=\max\{0,\ \sup\cY(\theta)\}$.
At each interior breakpoint $q(\varphi_\ell;\theta)$, the boundary contributions from the two adjacent intervals combine to give
\(
[g_\ell(1-\varphi_\ell)-g_{\ell+1}(1-\varphi_\ell)]q(\varphi_\ell;\theta),
\)
as the corresponding indicator terms cancel. Summing these contributions over all interior breakpoints yields the second term in \eqref{eq_rho}.

Along similar lines, the gradient of $\rho_g(\theta)$ can be obtained by differentiating both sides of Equation~\eqref{eq_rho_tmp} with respect to (w.r.t.) $\theta$ and applying the Leibniz integral rule. This yields:
\begin{equation}\label{eq_drho}
\begin{aligned}
    \nabla_\theta\rho_g(\theta)
    =&\sum\nolimits_{\ell=1}^L
    \underbrace{\int_{q(\varphi_{\ell-1};\theta)}^{q(\varphi_\ell;\theta)}
    g_\ell^{(1)}\big(\bar{F}(y;\theta)\big)
    \nabla_\theta \bar{F}(y;\theta)
    dy}_{\text{within intervals}}\\
    &\quad+\sum\nolimits_{\ell=1}^{L-1}
    \underbrace{\big[g_\ell(1-\varphi_\ell)
    -g_{\ell+1}(1-\varphi_\ell)\big]
    \nabla_\theta q(\varphi_\ell;\theta)
    }_{\text{at breakpoints}},
\end{aligned}
\end{equation}
where $g_\ell^{(1)}$ denotes the first order derivative of $g_\ell$.
It can be observed that the gradient $\nabla_\theta\rho_g(\theta)$ consists of the gradient(s) within the intervals where the distortion function $g$ is continuous, and the gradient(s) of VaRs $\nabla_\theta q(\varphi_\ell;\theta)$ at the breakpoints.

In the simulation setting considered here, however, $F$ is unknown. Consequently, neither the integrands in \eqref{eq_rho} and \eqref{eq_drho} nor the CDF and VaR gradients in \eqref{eq_drho} can be directly evaluated or estimated using methods that require explicit knowledge of the underlying distribution, such as those discussed in \cite{glynn2021computing}. Moreover, the integral components involve the unknown limits $q(\varphi_\ell;\theta)$, which must themselves be estimated. These features make a direct application of standard finite-difference methods \cite[cf., e.g.,][]{kiefer1952stochastic,spall1992multivariate} difficult.

\section{DRMs with Piecewise Polynomial Distortion Functions}\label{sec_DRM_poly}
A standard empirical estimator of a DRM is the L-statistic estimator, which replaces the underlying distribution by the empirical distribution and evaluates the resulting quantile functional using sample order statistics \citep{serfling1980approximation,shorackwellner1986empirical}. Let $\widehat\rho_g(\theta;m_k)$ denote this estimator based on $m_k$ independent simulation outputs at iteration $k$; its explicit formula is given in Section~F of the Appendix. Thus, a direct finite-difference estimator of the $j$th component of the DRM gradient is
\begin{equation*}
\frac{\widehat\rho_g(\theta+c_ke_j;m_k)-\widehat\rho_g(\theta-c_ke_j;m_k)}{2c_k},
\end{equation*}
where $e_j$ is the unit vector in the $j$th direction and $c_k>0$ is the perturbation size. For finite $m_k$, the empirical DRM estimator is generally biased, and the biases of the two perturbed evaluations propagate into the finite-difference gradient estimator. Consequently, when this gradient estimator is used in a gradient descent algorithm, the per-iteration sample size generally needs to increase, with $\lim_{k\to\infty}m_k=\infty$, to reduce this source of error as the optimization proceeds. Moreover, each iteration requires sorting the batches of simulation outputs. Under a limited simulation budget, the increasing per-iteration sample size therefore creates a tradeoff between estimation accuracy and the number of decision updates.

This motivates us to seek representations of the DRM and its gradient that lead to estimators requiring only a fixed number of simulation outputs and a fixed finite set of quantiles. To simplify notation, we suppress the dependence on $\theta$ in the following characterization. We say that a DRM $\rho_g[Y]$ of $Y$ admits an exact finite-sample output representation at the levels $\varphi_1,\ldots,\varphi_{L-1}$ if there exist an integer $m\ge1$ and a measurable function $T$, independent of the distribution of $Y$, such that, for every square-integrable $Y$,
\begin{equation}\label{eq_finite_rep}
\rho_g[Y]
=\bE\big[T\big(Y_1,\ldots,Y_m;
q(\varphi_1),\ldots,q(\varphi_{L-1})\big)\big],
\end{equation}
where $Y_1,\ldots,Y_m$ are independent and identically distributed (i.i.d.) copies of $Y$, and the expectation is finite.

The following result characterizes the class of distortion functions for which such a representation exists. It shows that the desired representation is possible precisely when the distortion function is piecewise polynomial; that is, when $g$ takes the form in \eqref{model}, with each $g_\ell$ given by 
$g_\ell(u)=\sum_{r=0}^{R_\ell}\kappa_{\ell,r}u^r$, where $R_\ell$ is the polynomial degree and $\kappa_{\ell,r}$ are constant coefficients. The proof is given in Section~A of the Appendix.

\begin{proposition}\label{prop_identification}
An exact finite-sample output representation of the form \eqref{eq_finite_rep} exists if and only if $g$ is piecewise polynomial. More precisely, if \eqref{eq_finite_rep} holds with $m$ output samples, then each $g_\ell$ has degree at most $m$. Conversely, if $R:=\max_{\ell=1,\ldots,L}R_\ell<\infty$, then \eqref{eq_finite_rep} holds with $m=\max\{1,R\}$.
\end{proposition}

\subsection{Estimation, Sensitivity Analysis and Optimization}
Proposition~\ref{prop_identification} provides the motivation for our focus on piecewise polynomial distortion functions. This class includes many widely used risk measures, such as the VaR- and Gini-type measures in Section~\ref{sec_exmp}, as well as spectral risk measures with piecewise polynomial spectral densities. In this section, we develop the estimation, sensitivity analysis, and optimization of the corresponding DRMs using only simulation outputs.

Substituting the polynomial form of \(g_\ell\) into Equation~\eqref{eq_rho} gives
\begin{align}\label{tmp}
    \rho_g(\theta)
    =&\sum_{\ell=1}^L
    \sum_{r=1}^{R_\ell} \kappa_{\ell,r}
    \int_{q(\varphi_{\ell-1};\theta)}^{q(\varphi_\ell;\theta)}
    yd\big[1-\bar{F}^r(y;\theta)\big]
    +\sum_{\ell=1}^{L-1}
    \big[g_\ell(1-\varphi_\ell)-g_{\ell+1}(1-\varphi_\ell)\big]
    q(\varphi_\ell;\theta).
\end{align}
Let \(Y_1(\theta),\ldots,Y_r(\theta)\) be \(r\) i.i.d. copies of \(Y(\theta)\), and define the sample minimum \(Y_{(1:r)}(\theta):=\min\{Y_1(\theta),\ldots,Y_r(\theta)\}\). By independence, we have
\begin{equation*}
\bP[Y_{(1:r)}(\theta)\le y]
    =1-\bP[\cap_{i=1}^r\{Y_i(\theta)>y\}]
    =1-\bar{F}^r(y;\theta),
\end{equation*}
so \(1-\bar{F}^r(y;\theta)\) is the CDF of \(Y_{(1:r)}(\theta)\). Each integral term in \eqref{tmp} can therefore be expressed as an expectation involving this sample minimum
\begin{align}\label{eq_Fr}
    &\int_{q(\varphi_{\ell-1};\theta)}^{q(\varphi_\ell;\theta)}
    yd\big[1-\bar{F}^r(y;\theta)\big]
    =\int_{q(\varphi_{\ell-1};\theta)}^{\infty}
    yd\bP[Y_{(1:r)}(\theta)\le y]
    -\int_{q(\varphi_{\ell};\theta)}^{\infty}
    yd\bP[Y_{(1:r)}(\theta)\le y] \nonumber\\
    &\qquad=\bE\big[Y_{(1:r)}(\theta)I\{Y_{(1:r)}(\theta)\ge q(\varphi_{\ell-1};\theta)\}-Y_{(1:r)}(\theta)I\{Y_{(1:r)}(\theta)\ge q(\varphi_{\ell};\theta)\}\big].
\end{align}
Combining \eqref{eq_Fr} and \eqref{tmp} yields
\begin{align}
    \rho_g(\theta)=&\sum\nolimits_{\ell=1}^L
    \sum\nolimits_{r=1}^{R_\ell} \kappa_{\ell,r}
    \underbrace{
    \bE\big[Y_{(1:r)}(\theta)I\{Y_{(1:r)}(\theta)\ge q(\varphi_{\ell-1};\theta)\}
    -Y_{(1:r)}(\theta)I\{Y_{(1:r)}(\theta)\ge q(\varphi_{\ell};\theta)\}\big]
    }_{\text{expectation of the (truncated) sample minimum}}\nonumber\\
    &\quad+\sum\nolimits_{\ell=1}^{L-1}
    \big[g_\ell(1-\varphi_\ell)-g_{\ell+1}(1-\varphi_\ell)\big]
    \underbrace{q(\varphi_\ell;\theta)}_{\text{VaR}}.\label{eq_rho_poly}
\end{align}
Equation~\eqref{eq_rho_poly} expresses the DRM exactly as a linear combination of expectations and VaRs. Note that these components can be estimated recursively using a fixed number of new simulation outputs at each iteration, rather than a batch size that grows with the iteration.

The same representation also yields a tractable formulation of the DRM gradient. Substituting the polynomial form of \(g_\ell\) into Equation~\eqref{eq_drho} and applying Taylor's theorem, the \(j\)th component of \(\nabla_{\theta}\rho_g(\theta)\) can be written as
\begin{align*}
 e_j^\intercal\nabla_\theta\rho_g(\theta)
    =&\sum\nolimits_{\ell=1}^L\sum\nolimits_{r=1}^{R_\ell}\kappa_{\ell,r}
    \int_{q(\varphi_{\ell-1};\theta)}^{q(\varphi_{\ell};\theta)}
    {\lim\nolimits_{c\to0}}
    \Big(\frac{\bar{F}^r(y;\theta+ce_j)
    -\bar{F}^r(y;\theta-ce_j)}{2c}\Big)dy\\
    &\quad+\sum\nolimits_{\ell=1}^{L-1}
    \big[g_\ell(1-\varphi_\ell)-g_{\ell+1}(1-\varphi_\ell)\big]
    e_j^\intercal\nabla_\theta q(\varphi_\ell;\theta),
\end{align*}
where the interchange of the derivative and integral can be justified by the dominated convergence theorem; sufficient conditions are provided in Section~\ref{sec_th}.

Equivalently, we can write
\begin{equation}\label{eq_drho_poly}
\begin{aligned}
    e_j^\intercal\nabla_\theta\rho_g(\theta)
    =&\sum\nolimits_{\ell=1}^L
    \sum\nolimits_{r=1}^{R_\ell} \kappa_{\ell,r}
    {\lim\nolimits_{c\to0}}
    \underbrace{\Big(\frac{\cH_{\ell,r}(\theta+ce_j)
    -\cH_{\ell,r}(\theta-ce_j)}{2c}
    \Big)}_{\text{finite difference of expectations}}\\
    &\quad+\sum\nolimits_{\ell=1}^{L-1}
    \big[g_\ell(1-\varphi_\ell)-g_{\ell+1}(1-\varphi_\ell)\big]
    \underbrace{e_j^\intercal\nabla_\theta q(\varphi_\ell;\theta)}_{\text{derivative of VaR}}.
\end{aligned}
\end{equation}
The terms in the difference quotient, with the quantile endpoints held fixed at \(\theta\), can be expressed as
\begin{align}\label{eq_H}
    \cH_{\ell,r}(\theta\pm ce_j)
    :=&\int_{q(\varphi_{\ell-1};\theta)}^{q(\varphi_\ell;\theta)}
    \big[\bar F^r(y;\theta\pm ce_j)\I{y\ge0}
    +(\bar F^r(y;\theta\pm ce_j)-1)\I{y<0}\big]dy\nonumber\\
    =&\bE[\max\{Y_{(1:r)}(\theta\pm ce_j),q(\varphi_{\ell-1};\theta)\}]
    -\bE[(Y_{(1:r)}(\theta\pm ce_j)-q(\varphi_\ell;\theta))^+]\nonumber\\
    &-q(\varphi_{\ell-1};\theta)\I{q(\varphi_{\ell-1};\theta)\ge0}
    -q(\varphi_\ell;\theta)\I{q(\varphi_\ell;\theta)<0}.
\end{align}
Here \(x^+:=\max\{x,0\}\). In the recursions below, we set \(q_k(\varphi_0)=-\infty\) and \(q_k(\varphi_L)=+\infty\), with \(\max\{Y_{(1:r)},-\infty\}=Y_{(1:r)}\) and each paired difference at \(+\infty\) interpreted as zero. The corresponding upper positive part and deterministic endpoint terms vanish. The deterministic terms cancel in the difference quotient. For finite upper endpoints, \((Y_{(1:r)}-q)^+=\max\{Y_{(1:r)},q\}-q\), giving the paired differences in \eqref{proc_drho}.

Equations~\eqref{eq_drho_poly}--\eqref{eq_H} express the DRM gradient through finite differences of expectations and gradients of VaRs. This decomposition leads to an SA procedure for the joint estimation, sensitivity analysis, and optimization of DRMs. At a given \(\theta\), we estimate \(q(\varphi_{\ell};\theta)\) and \(\nabla_{\theta}q(\varphi_{\ell};\theta)\) using the recursions of \cite{hu2025quantile}, given in \eqref{proc_q}--\eqref{proc_dq} of Algorithm~\ref{alg}. The expectation components in \eqref{eq_rho_poly} and \eqref{eq_drho_poly} are estimated using sample minima from fixed numbers of independent outputs. The resulting risk and gradient estimates are updated recursively, and the gradient estimate is incorporated into a decision update when optimization is required.

Our proposed algorithm is summarized in Algorithm~\ref{alg} below. At each iteration $k$, $\theta_k$ stands for the current decision iterate, while $q_k(\varphi)$ and $D_{q,k}(\varphi)$ denote the respective estimates of $\varphi$-VaR and its gradient. Additionally,
$\rho_k$ and $G_k$ represent the current estimates of $\rho_g(\theta)|_{\theta=\theta_k}$ and $\nabla_\theta\rho_g(\theta)|_{\theta=\theta_k}$, respectively.
At each iteration, $Y_{(1:r)}(\theta):=\min\{Y_1(\theta),\ldots,Y_r(\theta)\}$ denotes the first order statistic (minimum) of $r$ i.i.d. observations from the distribution $F(\cdot;\theta)$.
The operators \(\Pi_\cD\) and \(\Pi_\Theta\) denote Euclidean projections onto a fixed compact convex set \(\cD\) containing the true quantile gradients and onto the feasible region \(\Theta\), respectively.
Note that in sensitivity analysis, a straightforward way to estimate the gradient would involve separately applying (\ref{eq_drho_poly}) to each of the $d$ components of $\nabla_{\theta}\rho_g(\theta)$. However, this would result in a computational complexity that scales with the problem dimension, making the algorithm computationally expensive when
$d$ is large. Therefore, we instead adopt a simultaneous perturbation (SP)-based approach in (\ref{proc_drho}), following the method of \cite{spall1992multivariate}. This approach perturbs
all components of the vector $\theta_k$ simultaneously in random directions specified by a random vector $\Delta_k$, with division by $\Delta_k$ applied in an element-wise manner.

\begin{algorithm}[htb]
\caption{Estimation, sensitivity analysis, and optimization of DRMs.\protect\label{alg}}
\begin{algorithmic}[1]
\State \textbf{Input} initial estimates $\theta_0,~\rho_0,~G_0,~q_0,~D_{q,0}$; step-sizes $\{\gamma_k\}$, $\{\alpha_k\}$; perturbation-sizes $\{c_k\}$;
\State \textbf{Initialize} the iteration counter \(k\gets0\); assign each schedule undefined at \(k=0\) its \(k=1\) value;
\State \textbf{Iterate} until a stopping rule is satisfied:
\Indent
\For{$\varphi\in\{\varphi_1,\ldots,\varphi_{L-1}\}$}
{\small\begin{align}
    q_{k+1}(\varphi)=&q_k(\varphi)+\gamma_{k}
    (\varphi-I\{Y(\theta_k)\le q_k(\varphi)\}),\label{proc_q}\\
    D_{q,k+1}(\varphi)=&\Pi_\cD\Big[D_{q,k}(\varphi)
    +\frac{\gamma_k}{2c_k\Delta_k}
    \big(-I\{Y(\theta_k+c_k\Delta_k)\le q_k(\varphi)+c_k\Delta_k^\intercal D_{q,k}(\varphi)\} \nonumber\\
    &\hspace{50mm}+
    I\{Y(\theta_k-c_k\Delta_k)\le q_k(\varphi)-c_k\Delta_k^\intercal D_{q,k}(\varphi)\}\big)
    \Big];\label{proc_dq}
\end{align}}
\EndFor
\State \textbf{Estimation}
{\small\begin{align}
    \rho_{k+1}=&\rho_k+\gamma_k\Big\{\sum\nolimits_{\ell=1}^L
    \sum\nolimits_{r=1}^{R_\ell} \kappa_{\ell,r}
    \big[Y_{(1:r)}(\theta_k)I\{Y_{(1:r)}(\theta_k)\ge q_k(\varphi_{\ell-1})\}
    -Y_{(1:r)}(\theta_k)I\{Y_{(1:r)}(\theta_k)\ge q_k(\varphi_{\ell})\}\big]\nonumber\\
    &\quad\quad\quad\quad
    +\sum\nolimits_{\ell=1}^{L-1}
    \big[g_\ell(1-\varphi_\ell)-g_{\ell+1}(1-\varphi_\ell)\big]
    q_k(\varphi_\ell)
    -\rho_k\Big\};\label{proc_rho}
\end{align}}
\State \textbf{Sensitivity Analysis}
{\small\begin{align}
    G_{k+1}=&G_k+\gamma_k\Big\{
    \sum\nolimits_{\ell=1}^L
    \sum\nolimits_{r=1}^{R_\ell} \kappa_{\ell,r}
    \Big(\frac{
    \max\{Y_{(1:r)}(\theta_k+c_k\Delta_k),q_k(\varphi_{\ell-1})\}
    -\max\{Y_{(1:r)}(\theta_k-c_k\Delta_k),q_k(\varphi_{\ell-1})\}
    }{2c_k\Delta_k}\nonumber\\
    &\quad\quad\quad\quad-\frac{
    \max\{Y_{(1:r)}(\theta_k+c_k\Delta_k),q_k(\varphi_{\ell})\}
    -\max\{Y_{(1:r)}(\theta_k-c_k\Delta_k),q_k(\varphi_{\ell})\}
    }{2c_k\Delta_k}\Big)\nonumber\\
    &\quad\quad\quad\quad+\sum\nolimits_{\ell=1}^{L-1}
    \big[g_\ell(1-\varphi_\ell)-g_{\ell+1}(1-\varphi_\ell)\big]
    D_{q,k}(\varphi_\ell)
    -G_k\Big\};\label{proc_drho}
\end{align}}
\State \textbf{Optimization}
{\small\begin{align}
    \theta_{k+1}=&\Pi_\Theta\big[
    \theta_k-\alpha_kG_k
    \big],\label{proc_theta}\\
    k\gets&k+1.\nonumber
\end{align}}
\EndIndent
\end{algorithmic}
\end{algorithm}


A useful consequence of our reformulation is that, by exploiting the piecewise-polynomial structure of $g$, quantile-related quantities need only be tracked at the breakpoints $\varphi\in\cSd:=\{\varphi_1,\ldots,\varphi_{L-1}\}$.
Moreover, the quantile gradient at $\varphi_\ell$ is needed only when $g$ has a jump at that breakpoint, since its coefficient in \eqref{eq_drho_poly} is $g_\ell(1-\varphi_\ell)-g_{\ell+1}(1-\varphi_\ell)$.
For instance, in the analysis of a global polynomial distortion such as expectation or mean-Gini (see Examples~\ref{exmp_mean} and \ref{exmp_power} in Section~\ref{sec_exmp}), $\cSd$ is an empty set, meaning that recursions~\eqref{proc_q}--\eqref{proc_dq} are actually not required.
For $\varphi$-VaR, $\varphi$-CVaR, and $\varphi$-Gini shortfall (see Examples~\ref{exmp_VaR}, \ref{exmp_RVaR}, and \ref{exmp_Gini}), there is only one interior level to track. The quantile-gradient recursion \eqref{proc_dq} is needed for VaR, which has a jump, but can be omitted for the continuous CVaR and Gini-shortfall distortions. RVaR requires two quantile levels but no quantile gradients.

The partition-based approach \citep{jiang2024distortion, jiang2025stochastic} discussed in Section~\ref{intro} relies on the following analytical expressions for $\rho_g(\theta)$ and its gradient:
$\rho_g(\theta)=\int_0^1q(1-u;\theta)dg(u)
\text{~~and~~}
\nabla_\theta\rho_g(\theta)=\int_0^1\nabla_\theta q(1-u;\theta)dg(u)$. Both $\rho_g(\theta)$ and $\nabla_{\theta}\rho_g(\theta)$ are then approximated by discretizing the integral domains
into equally spaced intervals using a partition $\cSg=\{1/(2N),3/(2N),\ldots,(2N-1)/(2N)\}$ and  summing the values of the integrands at these grid points:
\begin{equation}\label{eq_partition}
\begin{aligned}
    {\rho_g}(\theta)\approx&
    \sum\nolimits_{n=1}^N
    \Big[g\Big(\frac{n}{N}\Big)-g\Big(\frac{n-1}{N}\Big)\Big]
    {q}\Big(1-\frac{2n-1}{2N};\theta\Big),\\
    {\nabla}_\theta \rho_g(\theta)\approx&
    \sum\nolimits_{n=1}^N
    \Big[g\Big(\frac{n}{N}\Big)-g\Big(\frac{n-1}{N}\Big)\Big]
    {\nabla}_\theta q\Big(1-\frac{2n-1}{2N};\theta\Big).
\end{aligned}
\end{equation}
Obviously, accurately estimating these quantities requires computing
both VaR and its gradient at a large number of increasingly dense grid points. As $N$ increases, this could become computationally demanding. Moreover, for large $N$, the values $q(1/(2N);\theta)$ and $q(1-1/(2N);\theta)$ correspond to extreme quantiles, which could be difficult to estimate reliably due to the scarcity of observations in the tails of the distribution.
In contrast, when the distortion function is piecewise polynomial, our approach leverages the structural properties of the function $g$ to estimate the integrals in their entirety. This, in effect, avoids the typical errors introduced by discretizing the integral, enabling the evaluation of both the DRM and its gradient without the computational overhead and accuracy concerns associated with dense grid discretization.

\subsection{Examples of Application}\label{sec_exmp}
\begin{figure}
    \FIGURE
    {\includegraphics[width=1\linewidth]{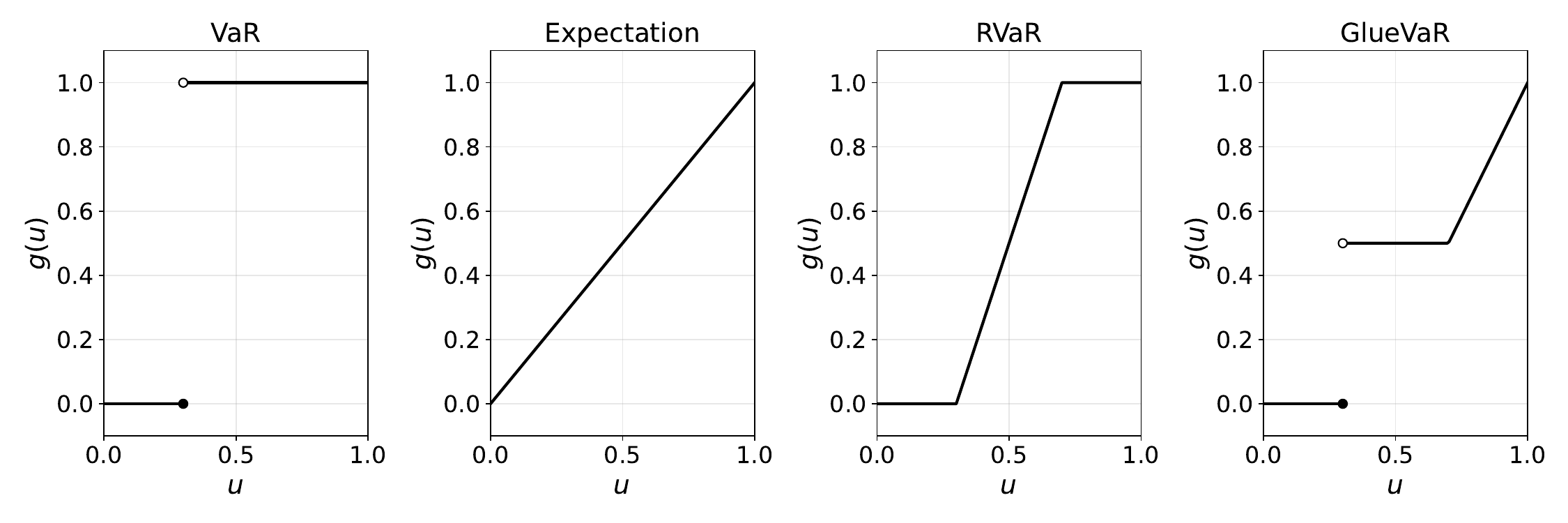}}
    {Examples of Piecewise Linear Distortion functions.\label{fig:distortion_function}}
    {}
\end{figure}

\begin{figure}
    \FIGURE
    {\includegraphics[width=1\linewidth]{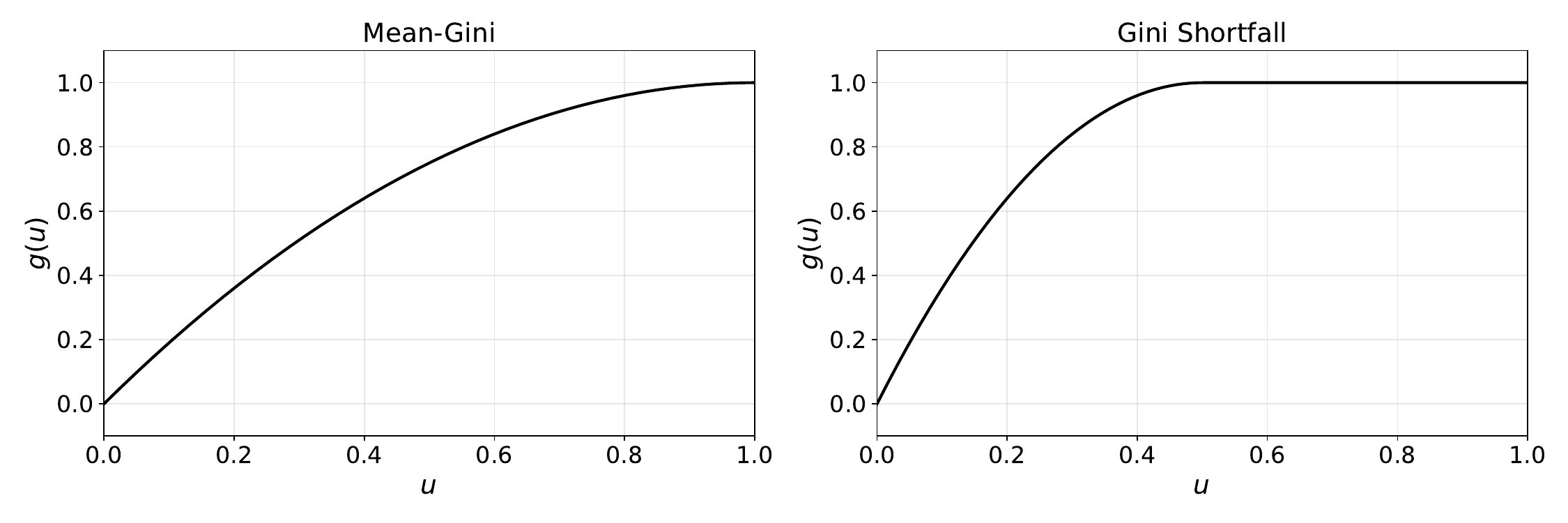}}
    {Examples of Nonlinear Distortion Functions.\label{fig:distortion_function_nonlinear}}
    {}
\end{figure}

Algorithm~\ref{alg} presents our proposed approach for DRMs with piecewise polynomial distortion functions.
To demonstrate the broad applicability of this approach, we specialize our discussion to several widely adopted DRMs, with their corresponding distortion functions shown in Figures~\ref{fig:distortion_function}--\ref{fig:distortion_function_nonlinear}. Because the optimization step (\ref{proc_theta}) is common across all risk measures, we restrict attention to the implementation of procedures~\eqref{proc_rho}--\eqref{proc_drho} for each specific measure.

\begin{example}[VaR]\label{exmp_VaR}
As mentioned at the beginning of Section~\ref{sec_prob}, when choosing
$g(u)=I\{u>1-\varphi\}$, the
$\varphi$-VaR can be recovered from (\ref{drm}) as a special case. This piecewise constant distortion function can be expressed in the form of \eqref{model} by setting $\varphi_1=\varphi$, $g_1(u)\equiv1$, $R_1=0$, $\varphi_2=1$, $g_2(u)\equiv0$, and $R_2=0$.
Therefore, procedures~\eqref{proc_rho}--\eqref{proc_drho} reduce to
$\rho_{k+1}=\rho_k+\gamma_k(q_k(\varphi)-\rho_k)\text{~~and~~}
G_{k+1}=G_k+\gamma_k(D_{q,k}(\varphi)-G_k),$
which are moving averages of the VaR and its gradient estimators, as proposed in \cite{hu2025quantile}.
\end{example}

\begin{example}[Expectation]\label{exmp_mean}
In the case of expectation, we have \(g(u)=u\). It immediately follows that procedures~\eqref{proc_rho}--\eqref{proc_drho} simplify to
\begin{equation*}
    \rho_{k+1}=\rho_k+\gamma_k(Y(\theta_k)-\rho_k)\text{~~and~~}
    G_{k+1}=G_k+\gamma_k\Big[
    \frac{Y(\theta_k+c_k\Delta_k)
    -Y(\theta_k-c_k\Delta_k)}{2c_k\Delta_k}
    -G_k\Big].
\end{equation*}
The first equation is a recursive SA procedure for estimating $\bE[Y(\theta_k) ]$, whereas the second equation is essentially the well-known SP technique \cite[cf.,][]{spall1992multivariate} applied to estimating the gradient $\nabla_{\theta}\bE[Y(\theta) ]|_{\theta=\theta_k}$.
\end{example}

\begin{example}[RVaR]\label{exmp_RVaR}
In the case of $(\varphi,\psi)$-RVaR \citep{cont2010robustness}, the distortion function is given by $g(u)=\max\{0,\allowbreak
\min\{[u-(1-\psi)]/(\psi-\varphi),1\}\}$ for $0<\varphi<\psi<1$; see Figure~\ref{fig:distortion_function}. 
The function can be written in the form of \eqref{model}
by setting
$\varphi_1=\varphi$, $g_1(u)\equiv1$,
$\varphi_2=\psi$, $g_2(u)=[u-(1-\psi)]/(\psi-\varphi)$,
$\varphi_3=1$, and $g_3(u)\equiv0$.
Consequently, procedures~\eqref{proc_rho}--\eqref{proc_drho} become:
\begin{equation*}
\begin{aligned}
    \rho_{k+1}=&\rho_k+\gamma_k\Big[
    \frac{Y(\theta_k)I\{Y(\theta_k)\ge q_k(\varphi)\}
    -Y(\theta_k)I\{Y(\theta_k)\ge q_k(\psi)\}
    }{\psi-\varphi}
    -\rho_k\Big],\\
    G_{k+1}=&G_k+\gamma_k\Big[
    \frac{\max\{Y(\theta_k+c_k\Delta_k),q_k(\varphi)\}
    -\max\{Y(\theta_k-c_k\Delta_k),q_k(\varphi)\}}{2(\psi-\varphi)c_k\Delta_k}\\
    &\hspace{15mm}
    -\frac{\max\{Y(\theta_k+c_k\Delta_k),q_k(\psi)\}
    -\max\{Y(\theta_k-c_k\Delta_k),q_k(\psi)\}}{2(\psi-\varphi)c_k\Delta_k}-G_k\Big].
\end{aligned}
\end{equation*}

A special case of interest occurs when $\psi=1$, in which case the $(\varphi,\psi)$-RVaR reduces to the $\varphi$-CVaR. Direct substitution of $\psi=1$ into the above recursions yields
\begin{align*}
    \rho_{k+1}=&\rho_k+\gamma_k\Big[
    \frac{Y(\theta_k)I\{Y(\theta_k)\ge q_k(\varphi)\}}{1-\varphi}
    -\rho_k\Big],\\
    G_{k+1}=&G_k+\gamma_k\Big[
    \frac{\max\{Y(\theta_k+c_k\Delta_k),q_k(\varphi)\}
    -\max\{Y(\theta_k-c_k\Delta_k),q_k(\varphi)\}}{2(1-\varphi)c_k\Delta_k}
    -G_k\Big].
\end{align*}
 The former equation provides a novel numerical procedure for   {\(\varphi\)} -CVaR estimation, while the latter becomes a moving average of the CVaR gradient estimator proposed in \cite{hu2025simulation}.
 \end{example}

\begin{example}[GlueVaR]\label{exmp_GlueVaR}
As another example, the $(\varphi,\psi,\varpi)$-GlueVaR \citep{belles2014beyond} can be viewed as  a DRM by defining the distortion function $g(u):=\varpi\I{1-\psi<u\le 1-\varphi}+[(1-\varpi)(u-1)+\varphi]/\varphi\I{u>1-\varphi}$ for $0<\varphi<\psi<1$ and $0\le\varpi\le1$, which is piecewise linear (see Figure~\ref{fig:distortion_function}). In the form of \eqref{model}, this corresponds to the choices:
$\varphi_1=\varphi$, $g_1(u)=[(1-\varpi)(u-1)+\varphi]/\varphi$,
$\varphi_2=\psi$, $g_2(u)\equiv\varpi$,
$\varphi_3=1$, and $g_3(u)\equiv0$.
It is not difficult to see that procedures~\eqref{proc_rho}--\eqref{proc_drho} become
\begin{align*}
    \rho_{k+1}=&\rho_k+\gamma_k\Big[
    (1-\varpi)\frac{Y(\theta_k)I\{Y(\theta_k)\le q_k(\varphi)\}}{\varphi}
    +\varpi q_k(\psi)-\rho_k\Big],\\
    G_{k+1}=&G_k+\gamma_k\Big[
    (1-\varpi)\frac{
    \min\{Y(\theta_k+c_k\Delta_k),q_k(\varphi)\}
    -\min\{Y(\theta_k-c_k\Delta_k),q_k(\varphi)\}
    }{2\varphi c_k\Delta_k}\\
    &\hspace{15mm}+\varpi D_{q,k}(\psi)-G_k\Big].
\end{align*}
These recursions provide a new computational approach for GlueVaR and its gradient estimation.
\end{example}

In addition to the above DRMs with piecewise linear distortion functions, the proposed approach can also be readily applied to analyze DRMs with nonlinear piecewise polynomial distortion functions; see Figure~\ref{fig:distortion_function_nonlinear}.

\begin{example}[Mean-Gini and Dual-Power]\label{exmp_power}
The distortion function for mean-Gini \citep{shalit1984mean} is given by $g(u)=u+\lambda(u-u^2),$ where $\lambda\in[0,1]$ is a risk-averse parameter. For its estimation and sensitivity analysis, \eqref{proc_rho}--\eqref{proc_drho} become
\begin{align*}
    \rho_{k+1}=&\rho_k+\gamma_k\big[
    Y(\theta_k)
    +\lambda(Y(\theta_k)-Y_{(1:2)}(\theta_k))
    -\rho_k\big],\\
    G_{k+1}=&G_k+\gamma_k\Big[
    (1+\lambda)\frac{Y(\theta_k+c_k\Delta_k)-Y(\theta_k-c_k\Delta_k)}{2c_k\Delta_k}\\
    &\hspace{15mm}-\lambda\frac{Y_{(1:2)}(\theta_k+c_k\Delta_k)-Y_{(1:2)}(\theta_k-c_k\Delta_k)}{2c_k\Delta_k}
    -G_k\Big].
\end{align*}
Note that in the procedures above, $Y(\theta_k)$ and $(Y(\theta_k)-Y_{(1:2)}(\theta_k))$ are the respective unbiased estimators of the mean $\bE[Y(\theta_k)]$ and the \textit{Gini deviation}:
\begin{equation*}
    \operatorname{GiniDev}
    :=\frac{1}{2}\bE[|Y_1(\theta_k)-Y_2(\theta_k)|]
    =\frac{1}{2}\bE[Y_{(2:2)}(\theta_k)-Y_{(1:2)}(\theta_k)]
    =\bE[Y(\theta_k)]-\bE[Y_{(1:2)}(\theta_k)],
\end{equation*}
where $Y_{(2:2)}(\theta_k):=\max\{Y_1(\theta_k),Y_2(\theta_k) \}$, and the last step holds because
$2\bE[Y(\theta_k)]=\bE[Y_1(\theta_k)+Y_2(\theta_k)]=\bE[Y_{(2:2)}(\theta_k)+Y_{(1:2)}(\theta_k)].$

For integer \(R\ge1\), the dual-power distortion function \citep{wang1998actuarial} \(g(u)=1-(1-u)^R\) is also a polynomial. This implies that $1-g(\bar{F}(y;\theta))=F^R(y;\theta)=\bP[Y_{(R:R)}(\theta)\le y]$, where $Y_{(R:R)}(\theta):=\max\{Y_1(\theta),\ldots,Y_R(\theta)\}$ denotes the maximum of $R$ i.i.d. copies of $Y(\theta)$.
As a consequence of Equations~\eqref{eq_rho_poly}--\eqref{eq_H}, the dual-power DRM and its sensitivity at a given $\theta_k$ can be estimated through
\begin{align*}
    \rho_{k+1}=&\rho_k+\gamma_k[
    Y_{(R:R)}(\theta_k)-\rho_k],\\
    G_{k+1}=&G_k+\gamma_k\Big[
    \frac{Y_{(R:R)}(\theta_k+c_k\Delta_k)-Y_{(R:R)}(\theta_k-c_k\Delta_k)}{2c_k\Delta_k}
    -G_k\Big],
\end{align*}
{where the resulting DRM coincides with the $R$-extremile and recovers its expected-maximum representation
\citep{daouia2019extremiles, mehta2023stochastic}.}
The mean-Gini and dual-power DRMs become equivalent when selecting $\lambda=1$ and $R=2$.
\end{example}

\begin{example}[Gini Shortfall]\label{exmp_Gini}
The distortion function for the $\varphi$-Gini shortfall \citep{furman2017gini} is given by $g(u)=u/(1-\varphi)+\lambda[(1-\varphi) u-u^2]/(1-\varphi)^2$ for $u\in[0,1-\varphi]$ and $g(u)=1$ for $u\in(1-\varphi,1]$, where $\varphi\in[0,1)$ and $\lambda\in[0,1]$ are given tail and risk-averse parameters.
Applying procedures~\eqref{proc_rho}--\eqref{proc_drho}, the Gini shortfall and its sensitivity at $\theta_k$ can be estimated through
\begin{align*}
    \rho_{k+1}=&\rho_k+\gamma_k\Big[
    \frac{Y(\theta_k)I\{Y(\theta_k)\ge q_k(\varphi)\}}{1-\varphi}\\
    &\hspace{15mm}
    +\lambda\Big(
    \frac{Y(\theta_k)I\{Y(\theta_k)\ge q_k(\varphi)\}}{1-\varphi}
    -\frac{Y_{(1:2)}(\theta_k)I\{Y_{(1:2)}(\theta_k)\ge q_k(\varphi)\}}{(1-\varphi)^2}
    \Big)-\rho_k\Big],\\
    G_{k+1}=&G_k+\gamma_k\Big[
    \frac{1+\lambda}{1-\varphi}
    \frac{\max\{Y(\theta_k+c_k\Delta_k),q_k(\varphi)\}
    -\max\{Y(\theta_k-c_k\Delta_k),q_k(\varphi)\}}{2c_k\Delta_k}\\
    &\hspace{5mm}
    -\frac{\lambda}{(1-\varphi)^2}
    \frac{\max\{Y_{(1:2)}(\theta_k+c_k\Delta_k),q_k(\varphi)\}
    -\max\{Y_{(1:2)}(\theta_k-c_k\Delta_k),q_k(\varphi)\}}{2c_k\Delta_k}
    -G_k\Big].
\end{align*}
Analogous to mean-Gini, the Gini shortfall comprises the $\varphi$-CVaR and the $\varphi$-\textit{tail Gini deviation}:
\begin{align*}
    \operatorname{Tail-GiniDev}_\varphi
    &:=\frac{1}{2}\bE[|Y_1^*-Y_2^*|]
    =\bE[Y^*]-\bE[Y_{(1:2)}^*]\\
    &=\frac{\bE[Y(\theta_k)I\{Y(\theta_k)\ge q(\varphi;\theta_k)\}]}{1-\varphi}
    -\frac{\bE[Y_{(1:2)}(\theta_k)I\{Y_{(1:2)}(\theta_k)\ge q(\varphi;\theta_k)\}]}{(1-\varphi)^2},
\end{align*}
where $Y_1^*,Y_2^*,Y^*\sim(Y(\theta_k)|Y(\theta_k)\ge q(\varphi;\theta_k))$ are i.i.d.
As $\varphi\to0$, the $\varphi$-Gini shortfall degenerates to mean-Gini, with the $\varphi$-CVaR and $\varphi$-tail Gini deviation recovering to expectation and Gini deviation, respectively.
\end{example}

\section{Theoretical Analysis}\label{sec_th}
In this section, we establish the strong consistency of the estimators proposed in Algorithm 1 and characterize their convergence rates.
Let \((\Omega,\cF,\bP)\) be the probability space on which the algorithm is defined. The filtration \(\cF_k\) contains the initial values and all directions, simulator outputs, and iterates observed before the new direction and simulator calls of iteration \(k\). Conditional on \(\cF_k\) and \(\Delta_k\), the fresh simulator replications at the selected decisions have the corresponding laws \(F(\cdot;\theta)\) and are independent within each sample minimum. Replications may be reused across risk components at the same decision.
Let $R=\max\{R_1,\ldots,R_L\}$ be the highest degree of $g_\ell$'s.
Throughout the paper, we assume the feasible region $\Theta$ is a compact convex set taking the form $\Theta=\{\theta\in\Re^d:h_j(\theta)\le0,j=1,\ldots,m\}$, for some continuously differentiable functions $h_j(\theta)$ with their gradients satisfying $\nabla_\theta h_j(\theta)\ne0$ whenever $h_j(\theta)=0$.
This formulation covers many commonly encountered feasible regions, such as hyperballs and hyperrectangles \cite[cf., e.g.,][]{kushner2003stochastic, hu2025quantile}.
For two positive functions $a(k)$ and $b(k)$, we write $a(k) = O(b(k))$ if $\limsup_{k\to \infty}a(k)/b(k)<\infty$, $a(k) = o(b(k))$ if $\lim_{k\to \infty}a(k)/b(k)=0$, and $a(k)=\mathsf{\Omega}(b(k))$ if $\liminf_{k\to\infty}a(k)/b(k)>0$.
For a given vector $\nu\in\Re^d$, let $\|\nu\|$ be its Euclidean norm and $\nu_j$ be its $j$th element.
For a given third-order tensor $\cT\in\Re^{d\times d\times d}$, let $\cT_{ijk}$ be its $(i,j,k)$-th element.
For a given set $A$, let $\interior~A$ denote its interior.
For a $J$-times continuously differentiable univariate function $g$, we denote its $J$th order derivative by $g^{(J)}$, with the convention that $g^{(0)}=g$.

\subsection{Strong Consistency of Estimators}\label{sec_cvg}
We make the following assumptions in our consistency analysis:

\begin{assumption}\label{A_distr}
Let \(\cN\) be an open neighborhood of \(\Theta\) containing all perturbed decisions used by the algorithm.
The distortion function $g$ and the CDF $F$ satisfy the following conditions:
\begin{enumerate}[label=\upshape(\alph*),ref=\theassumption (\alph*)]
\item\label{A_moment} If $R\ge1$, then
\(\sup_{\theta\in\cN}\bE[|Y(\theta)|^2]<\infty.\)
\item\label{A_conti} {For every \(y\in\Re\), the function \(F(y;\cdot)\) is three-times continuously differentiable on \(\cN\), with its first three parameter derivatives uniformly bounded over \(y\in\Re\) and \(\theta\in\cN\). For each interior breakpoint, \(F\) is also three-times continuously differentiable jointly in \((y,\theta)\), with uniformly bounded mixed derivatives up to order three, on a deterministic open region containing the true quantiles and the line segments between \((q_k(\varphi_\ell),\theta_k)\) and \((q_k(\varphi_\ell)\pm c_k\Delta_k^\intercal D_{q,k}(\varphi_\ell),\theta_k\pm c_k\Delta_k)\) used in \eqref{proc_dq};}
\item\label{A_bound}
There exists an integrable function $b(y)$ such that $\sup_{\theta\in\cN}
|\nabla_\theta^j F(y;\theta)_{i_1,\ldots,i_j}|
\le b(y)$, for \(j=1,2,3\) and all corresponding component indices.
\item\label{A_density}\label{A_density2}There exists \(\varepsilon>0\) such that, with probability one (w.p.1), \(f(y;\theta)\ge\varepsilon\) for every \(k\ge0\), \(\theta\in\cN\), and \(\ell=1,\ldots,L-1\), and all \(y\) between \(q_k(\varphi_\ell)\) and \(q(\varphi_\ell;\theta)\).
\end{enumerate}
\end{assumption}

\begin{assumption}\label{A_step}
The step- and perturbation-sizes satisfy the following conditions:
\begin{enumerate}[label=\upshape(\alph*),ref=\theassumption (\alph*)]
\item \(\gamma_k>0\), \(c_k>0\),
\(\sum_k\gamma_k^2/c_k^2<\infty\), \(\lim_{k\to\infty}c_k=0\);
\item \(\alpha_k>0\),
\(\sum_k\alpha_k=\infty\);
\item\label{A_multi} \(\alpha_k=o(\gamma_k)\).
\end{enumerate}
\end{assumption}

\begin{assumption}\label{A_direct}
    For each \(k\ge0\), \(\Delta_k\) is independent of $\cF_k$. The components of $\Delta_k$ are i.i.d. with the Bernoulli distribution $\bP[\Delta_{k,j}=1]=\bP[\Delta_{k,j}=-1]=1/2~$ for $j\in\{1,\ldots,d\}$.
\end{assumption}

Assumption~\ref{A_moment} is primarily required when a DRM involves expectation(s), ensuring that $\rho_g(\cdot)$ is well-defined over $\Theta$, which is a standard moment condition used in the literature \citep{spall1992multivariate, hong2009simulating, hu2025simulation}.
{Assumption~\ref{A_conti} supplies the third-order expansions for both the SP estimator of the sample-minimum terms and the breakpoint-quantile gradient recursion \cite[cf., e.g.,][]{spall1992multivariate, hu2025quantile}.}
Assumption~\ref{A_bound} is the regularity condition imposed for the application of Lebesgue's dominated convergence theorem and Fubini's theorem, allowing the interchange of integrals and limits.
Since \(0\le F(y;\theta)\le1\), Assumption~\ref{A_bound}, together with the uniform boundedness of the derivatives in Assumption~\ref{A_conti}, also provides the required integrable bounds for the derivatives of \(F^r(y;\theta)\) and \(\bar F^r(y;\theta)\), for every \(r=1,\ldots,R\).
Assumption~\ref{A_density} ensures contraction of the quantile and quantile-gradient recursions. As discussed in \cite{hu2025quantile}, it can be ensured by projecting the quantile iterates onto a fixed compact interval containing the relevant true quantiles in its interior, provided \(f\) is continuous and strictly positive on the product of this interval and a compact parameter neighborhood containing \(\cN\).
Together, conditions \ref{A_moment}--\ref{A_density} ensure the existence and boundedness of the DRM gradient $\nabla \rho_g(\cdot)$ on $\Theta$.
Assumption~\ref{A_step} consists of standard conditions on step-sizes \cite[cf., e.g.,][]{kushner2003stochastic, borkar2008stochastic, chau2014overview}. In particular,
Assumption~\ref{A_multi} specifies the two-timescale scheme of Algorithm~\ref{alg} through the selection of step-sizes \cite[cf., e.g.,][]{borkar2008stochastic, hu2022stochastic, hu2025quantile}.
Under this assumption, the decision recursion evolves on a
slower timescale than the estimation recursions so that the DRM and gradient estimators can track their changing targets as $\theta_k$ evolves.
Finally, the Bernoulli random directions specified in Assumption~\ref{A_direct} are commonly adopted in SP-based gradient estimators \citep{spall1992multivariate, cao2025blackbox, hu2025quantile}.

We present some preliminary results that facilitate the consistency analysis of Algorithm~\ref{alg}. Proofs of these results are provided in Section~C of the Appendix. The following result
quantifies the bias of the gradient estimators in procedure~\eqref{proc_drho}, specifically the SP technique based on (truncated) sample minimums.
\begin{lemma}\label{lm_bias}
    If Assumptions~\ref{A_moment}--\ref{A_bound} and \ref{A_step}--\ref{A_direct} hold, then, for \(r=1,\ldots,R\) and any \(\cF_k\)-measurable threshold \(q\in\Re\cup\{-\infty\}\), w.p.1,
    \begin{align*}
        \Big\|\bE\Big[\frac{\max\{Y_{(1:r)}(\theta_k+c_k\Delta_k),q\}
        -\max\{Y_{(1:r)}(\theta_k-c_k\Delta_k),q\}}{2c_k\Delta_k}\Big|\cF_k\Big]\hspace{5mm}&\nonumber\\
        -\int_q^{+\infty}\nabla_\theta \bar{F}^r(y;\theta_k)dy\Big\|
        &=O(c_k^2).
    \end{align*}
\end{lemma}
The next result establishes that both the DRM and its gradient estimators remain bounded almost surely (with probability 1) and in the second moment.
\begin{lemma}\label{lm_bound}
    If Assumptions~\ref{A_distr}--\ref{A_direct} hold, then $\{(\rho_k,G_k)\}$ generated by \eqref{proc_rho}--\eqref{proc_drho} satisfies:
    (a) $\sup_k\bE[\rho_k^2]<\infty$, $\sup_k\bE[\|G_k\|^2]<\infty$;
    (b) $\sup_k|\rho_k|<\infty$, $\sup_k\|G_k\|<\infty$ w.p.1.
\end{lemma}

Proposition~\ref{prop_cvg} below establishes the strong convergence of the DRM and its gradient estimators $\{(\rho_k,G_k)\}$.
  {For } DRM estimation and sensitivity analysis   at a fixed input parameter   {\(\theta\in\Theta\), we hold \(\theta_k\equiv\theta\) and omit the decision update~} \eqref{proc_theta}. In this   {fixed parameter } setting, the conditions involving \(\{\alpha_k\}\) in Assumption~\ref{A_step} are not needed{. Under the remaining conditions \(\gamma_k>0\), \(\sum_k\gamma_k=\infty\), \(\sum_k\gamma_k^2/c_k^2<\infty\), and \(c_k\to0\)} ,   the same tracking argument underlying Proposition~\ref{prop_cvg} establishes the strong consistency of \(\{(\rho_k,G_k)\}\).

\begin{proposition}[Estimation and sensitivity analysis]\label{prop_cvg}
    If the distortion function $g$ is piecewise polynomial and Assumptions~\ref{A_distr}--\ref{A_direct} hold, then w.p.1 we have
    \begin{align*}
        {\lim\nolimits_{k\to\infty}}|\rho_k-\rho_g(\theta_k)|=0
        \text{~~and~~}
        {\lim\nolimits_{k\to\infty}}\|G_k-\nabla_\theta\rho_g(\theta)|_{\theta=\theta_k}\|=0.
    \end{align*}
\end{proposition}

We next consider the optimization recursion.
Theorem~\ref{th_cvg} characterizes the limit points of \(\{\theta_k\}\) through the projected ordinary differential equation (ODE) in \eqref{ode_theta}.
When $\rho_g(\cdot)$ is strictly convex, the result can be strengthened to guarantee the strong convergence of the sequence to the unique optimal solution $\theta^*$ of Problem~\eqref{eq_prob}.

\begin{theorem}[Optimization]\label{th_cvg}
    {If the distortion function \(g\) is piecewise polynomial, Assumptions~\ref{A_distr}--\ref{A_direct} hold, and the set of objective values at equilibria of \eqref{ode_theta} has empty interior, then w.p.1 every limit point of \(\{\theta_k\}\) is an equilibrium point of the projected ODE}
    \begin{equation}\label{ode_theta}
        \dot{\theta}(t)=-\nabla_\theta\rho_g(\theta)|_{\theta=\theta(t)}
        +z(t),~~~~t\ge0,
    \end{equation}
    {where \(z(t)\) is the minimum normal force needed to keep the trajectory in \(\Theta\).}
    In addition, if $\rho_g(\cdot)$ is strictly convex on $\Theta$, then we further have
    \[
        {\lim\nolimits_{k\to\infty}}\|\theta_k-\theta^*\|=0,
        ~~~~w.p.1.
    \]
\end{theorem}

\subsection{Rates of Convergence}\label{sec_rate}
To further characterize the convergence rates of the proposed estimators, we
conduct a finite-time analysis by deriving error bounds on the RMSEs of the iterates generated by Algorithm \ref{alg}. Proofs of all results presented in this section are provided
in Section~D of the Appendix.

In addition to the conditions required for convergence,
we impose the following regularity conditions:
\begin{assumption}\label{A_rate}~\
\begin{enumerate}[label=\upshape(\alph*),ref=\theassumption (\alph*)]

\item\label{A_quantile_smooth}On a deterministic open region containing the line segments in Assumption~\ref{A_density}, \(f(y;\theta)\) is uniformly bounded above and jointly Lipschitz in \((y,\theta)\), and \(\nabla_\theta F(y;\theta)\) is uniformly Lipschitz in \(y\). The corresponding true quantiles \(q(\varphi_\ell;\theta)\) are twice continuously differentiable on \(\cN\), with uniformly bounded Hessians, for \(\ell=1,\ldots,L-1\).
\item\label{A_step2} The step- and perturbation-sizes follow the standard forms: $\gamma_k=C_1/k^\gamma$, $\alpha_k=C_2/k^{\alpha}$, and $c_k=C_3/k^\eta$ for constants $C_1,C_2,C_3>0$ and $\gamma,\alpha\in(1/2,1),\eta\in(0,\gamma-1/2)$.
\end{enumerate}
\end{assumption}

{Assumptions~\ref{A_density} and~\ref{A_quantile_smooth} control quantile and quantile-gradient tracking. They apply only when the distortion has interior breakpoints. The two recursions share a step size here; Lemma~24 in the Appendix proves the needed joint bounds directly, using distributional conditions comparable to those in \cite{hu2025quantile}.}
The step-size and perturbation-size conditions in Assumption~\ref{A_step2} are standard in both practice and the convergence analysis of SA algorithms
\cite[cf., e.g.,][]{spall1992multivariate}.

The following result characterizes the convergence rates of the DRM and its gradient estimators in terms of the step-size and perturbation-size parameters.
\begin{proposition}[Estimation and sensitivity analysis]\label{prop_rate}
    If the distortion function $g$ is piecewise polynomial and Assumptions~\ref{A_distr}--\ref{A_rate} hold, then $\{(\rho_k,G_k)\}$ generated by Procedures~\eqref{proc_rho}--\eqref{proc_drho} satisfies
    \begin{align*}
        \sqrt{\bE[|\rho_k-\rho_g(\theta_k)|^2]}
        =&O\bigg(\frac{\alpha_k}{\gamma_k}\bigg)
        +O(\gamma_k^{1/2}),\\
        \sqrt{\bE[\|G_k-\nabla_\theta\rho_g(\theta)|_{\theta=\theta_k}\|^2]}
        =&O\bigg(\frac{\alpha_k}{\gamma_k}\bigg)
        +O\bigg(\frac{\gamma_k^{1/2}}{c_k}\bigg)
        +O(c_k^2).
    \end{align*}
\end{proposition}

For estimation and sensitivity analysis at a fixed parameter \(\theta\in\Theta\), we again hold \(\theta_k\equiv\theta\) and omit the decision update \eqref{proc_theta}. The moving target term \(O(\alpha_k/\gamma_k)\) in Proposition~\ref{prop_rate} then vanishes, and the restrictions involving \(\alpha\) in Assumption~\ref{A_step2} are no longer relevant.  Proposition~\ref{prop_rate} therefore provides direct guidance for selecting \(\{\gamma_k\}\) and \(\{c_k\}\) in this setting .
Specifically, this result, together with the remaining conditions in Assumption~\ref{A_step2}, implies that an upper bound on the RMSE of $\{\rho_k\}$ diminishes at a rate that can be made arbitrarily close to $O(k^{-1/2})$, which is approximately attained by taking $\gamma\approx1$.
Similarly, an upper bound on the RMSE of $\{G_k\}$ diminishes at a rate that can be made arbitrarily close to $O(k^{-1/3})$ by balancing the $O(\gamma_k^{1/2}/c_k)$ and $O(c_k^2)$ terms in Proposition~\ref{prop_rate}, which is approximately attained by taking $\gamma\approx1$ and $\eta=1/6$.

For DRM optimization, we further impose the conditions and derive the convergence rate result as follows.

\begin{assumption}\label{A_rate_opt}~\
\begin{enumerate}[label=\upshape(\alph*),ref=\theassumption (\alph*)]
\item\label{A_int} Problem~\eqref{eq_prob} has a unique optimal solution $\theta^*\in\interior~\Theta$.
\item\label{A_Hessian} Let $\lambda_1(\theta)$ denote the smallest eigenvalue of the Hessian matrix $\nabla_\theta^2\rho_g(\theta)$. There exists a constant $\varepsilon_\lambda>0$ such that $\lambda_1(\theta)\ge\varepsilon_\lambda$ for all $\theta$ on the line segments connecting $\theta_k$ and $\theta^*$.
\end{enumerate}
\end{assumption}

\begin{theorem}[Optimization]\label{th_rate}
     Let the distortion function $g$ be piecewise polynomial, and suppose Assumptions~\ref{A_distr}--\ref{A_rate_opt} hold. Then the sequence $\{\theta_k\}$ generated by Algorithm~\ref{alg} satisfies
    \[
        \sqrt{\bE[\|\theta_k-\theta^*\|^2]}
        =O\bigg(\frac{\alpha_k}{\gamma_k}\bigg)
        +O\bigg(\frac{\gamma_k^{1/2}}{c_k}\bigg)
        +O(c_k^2).
    \]
\end{theorem}

Assumption~\ref{A_int} ensures that the optimal solution $\theta^*$ satisfies the first-order optimality condition $\nabla_\theta\rho_g(\theta)|_{\theta=\theta^*}=0$.
Assumption~\ref{A_Hessian} holds when $\rho_g(\theta)$ is strongly convex on $\Theta$, a common condition adopted in the rate analysis of gradient-based algorithms \cite[cf., e.g.,][]{hu2025convergence}.
As a consequence of Theorem~\ref{th_rate} and the specific forms of the algorithm parameters in Assumption~\ref{A_step2}, by balancing the three terms in the error bound, we find that an upper bound on the RMSE of $\{\theta_k\}$ diminishes at a rate that can be made arbitrarily close to $O(k^{-1/4})$.
This rate is approximately attained by taking $\alpha\approx1$, $\gamma=3/4$, and $\eta=1/8$.

\section{Extension to Non-Polynomial Distortion Functions}\label{sec_ext}
The methodology developed in Section~\ref{sec_DRM_poly} relies on the piecewise polynomial structure of the distortion function to obtain tractable representations of the DRM and its gradient. We now extend this approach to more general distortion functions that are sufficiently smooth on each subinterval in \eqref{model}. The main idea is to approximate each $g_\ell$ by a polynomial and then apply the proposed estimation and optimization procedures to the resulting approximate DRM.

One approach is to construct an approximation of $g$ using a \textit{(piecewise) Bernstein approximation with Taylor series correction} \cite[cf., e.g.,][]{orucc1999generalization, duman2025faster}:
\begin{equation*}
    \tilde{g}(u;M,J)
    :=\sum\nolimits_{\ell=1}^L\tilde{g}_\ell(u;M,J)\
    \I{1-\varphi_\ell<u\le 1-\varphi_{\ell-1}},
\end{equation*}
and on each subinterval, the local approximation is given by
\begin{equation}
    \tilde{g}_\ell(u;M,J)
    :=\sum\nolimits_{m=0}^{M}\Big[
    \sum\nolimits_{j=0}^{J}
    \frac{g_\ell^{(j)}(b_{\ell,m})}{j!}(u-b_{\ell,m})^j
    \Big]\binom{M}{m}
    \frac{(u-b_{\ell,0})^m
    (b_{\ell,M}-u)^{M-m}}{(b_{\ell,M}-b_{\ell,0})^M}, \label{eq_Bernstein}
\end{equation}
where \(M\) denotes the degree of the Bernstein basis,
$J$ denotes the degree of Taylor polynomial,
and the partition points are defined as $b_{\ell,m}:=(1-\varphi_{\ell})+(m/M)(\varphi_{\ell}-\varphi_{\ell-1})$ for $m=0,\ldots,M$.
When \(J=0\), the approximation in \eqref{eq_Bernstein} reduces to
the standard Bernstein approximation, which preserves monotonicity
and convexity or concavity; see Proposition~\ref{prop_Bernstein}.
For \(J\ge1\), the Taylor--Bernstein correction incorporates local
derivative information and can provide greater approximation accuracy for a given polynomial degree.
This increased approximation accuracy comes at the cost of the shape-preserving property of the standard Bernstein approximation, as monotonicity is not guaranteed in general.
A piecewise-linear approximation is included as a special case of this construction by applying the standard Bernstein approximation with \(M=1\) on each subinterval.
Increasing its approximation accuracy requires introducing additional breakpoints, whereas higher-order polynomial approximations can improve the approximation without increasing the number of pieces.

\begin{figure}
    \FIGURE
    {\includegraphics[width=1\linewidth]{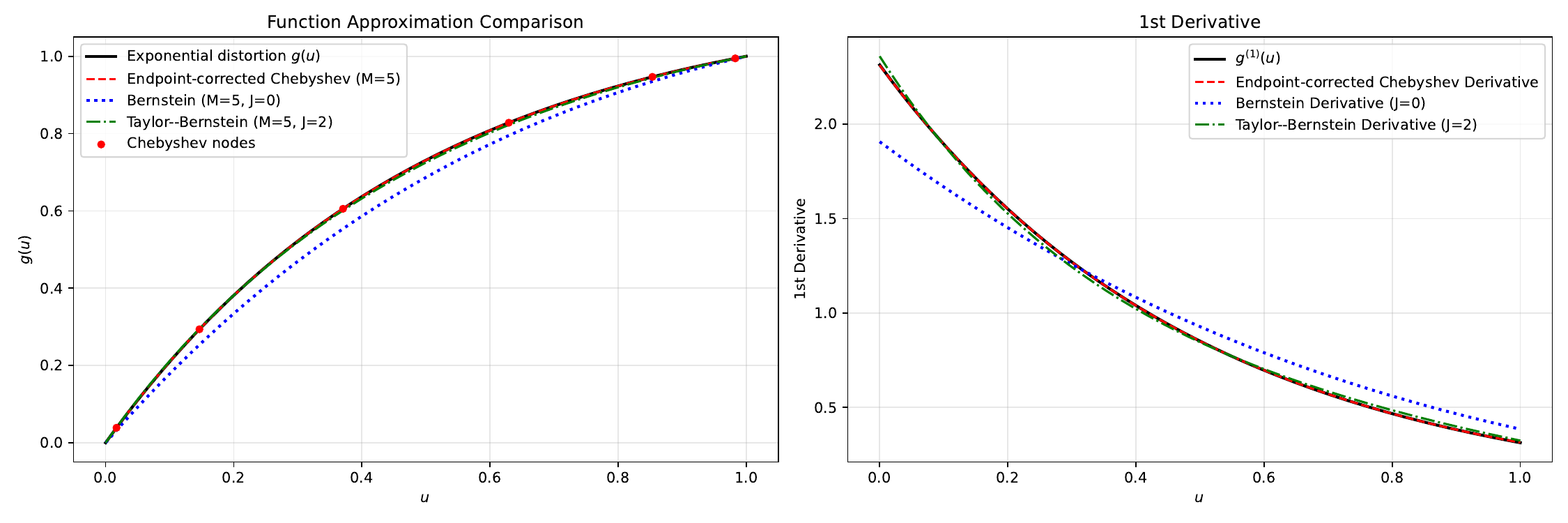}}
    {Illustration of Polynomial Approximation Methods.\label{fig:function_approx}}
    {}
\end{figure}

The Taylor--Bernstein approximation in \eqref{eq_Bernstein} can be readily replaced by other polynomial approximation techniques. Another natural alternative is the (piecewise) Chebyshev approximation \cite[cf.,][]{burden2010numerical}.
We use the following endpoint-corrected version:
\begin{align}
    \Breve{g}(u;M)
    :=&\sum\nolimits_{\ell=1}^L\Breve{g}_\ell(u;M)
    \I{1-\varphi_\ell<u\le 1-\varphi_{\ell-1}},\nonumber\\
    \Breve{g}_\ell(u;M)
    :=&\sum\nolimits_{m=0}^M
    g_\ell(b_{\ell,m}^{\rm C})L_{\ell,m}(u)
    +\frac{1-\varphi_{\ell-1}-u}{\varphi_\ell-\varphi_{\ell-1}}
    \big[
    g_\ell(1-\varphi_\ell)
    -\sum\nolimits_{m=0}^M
    g_\ell(b_{\ell,m}^{\rm C})L_{\ell,m}(1-\varphi_\ell)
    \big]\nonumber\\
    &+\frac{u-(1-\varphi_\ell)}{\varphi_\ell-\varphi_{\ell-1}}
    \big[
    g_\ell(1-\varphi_{\ell-1})
    -\sum\nolimits_{m=0}^M
    g_\ell(b_{\ell,m}^{\rm C})L_{\ell,m}(1-\varphi_{\ell-1})
    \big],
    \label{eq_Chebyshev}
\end{align}
where
\(
    L_{\ell,m}(u)
    :=\prod\nolimits_{j=0,j\ne m}^M
    (u-b_{\ell,j}^{\rm C})/(b_{\ell,m}^{\rm C}-b_{\ell,j}^{\rm C})
\)
and
\(
    b_{\ell,m}^{\rm C}
    :=(\varphi_\ell-\varphi_{\ell-1})/2
    \cos((2m+1)\pi/(2M+2))
    +(2-\varphi_\ell-\varphi_{\ell-1})/2,
\)
for \(m=0,\ldots,M\).
The two correction terms ensure
\(
    \Breve{g}_\ell(1-\varphi_\ell;M)
    =g_\ell(1-\varphi_\ell)
\)
and
\(
    \Breve{g}_\ell(1-\varphi_{\ell-1};M)
    =g_\ell(1-\varphi_{\ell-1}).
\)
Figure~\ref{fig:function_approx} illustrates these approximation
methods for the exponential distortion \(g(u)=\{1-\exp(-\lambda u)\}/\{1-\exp(-\lambda)\}\), where \(\lambda>0\). This distortion is non-polynomial and smooth on \([0,1]\). We set \(\lambda=2\) {and \(M=5\) } in the figure.

Let \(\gP\) denote any of the polynomial approximations above. Its associated risk functional is
\begin{equation*}
   \rho_{\gP}(\theta)
    :=\int_0^\infty \gP\big(\bar{F}(y;\theta)\big)dy
    -\int_0^\infty\big[1-\gP\big(\bar{F}(-y;\theta)\big)\big]dy.
\end{equation*}
In Propositions~\ref{prop_Bernstein}–\ref{prop_Bernstein_high} below, we characterize the approximation error and shape-preserving properties of the (piecewise) Bernstein polynomials defined in \eqref{eq_Bernstein}. Similarly, in Proposition~\ref{prop_Chebyshev}, we quantify the approximation error of the (piecewise) Chebyshev interpolation given in \eqref{eq_Chebyshev}. The subscripts \(\tilde g\) and \(\Breve g\) specify the chosen approximation, with its parameters suppressed. Proofs of these results are provided in Section~E of the Appendix.

\begin{proposition}[Bernstein approximation]\label{prop_Bernstein}
     Suppose that Assumptions~\ref{A_moment}--\ref{A_bound} hold and that $g_\ell$ is $(J+1)$-times continuously differentiable with bounded derivatives for $J\in\{1,2\}$. Then we have
    \begin{enumerate}[label=\upshape(\alph*),ref=(\alph*)]
    \item\label{prop_error} (Voronovskaja's Theorem, cf., \citealp{lorentz2012bernstein}) \\
    $\sup_{u\in(1-\varphi_\ell,1-\varphi_{\ell-1}]}
    |\tilde{g}_\ell(u;M,0)-g_\ell(u)|=O(1/M)$, for each $\ell=1,\ldots,L$;
    \item\label{prop_error_rho} $\sup_{\theta\in\Theta}|\rho_{\tilde{g}}(\theta)-\rho_g(\theta)|=O(M^{-J/2})$,
    $\sup_{\theta\in\Theta}\|\nabla_\theta\rho_{\tilde{g}}(\theta)-\nabla_\theta\rho_g(\theta)\|=O(M^{-J/2})$.
    \item\label{prop_mono} $\tilde{g}(\cdot;M,0)$ is non-decreasing on $[0,1]$;
    \item\label{prop_conc} If $g$ is either concave or convex on $[0,1]$, then $\tilde{g}$ preserves this property.
    \end{enumerate}
\end{proposition}

\begin{proposition}[Taylor--Bernstein approximation]\label{prop_Bernstein_high}
     Suppose Assumptions~\ref{A_moment}--\ref{A_bound} hold and that $g_\ell$'s are $(J+1)$-times continuously differentiable with bounded derivatives for $J\ge1$. Then we have
    \begin{enumerate}[label=\upshape(\alph*)]
    \item $\sup_{u\in(1-\varphi_\ell,1-\varphi_{\ell-1}]}
    |\tilde{g}_\ell(u;M,J)-g_\ell(u)|=O(M^{-(J+1)/2})$;
    \item\label{prop_error_rho_high} $\sup_{\theta\in\Theta}|\rho_{\tilde{g}}(\theta)-\rho_g(\theta)|
    =O(M^{-J/2})$,
    $\sup_{\theta\in\Theta}\|\nabla_\theta\rho_{\tilde{g}}(\theta)-\nabla_\theta\rho_g(\theta)\|=O(M^{-J/2})$.
    \end{enumerate}
\end{proposition}

\begin{proposition}[Chebyshev approximation]\label{prop_Chebyshev}
Suppose \(M\ge1\), Assumptions~\ref{A_moment}--\ref{A_bound} hold, and \(g_\ell\) is \(M+2\)-times continuously differentiable with bounded derivatives for each \(\ell=1,\ldots,L\). Then we have
\begin{enumerate}[label=\upshape(\alph*)]

\item\label{error_Chebyshev}
For each $\ell=1,\ldots,L$,
\(
    \sup_{u\in(1-\varphi_\ell,1-\varphi_{\ell-1}]}
    |\Breve{g}_\ell(u;M)-g_\ell(u)|
    =O(\|g_\ell^{(M+1)}\|_\infty/
    [(M+1)!2^{2M}]),
\)
where \(\|\cdot\|_\infty\) denotes the supremum norm on the corresponding subinterval.
\item\label{error_Chebyshev_rho}
\begin{align*}
    {\sup\nolimits_{\theta\in\Theta}}
    |\rho_{\Breve{g}}(\theta)-\rho_g(\theta)|
    =&O\Big(
    \frac{
    {\max\nolimits_{\ell=1,\ldots,L}}
    \{\|g_\ell^{(M+1)}\|_\infty
    +\|g_\ell^{(M+2)}\|_\infty\}}{(M-1)!2^{2M}}\Big),\\
    {\sup\nolimits_{\theta\in\Theta}}
    \|\nabla_\theta\rho_{\Breve{g}}(\theta)
    -\nabla_\theta\rho_g(\theta)\|
    =&O\Big(
    \frac{
    {\max\nolimits_{\ell=1,\ldots,L}}
    \{\|g_\ell^{(M+1)}\|_\infty
    +\|g_\ell^{(M+2)}\|_\infty\}}{(M-1)!2^{2M}}\Big).
\end{align*}
\end{enumerate}
\end{proposition}

By Propositions~\ref{prop_Bernstein}\ref{prop_error_rho}, \ref{prop_Bernstein_high}\ref{prop_error_rho_high}, and \ref{prop_Chebyshev}\ref{error_Chebyshev_rho}, we have
$${\sup\nolimits_{\theta\in\Theta}}|\rho_{\gP}(\theta)-\rho_g(\theta)|=O(\cR)~\mbox{~and~~} {\sup\nolimits_{\theta\in\Theta}}\|\nabla_\theta\rho_{\gP}(\theta)-\nabla_\theta\rho_g(\theta)\|=O(\cR),$$ where $\cR$ denotes the approximation error associated with the chosen technique.
{The Taylor--Bernstein and Chebyshev polynomials need not be nondecreasing or take values in \([0,1]\). Their functionals remain well-defined under the stated conditions but are DRMs only if the polynomials satisfy the distortion requirements. Because \(\gP\) is piecewise polynomial, Algorithm~\ref{alg} can be applied to the surrogate optimization problem \(\min_{\theta\in\Theta}\rho_{\gP}(\theta)\).}
Corollary~\ref{cor_Bernstein} below quantifies the performance of this approach by providing explicit error bounds on the RMSEs of the resulting algorithm iterates, showing how the approximation error $\cR$ propagates to the accuracy of both the DRM and its gradient estimates, and affects the overall optimization process.

\begin{corollary}\label{cor_Bernstein}
Let Assumptions~\ref{A_distr}--\ref{A_rate} hold. Suppose the distortion function \(g\) is approximated using one of the aforementioned techniques with its approximation parameters held fixed, and each \(g_\ell\) satisfies the corresponding smoothness conditions in Propositions~\ref{prop_Bernstein}--\ref{prop_Chebyshev}. Then
\begin{align*}
    \sqrt{\bE[|\rho_k-\rho_g(\theta_k)|^2]}
    =&O\bigg(\frac{\alpha_k}{\gamma_k}\bigg)
    +O(\gamma_k^{1/2})+O(\cR),\\
    \sqrt{\bE[\|G_k-\nabla_\theta\rho_g(\theta)|_{\theta=\theta_k}\|^2]}
    =&O\bigg(\frac{\alpha_k}{\gamma_k}\bigg)
    +O\bigg(\frac{\gamma_k^{1/2}}{c_k}\bigg)
    +O(c_k^2)+O(\cR).
\end{align*}
If, in addition, Assumption~\ref{A_rate_opt} holds, then
\begin{align*}
    \sqrt{\bE[\|\theta_k-\theta^*\|^2]}
    =O\bigg(\frac{\alpha_k}{\gamma_k}\bigg)
    +O\bigg(\frac{\gamma_k^{1/2}}{c_k}\bigg)
    +O(c_k^2)+O(\cR).
\end{align*}
\end{corollary}

\section{Numerical Experiments}\label{sec_num}
We conduct numerical experiments to evaluate the performance of Algorithm~\ref{alg}, using a set of artificially constructed black-box functions, a financial option example, and a portfolio optimization problem. Since VaR and CVaR are special cases of our framework and have been extensively investigated in \cite{hu2025quantile} and \cite{hu2025simulation}, respectively, we focus our experiments on other DRMs to illustrate the broader applicability of the proposed approach.

In all cases, samples are collected sequentially during algorithm execution.
For the synthetic test cases and the financial option example, we compare Algorithm~\ref{alg} with a partition-based method, which combines the discretization scheme \eqref{eq_partition} with Procedures~\eqref{proc_q}--\eqref{proc_dq} using
the partition $\cSg=\{1/(2N),3/(2N),\ldots,(2N-1)/(2N)\}$ with $N=100$ partition points. Specifically, at each iteration of the method, both \eqref{proc_q}--\eqref{proc_dq} are executed
at all $100$ partition points to iteratively approximate the VaRs $q$ and their gradients $\nabla_{\theta}q$ in (\ref{eq_partition}).
For the portfolio optimization experiment in Section~\ref{sec_portfolio}, we compare Algorithm~\ref{alg} with an L-statistic benchmark that uses increasing batches of simulation outputs; implementation details for the benchmark are provided in Section~F of the Appendix. For the implementation of Algorithm~\ref{alg}, step-size and perturbation parameters are chosen according to the theoretical results in Section~\ref{sec_rate}: for DRM optimization in Section~\ref{sec_toy}, we set
$\alpha_k=O(k^{-0.999})$, $\gamma_k=O(k^{-3/4})$, and $c_k=O(k^{-1/8})$;
for DRM estimation and sensitivity analysis at fixed input parameters in Section~\ref{sec_option}, we set $\alpha_k\equiv0$, $\gamma_k=O(k^{-0.999})$, and $c_k=O(k^{-1/6})$.
All computations are performed on a server equipped with 32 CPU cores, 64 GB of RAM, and 500 GB of local storage.

\subsection{Synthetic Test Cases}\label{sec_toy}
We begin with the two artificially constructed test functions:
\begin{enumerate}
\item[]Case 1. $Y(\theta)=[\frac{1}{d}\sum_{i=1}^d(\theta_i-i)^2+1]\xi$, for $\theta_i\in[i-5,i+5]$ and $d=10$,
\item[]Case 2. $Y(\theta)=[\frac{1}{d}\sum_{i=1}^d(\theta_i-\frac{i}{2})^2+1]\xi+\frac{1}{d}\sum_{i=1}^d(\theta_i-i)\theta_i$, for $\theta_i\in[\frac{i}{2}-4,\frac{i}{2}+4]$, $d=20$.
\end{enumerate}
For each case, we consider two unknown noise distributions: $\xi\sim \exp(1)$ (the exponential distribution with the rate parameter $1$) and $\xi\sim t_3$ (the $t$-distribution with 3 degrees of freedom).
We examine \((0.5,0.75)\)-RVaR (Example~\ref{exmp_RVaR}), the mean-Gini (Example~\ref{exmp_power}) with $\lambda_{\operatorname{MG}}=1$, the Gini shortfall with $\varphi=0.75$ and $\lambda_{\operatorname{GS}}=1$ (Example~\ref{exmp_Gini}), and the DRM with the Wang transform $g(u)=\Phi[\Phi^{-1}(u)+\lambda_{\operatorname{WT}}]$ where $\Phi$ denotes the standard normal CDF and $\lambda_{\operatorname{WT}}=0.5$ \cite[cf.,][]{wang2000class}, yielding \(16\) combinations of test functions, noise distributions, and risk measures. The Wang transform provides a numerical assessment of the polynomial approximation procedure in Section~\ref{sec_ext}.
Tables~\ref{tab:transformation_rules}--\ref{tab:base_risks} list the DRMs as functions of decision variables $\theta$ across all scenarios for comparison.

In the implementation of Algorithm~\ref{alg}, we set step- and perturbation-sizes as $\alpha_k=10/k^{0.999}$, $\gamma_k=100/k^{3/4}$, and $c_k=1/(k+10)^{1/8}$. The initial values are
chosen as $\rho_0=q_0(\cdot)=0$, $G_0=D_{q,0}(\cdot)=(0,\ldots,0)^\intercal\in\Re^d$, $\theta_0=(-2,-1,\ldots,7)^\intercal\in\Re^{10}$ in Case 1 and $\theta_0=(-2.5,-2,\ldots,7)^\intercal\in\Re^{20}$ in Case 2.
For the Wang transform, we employ both the standard Bernstein approximation and Chebyshev interpolation with polynomial degree 3.
A projection operator is applied in Procedures~\eqref{proc_dq} and \eqref{proc_drho} to ensure that the gradient estimators remain within the bounded set $\{x\in\Re^d:\|x\|\le50\}$ \cite[cf., e.g.,][]{hu2022stochastic, hu2025quantile}.
The experiments are conducted with the number of iterations ranging from $10^2$ to $10^5$, and the performance is evaluated based on $40$ independent replications.

\begin{table}
\TABLE
{Distortion Risk Measures considered in the test scenarios.\label{tab:transformation_rules}}
{\begin{tabular*}{\textwidth}{@{\extracolsep{\fill}}lcc}
\toprule
Risk measure & Case 1 & Case 2 \\
\midrule
RVaR & \(c_1(\theta)\rho_{\operatorname{RVaR}}[\xi]\) & \(c_2(\theta)+c_3(\theta)\rho_{\operatorname{RVaR}}[\xi]\) \\
Mean-Gini & \(c_1(\theta)\rho_{\operatorname{MG}}[\xi]\) & \(c_2(\theta)+c_3(\theta)\rho_{\operatorname{MG}}[\xi]\) \\
Gini Shortfall & \(c_1(\theta)\rho_{\operatorname{GS}}[\xi]\) & \(c_2(\theta)+c_3(\theta)\rho_{\operatorname{GS}}[\xi]\) \\
Wang Transform & \(c_1(\theta)\rho_{\operatorname{WT}}[\xi]\) & \(c_2(\theta)+c_3(\theta)\rho_{\operatorname{WT}}[\xi]\) \\
\bottomrule
\end{tabular*}}
{\textit{Note.} $c_1(\theta):=(1/d)\sum_{i=1}^d(\theta_i-i)^2+1$ and $c_3(\theta):=(1/d)\sum_{i=1}^d(\theta_i-i/2)^2+1$ are the scale factors.
$c_2(\theta):=(1/d)\sum_{i=1}^d(\theta_i-i)\theta_i$ represents the shift factor.}
\end{table}

\begin{table}
\TABLE
{Analytical Expressions of Base Distortion Risk Measures $\rho_g[\xi]$.\label{tab:base_risks}}
{\begin{tabular}{lcc}
\toprule
Risk measure & \(\exp(1)\) & \(t_3\) \\
\midrule
\((0.5,0.75)\)-RVaR & \(1\) & \(4\int_{0.5}^{0.75}F_{t_3}^{-1}(u)\,du\) \\
Mean-Gini & \(1+\lambda/2\)
    & \(\int_0^1 F_{t_3}^{-1}(u)(1-\lambda+2\lambda u)\,du\) \\
Gini Shortfall & \(1+\lambda/2-\ln(1-\varphi)\)
    & \(\int_\varphi^1 F_{t_3}^{-1}(u)[1/(1-\varphi)+\lambda(2u-1-\varphi)/(1-\varphi)^2]\,du\) \\
Wang Transform & \(\int_{-\infty}^{\infty}[-\ln(1-\Phi(y))]\phi(y-\lambda)\,dy\)
    & \(\int_{-\infty}^{\infty}F_{t_3}^{-1}(\Phi(y))\phi(y-\lambda)\,dy\) \\
\bottomrule
\end{tabular}}
{\textit{Note.} $\Phi(\cdot)$ and $\phi(\cdot)$ are the CDF and the density function of the standard normal distribution. $F_{t_3}^{-1}(\cdot)$ is the quantile function of the $t_3$ distribution.
}
\end{table}

Figure~\ref{fig:conv_Toy} illustrates the convergence behavior of Algorithm~\ref{alg} and the partition-based method across different scenarios by plotting the average objective values at the estimated solutions, i.e., $(1/40)\sum_{r=1}^{40}\rho_g(\theta_k(\omega_r))$, as functions of the iteration number $k$.
Table~\ref{tab:toy} reports the average running time (in seconds) and descriptive statistics of the final solutions $\theta_K$ obtained by the comparison algorithms, based on $40$ independent replications.
For RVaR, the two methods achieve comparable final accuracy under exponential noise, whereas Algorithm~\ref{alg} yields smaller RMSEs of both \(\rho_g(\theta_K)\) and \(\theta_K\) under \(t_3\) noise. Across the four RVaR scenarios, Algorithm~\ref{alg} requires 57--68 seconds, compared with 1,206--1,272 seconds for the partition-based method. For the mean-Gini and Gini shortfall, Algorithm~\ref{alg} achieves smaller RMSEs in all scenarios. Algorithm~\ref{alg} is also more computationally efficient throughout Table~\ref{tab:toy}, as it avoids evaluating VaRs and their gradients over a dense set of partition points. For the Wang transform, both Bernstein and Chebyshev approximations yield smaller RMSEs than the partition-based method.
Figure~\ref{fig:rate_Toy} compares the empirical convergence behavior through log-log plots of the decision RMSEs against the iteration number \(k\). For RVaR, the fitted slopes for Algorithm~\ref{alg} range from \(-0.66\) to \(-0.44\) over the tested iteration range.

\newsavebox{\boxConvToy}
\sbox{\boxConvToy}{%
    \begin{minipage}{\textwidth}
        \centering

        \begin{minipage}{0.23\textwidth}
            \centering
            \includegraphics[width=\linewidth]{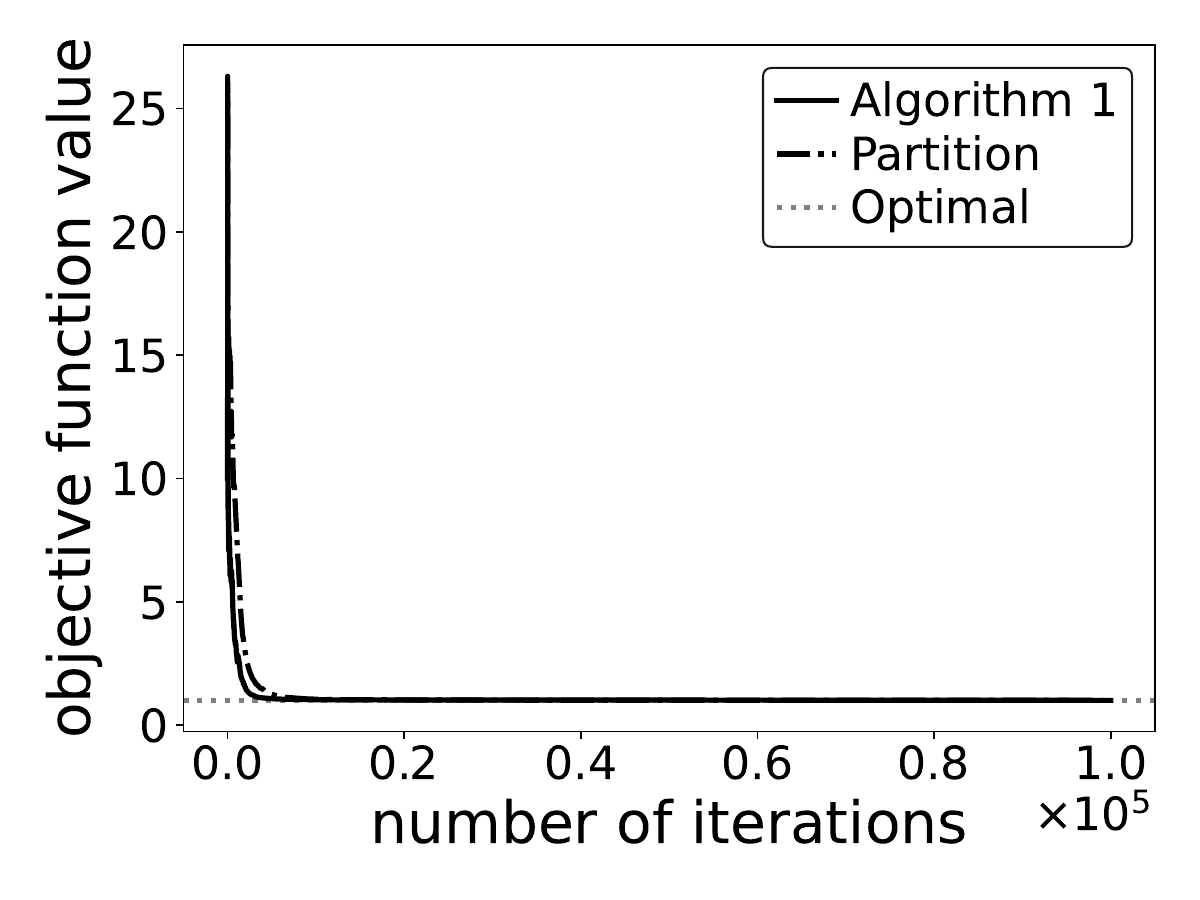}\par
            \vspace{3pt}\footnotesize (a) Case 1, \(\exp\), RVaR
        \end{minipage}\hfill
        \begin{minipage}{0.23\textwidth}
            \centering
            \includegraphics[width=\linewidth]{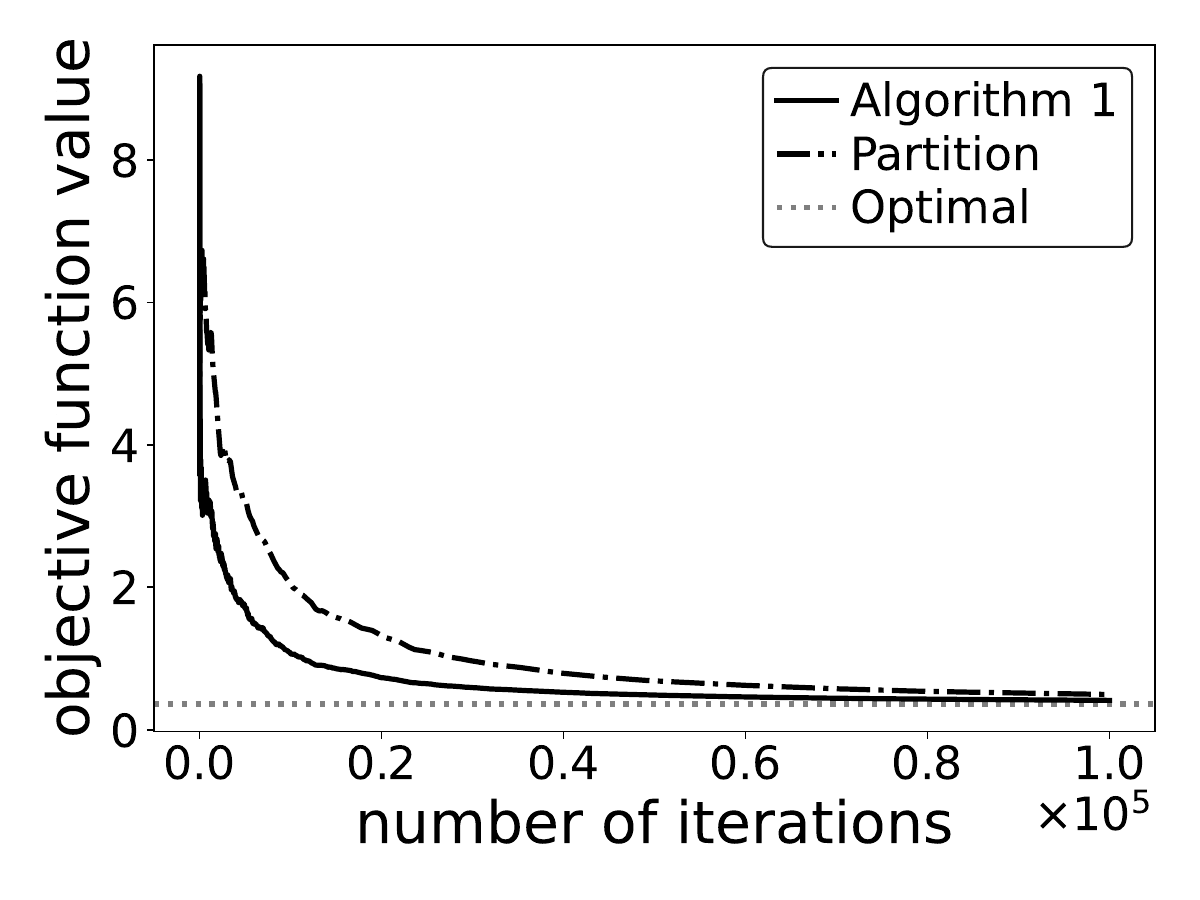}\par
            \vspace{3pt}\footnotesize (b) Case 1, \(t\), RVaR
        \end{minipage}\hfill
        \begin{minipage}{0.23\textwidth}
            \centering
            \includegraphics[width=\linewidth]{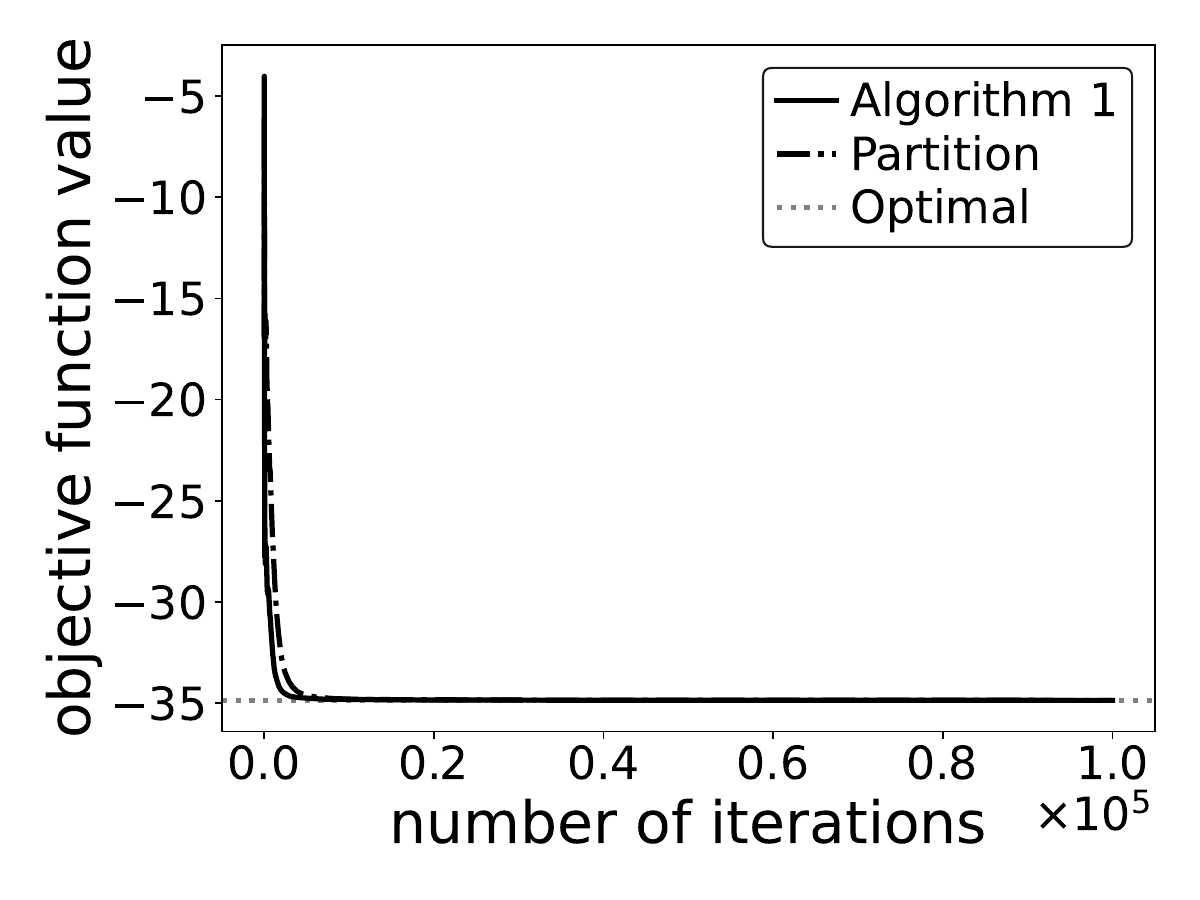}\par
            \vspace{3pt}\footnotesize (c) Case 2, \(\exp\), RVaR
        \end{minipage}\hfill
        \begin{minipage}{0.23\textwidth}
            \centering
            \includegraphics[width=\linewidth]{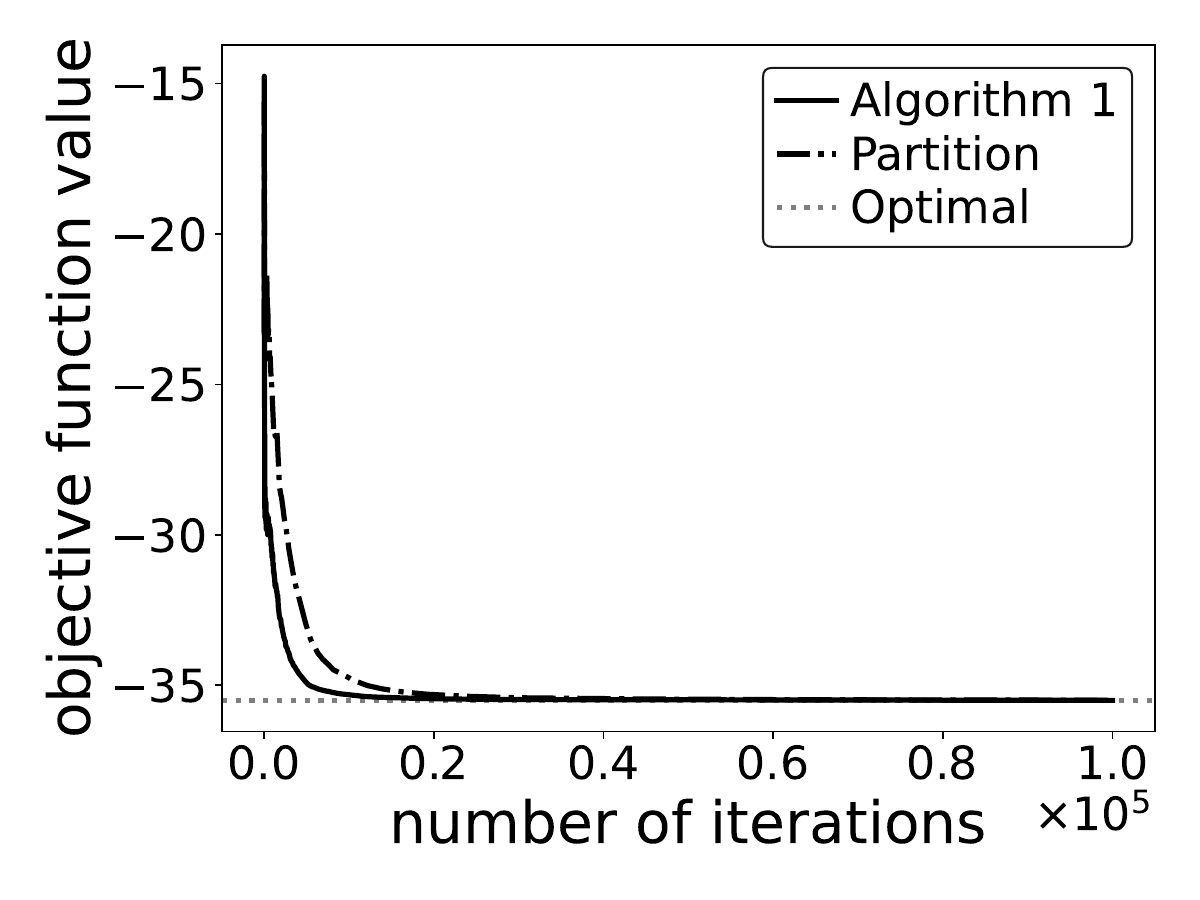}\par
            \vspace{3pt}\footnotesize (d) Case 2, \(t\), RVaR
        \end{minipage}

        \par\vspace{3ex}

        \begin{minipage}{0.23\textwidth}
            \centering
            \includegraphics[width=\linewidth]{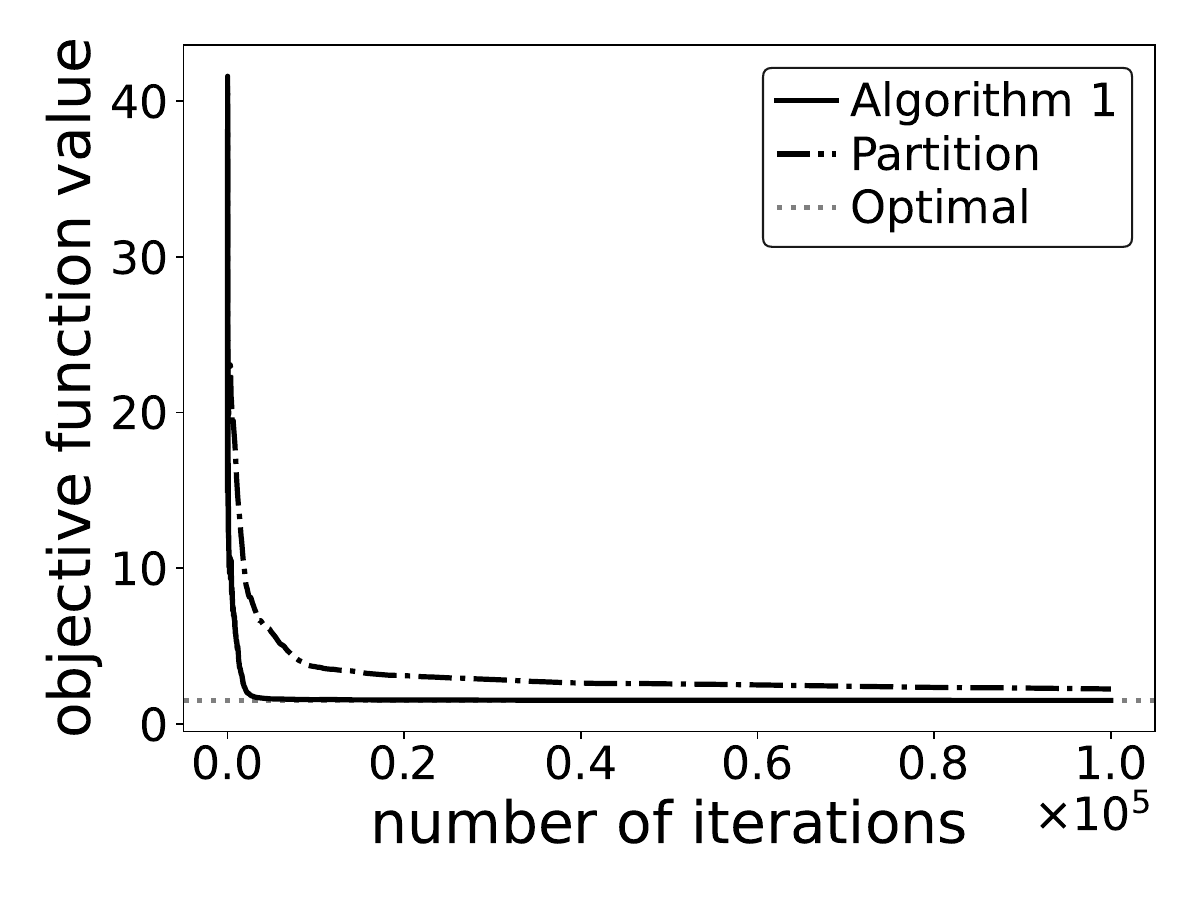}\par
            \vspace{3pt}\footnotesize (e) Case 1, $\exp$, M-G
        \end{minipage}\hfill
        \begin{minipage}{0.23\textwidth}
            \centering
            \includegraphics[width=\linewidth]{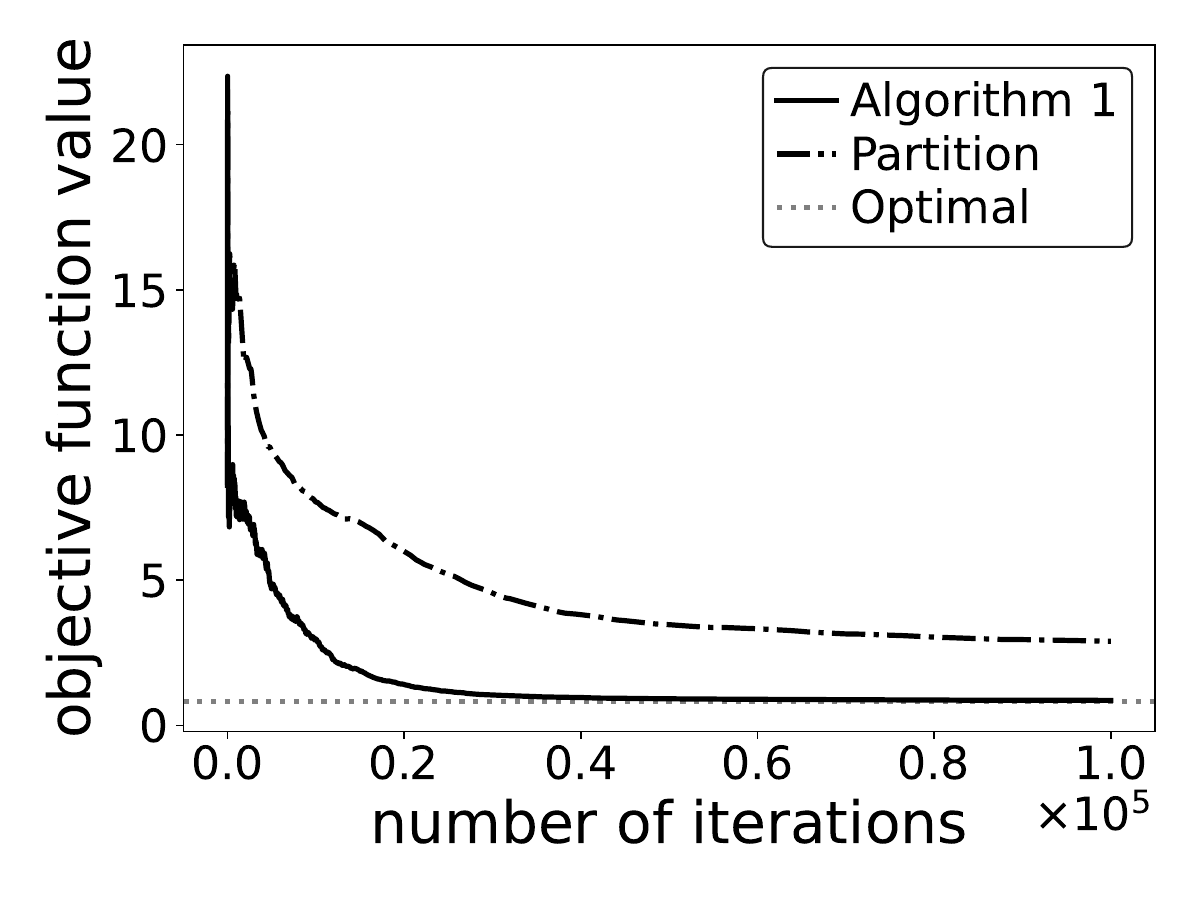}\par
            \vspace{3pt}\footnotesize (f) Case 1, $t$, M-G
        \end{minipage}\hfill
        \begin{minipage}{0.23\textwidth}
            \centering
            \includegraphics[width=\linewidth]{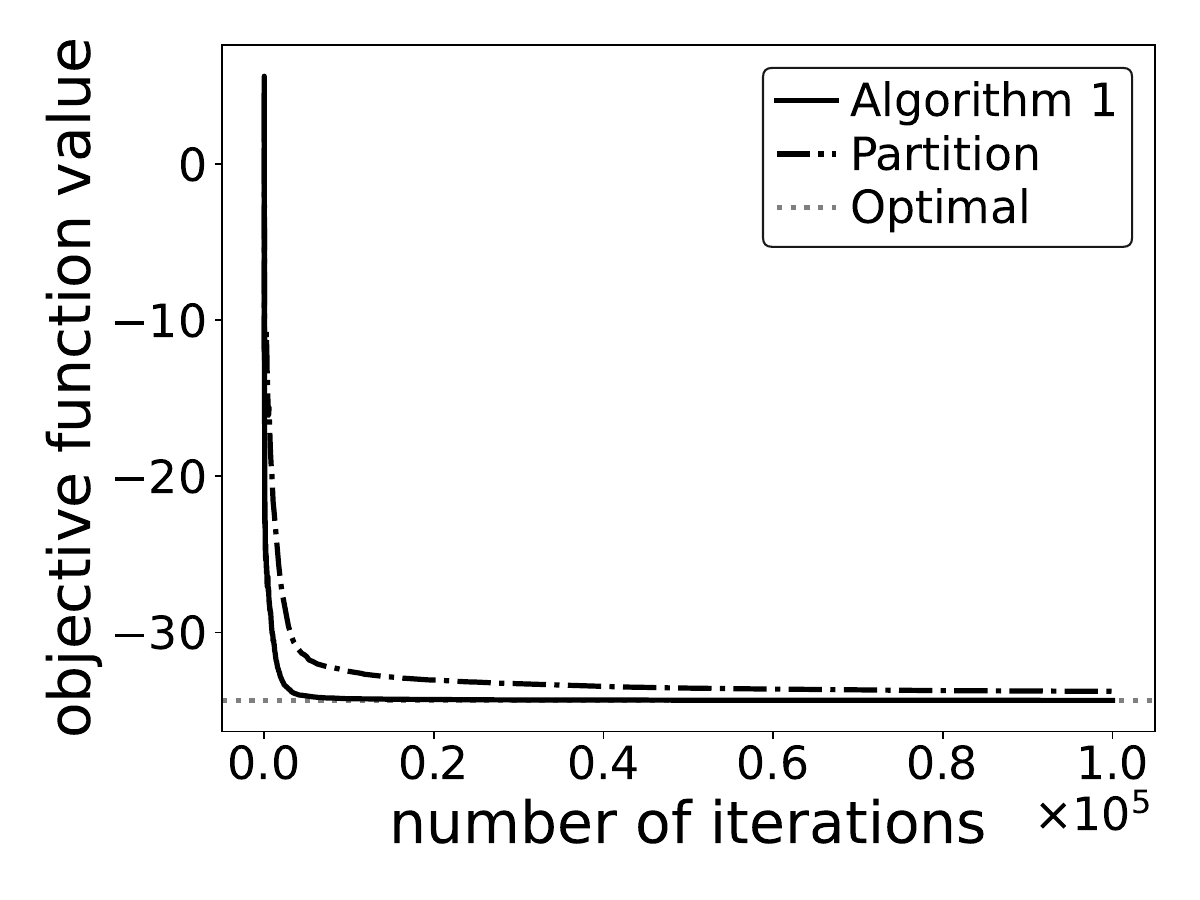}\par
            \vspace{3pt}\footnotesize (g) Case 2, $\exp$, M-G
        \end{minipage}\hfill
        \begin{minipage}{0.23\textwidth}
            \centering
            \includegraphics[width=\linewidth]{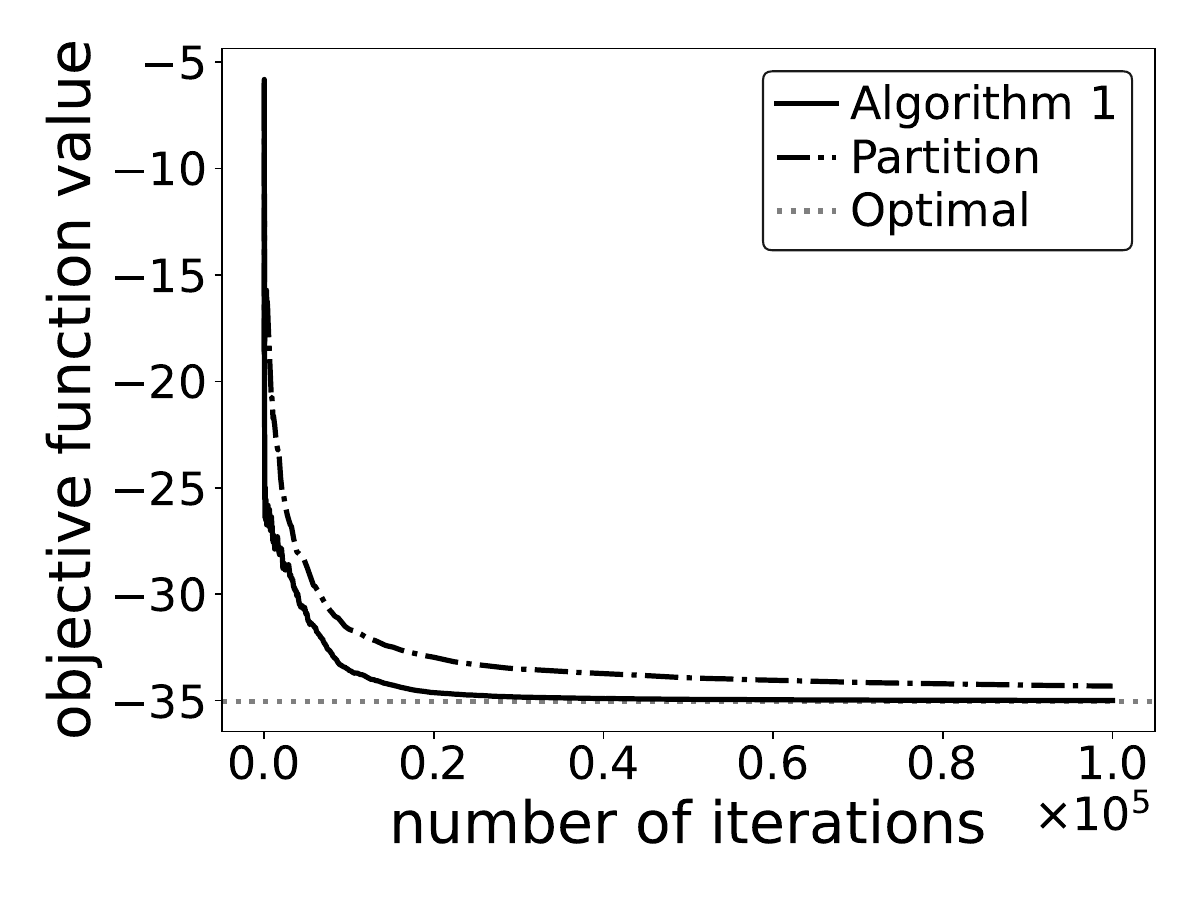}\par
            \vspace{3pt}\footnotesize (h) Case 2, $t$, M-G
        \end{minipage}

        \par\vspace{3ex}

        \begin{minipage}{0.23\textwidth}
            \centering
            \includegraphics[width=\linewidth]{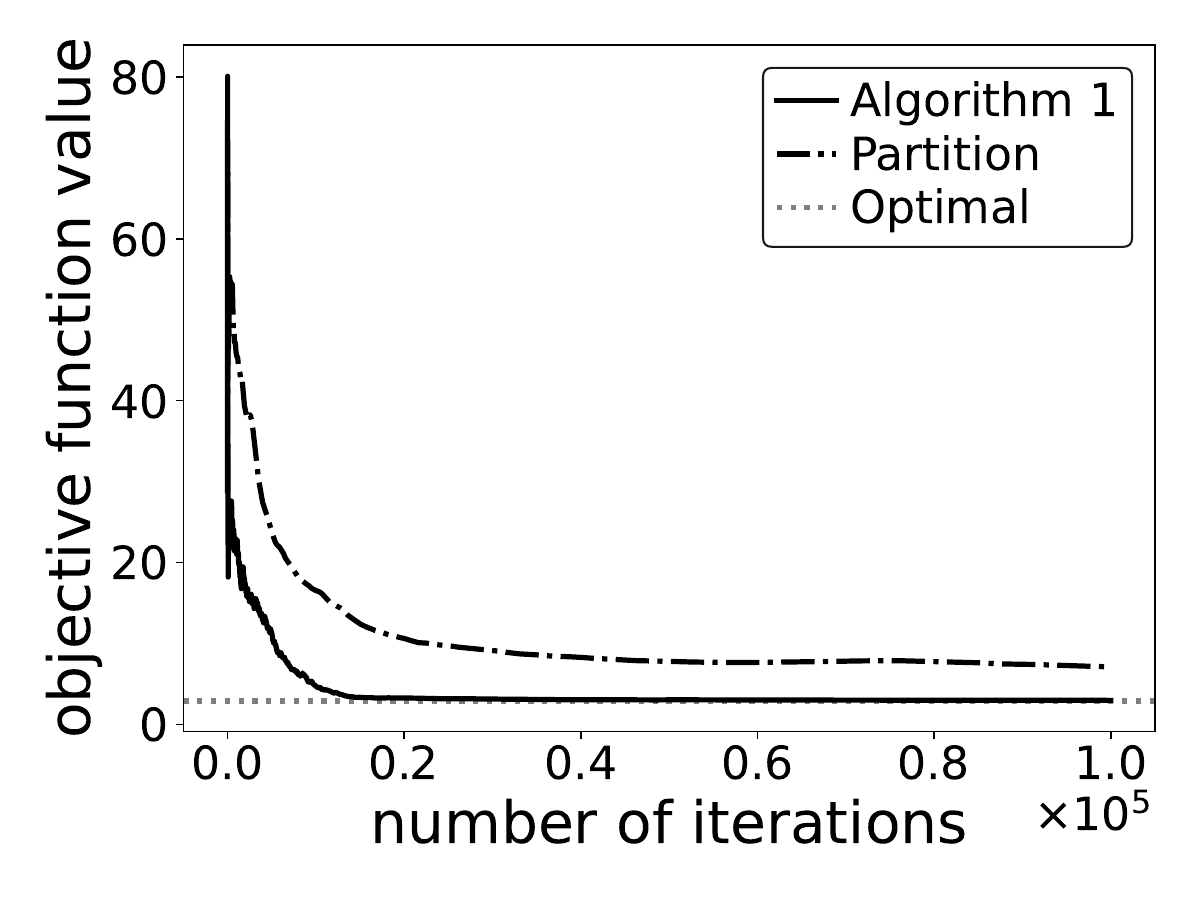}\par
            \vspace{3pt}\footnotesize (i) Case 1, $\exp$, G-S
        \end{minipage}\hfill
        \begin{minipage}{0.23\textwidth}
            \centering
            \includegraphics[width=\linewidth]{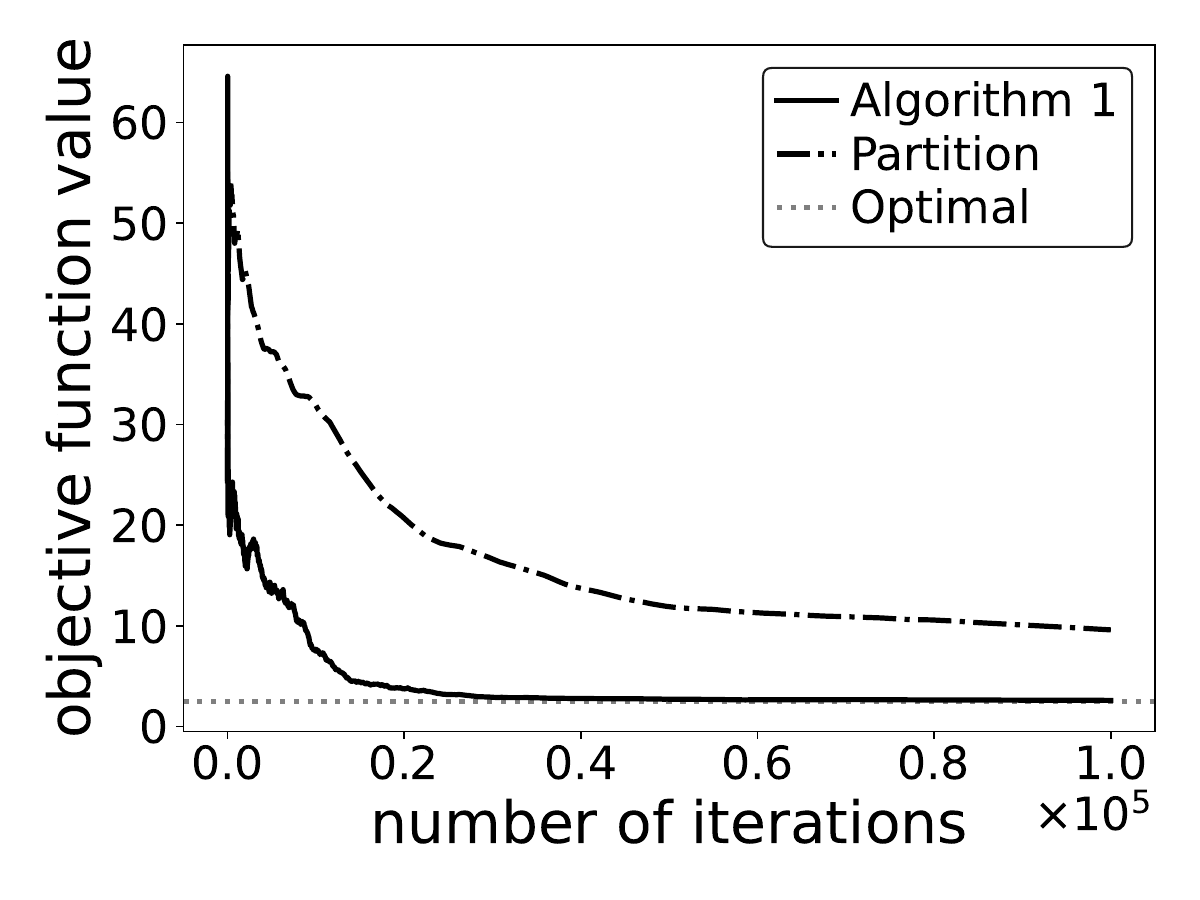}\par
            \vspace{3pt}\footnotesize (j) Case 1, $t$, G-S
        \end{minipage}\hfill
        \begin{minipage}{0.23\textwidth}
            \centering
            \includegraphics[width=\linewidth]{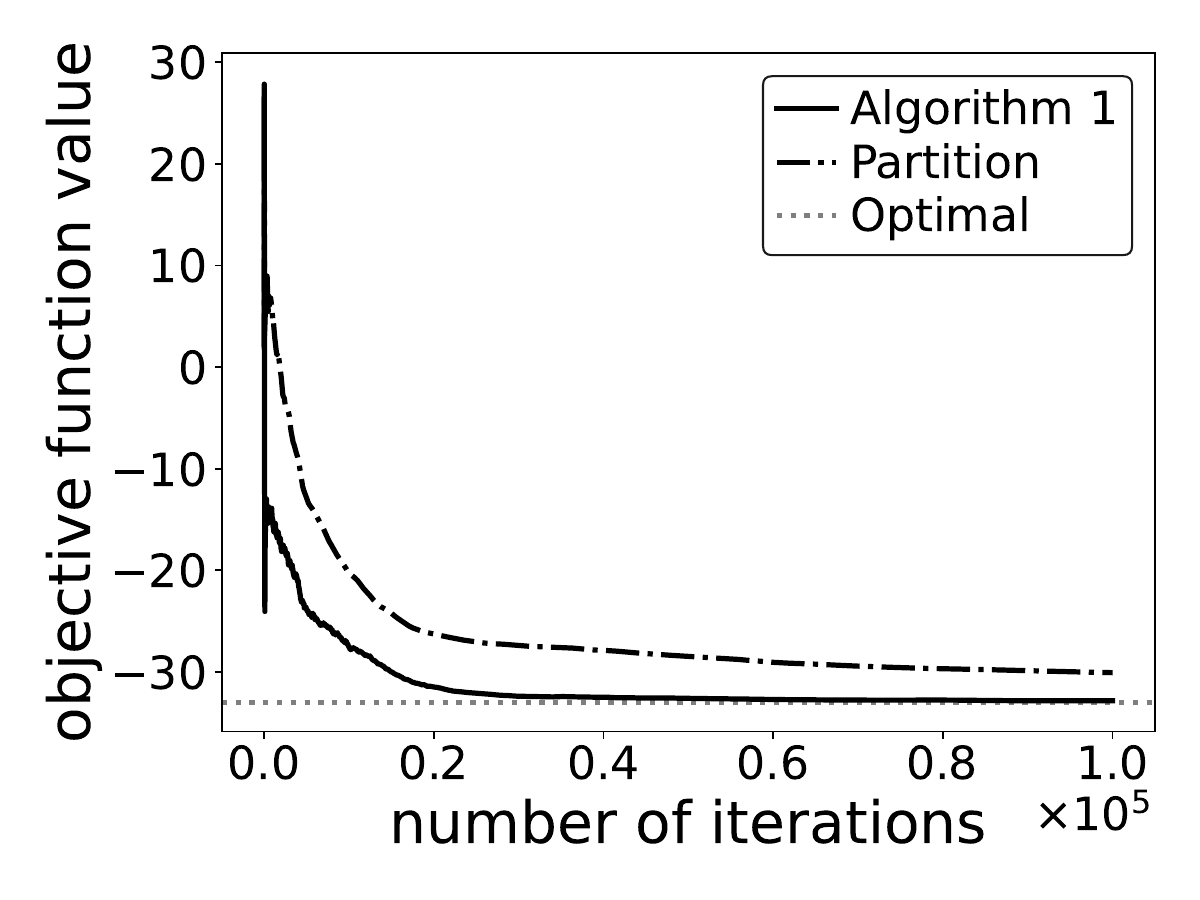}\par
            \vspace{3pt}\footnotesize (k) Case 2, $\exp$, G-S
        \end{minipage}\hfill
        \begin{minipage}{0.23\textwidth}
            \centering
            \includegraphics[width=\linewidth]{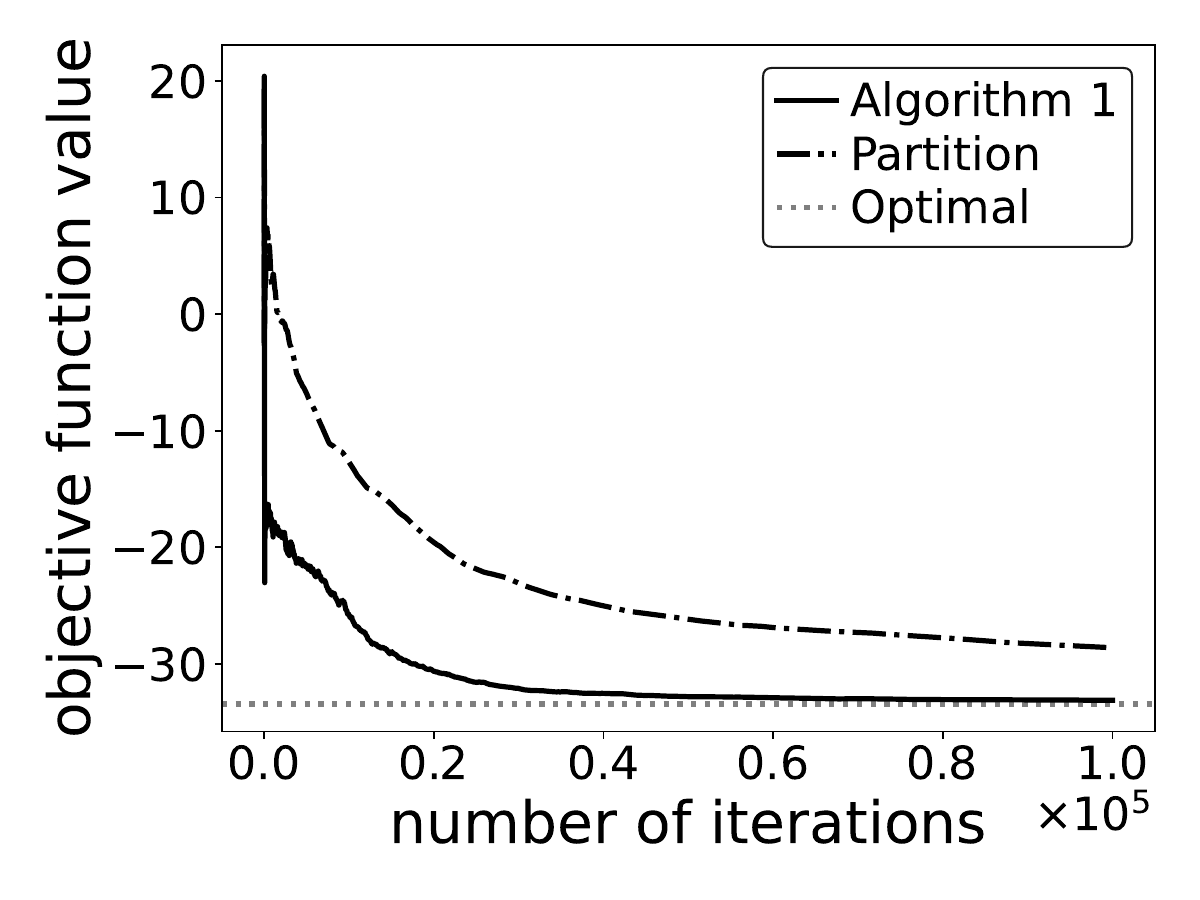}\par
            \vspace{3pt}\footnotesize (l) Case 2, $t$, G-S
        \end{minipage}

        \par\vspace{3ex}

        \begin{minipage}{0.23\textwidth}
            \centering
            \includegraphics[width=\linewidth]{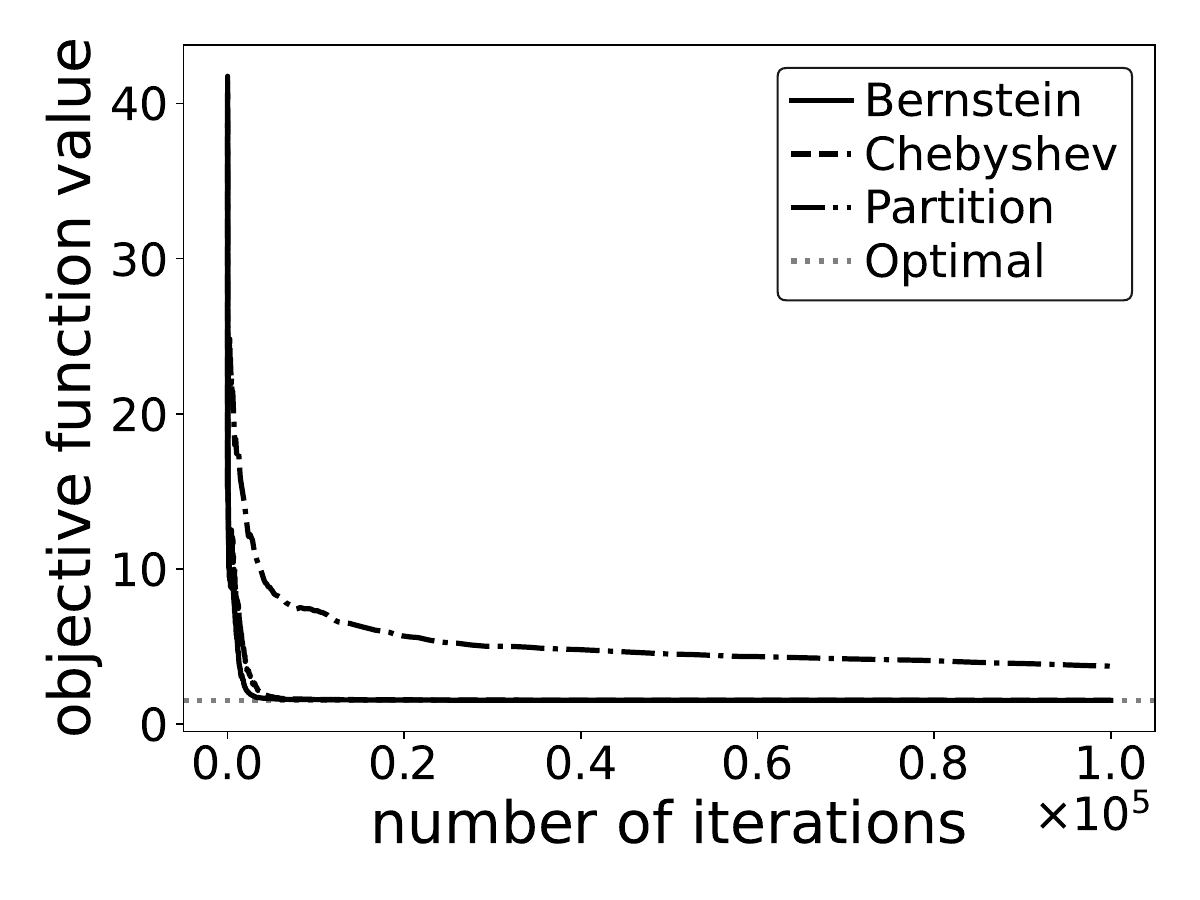}\par
            \vspace{3pt}\footnotesize (m) Case 1, $\exp$, Wang
        \end{minipage}\hfill
        \begin{minipage}{0.23\textwidth}
            \centering
            \includegraphics[width=\linewidth]{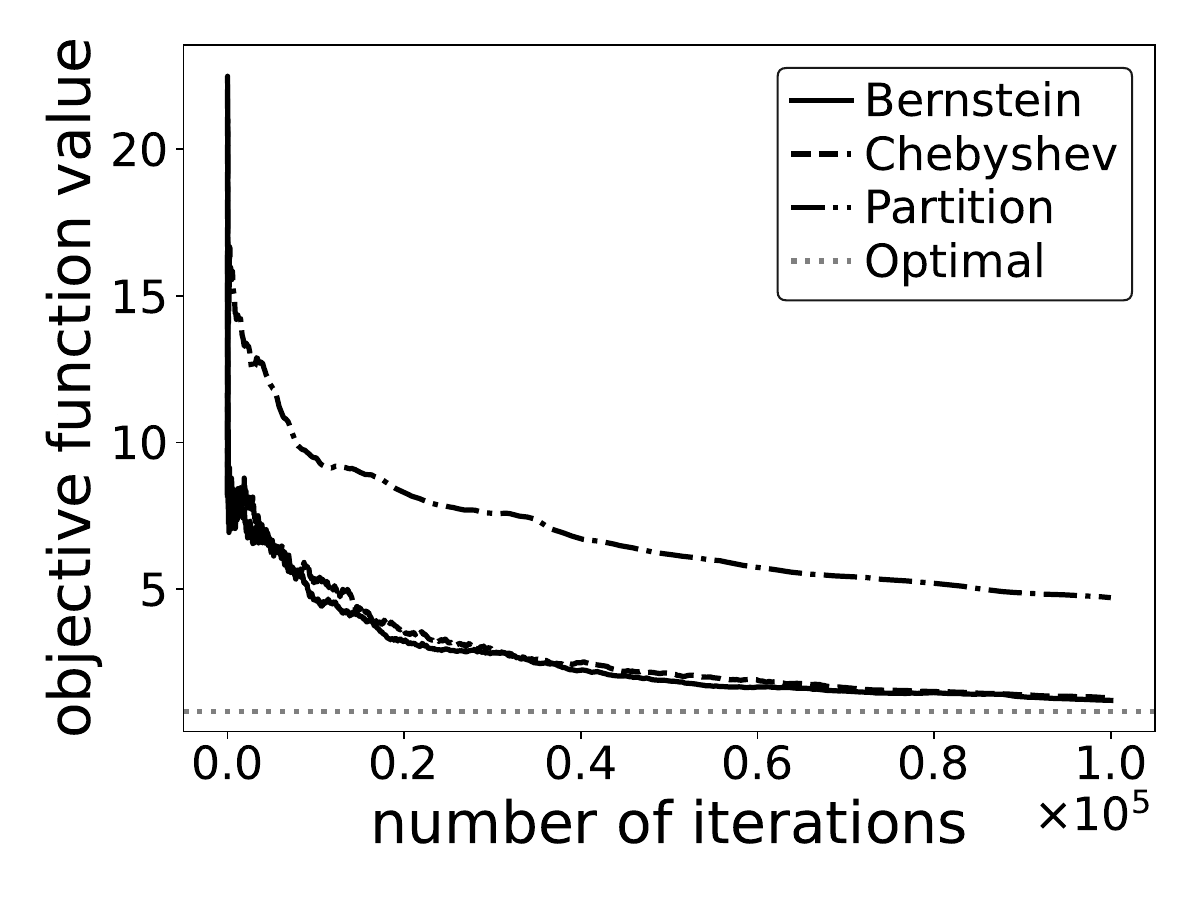}\par
            \vspace{3pt}\footnotesize (n) Case 1, $t$, Wang
        \end{minipage}\hfill
        \begin{minipage}{0.23\textwidth}
            \centering
            \includegraphics[width=\linewidth]{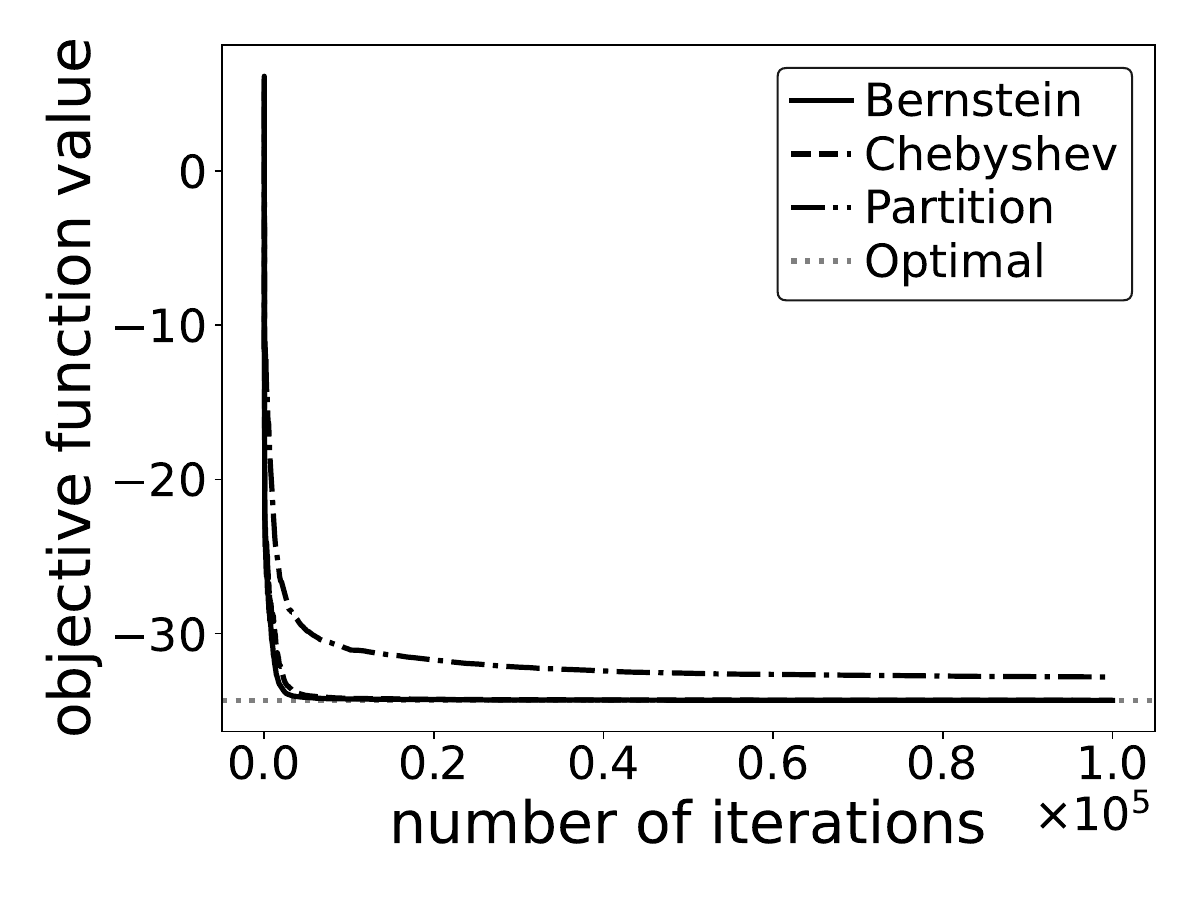}\par
            \vspace{3pt}\footnotesize (o) Case 2, $\exp$, Wang
        \end{minipage}\hfill
        \begin{minipage}{0.23\textwidth}
            \centering
            \includegraphics[width=\linewidth]{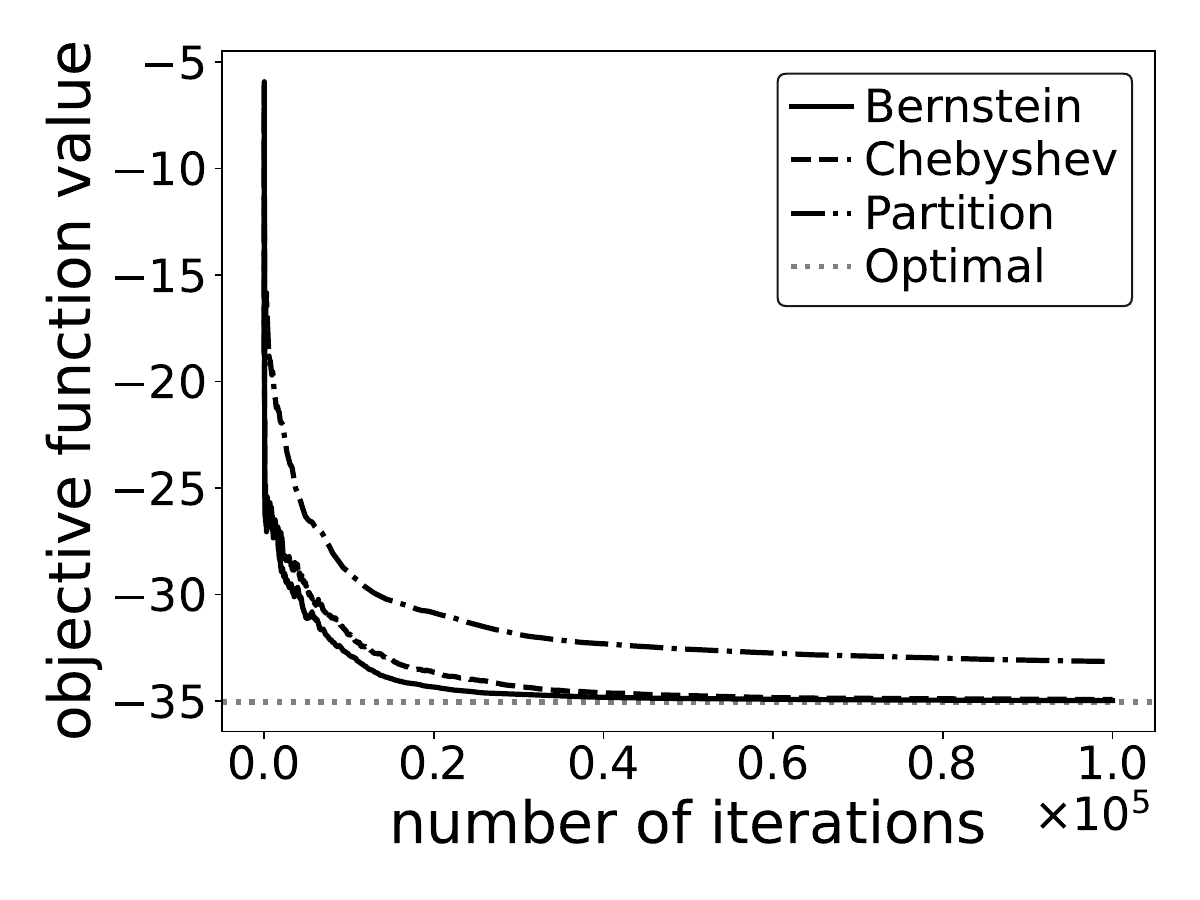}\par
            \vspace{3pt}\footnotesize (p) Case 2, $t$, Wang
        \end{minipage}
    \end{minipage}%
}

\begin{figure}
    \FIGURE
    {\usebox{\boxConvToy}}
    {Convergence Behavior of Algorithm 1 and the Partition-Based Method on Synthetic Test Cases.\label{fig:conv_Toy}}
    {Rows from top to bottom represent RVaR, mean-Gini (M-G), Gini shortfall (G-S), and the Wang transform (Wang), respectively.}
\end{figure}

\begin{table}
\TABLE
{Comparison Between Algorithm~\ref{alg} and the Partition-Based Method on Toy Examples.\label{tab:toy}}
{\setlength{\tabcolsep}{3.5pt}
\begin{tabular}{l *{9}{c}}
\toprule
& \multicolumn{2}{c}{RVaR} & \multicolumn{2}{c}{M-G} & \multicolumn{2}{c}{G-S} & \multicolumn{3}{c}{Wang} \\
\cmidrule(lr){2-3} \cmidrule(lr){4-5} \cmidrule(lr){6-7} \cmidrule(lr){8-10}
& A1 & P & A1 & P & A1 & P & A1-B & A1-C & P \\
\midrule
\textit{Case 1, Exponential} & & & & & & & & & \\
ART (s) & $\mathbf{57}$ & $1206$ & $\mathbf{27}$ & $477$ & $\mathbf{38}$ & $723$ & $\mathbf{40}$ & $\mathbf{40}$ & $852$ \\
Avg. $\rho_g(\theta_K)$ & $\mathbf{1.005}$ & $\mathbf{1.005}$ & $\mathbf{1.509}$ & $2.240$ & $\mathbf{2.964}$ & $7.097$ & $\mathbf{1.539}$ & $1.542$ & $3.742$ \\
RMSE of $\rho_g(\theta_K)$ & $\mathbf{0.005}$ & $\mathbf{0.005}$ & $\mathbf{0.009}$ & $1.225$ & $\mathbf{0.085}$ & $5.003$ & $\mathbf{0.010}$ & $0.013$ & $2.814$ \\
RMSE of $\theta_K$ & $0.224$ & $\mathbf{0.220}$ & $\mathbf{0.245}$ & $2.221$ & $\mathbf{0.518}$ & $3.820$ & $\mathbf{0.237}$ & $0.276$ & $3.803$ \\
\midrule
\textit{Case 1, $t$-distribution} & & & & & & & & & \\
ART (s) & $\mathbf{58}$ & $1210$ & $\mathbf{27}$ & $478$ & $\mathbf{37}$ & $716$ & $42$ & $\mathbf{41}$ & $837$ \\
Avg. $\rho_g(\theta_K)$ & $\mathbf{0.412}$ & $0.494$ & $\mathbf{0.861}$ & $2.894$ & $\mathbf{2.591}$ & $9.592$ & $\mathbf{1.205}$ & $1.289$ & $4.707$ \\
RMSE of $\rho_g(\theta_K)$ & $\mathbf{0.064}$ & $0.179$ & $\mathbf{0.038}$ & $2.655$ & $\mathbf{0.167}$ & $8.604$ & $\mathbf{0.653}$ & $0.869$ & $4.468$ \\
RMSE of $\theta_K$ & $\mathbf{1.201}$ & $1.933$ & $\mathbf{0.644}$ & $5.000$ & $\mathbf{0.787}$ & $5.415$ & $\mathbf{2.165}$ & $2.389$ & $6.881$ \\
\midrule
\textit{Case 2, Exponential} & & & & & & & & & \\
ART (s) & $\mathbf{64}$ & $1214$ & $\mathbf{29}$ & $482$ & $\mathbf{38}$ & $719$ & $\mathbf{44}$ & $\mathbf{44}$ & $846$ \\
Avg. $\rho_g(\theta_K)$ & $\mathbf{-34.866}$ & $\mathbf{-34.866}$ & $\mathbf{-34.356}$ & $-33.778$ & $\mathbf{-32.826}$ & $-30.067$ & $\mathbf{-34.329}$ & $-34.325$ & $-32.813$ \\
RMSE of $\rho_g(\theta_K)$ & $\mathbf{0.009}$ & $0.010$ & $\mathbf{0.020}$ & $1.134$ & $\mathbf{0.174}$ & $3.433$ & $\mathbf{0.017}$ & $0.021$ & $1.810$ \\
RMSE of $\theta_K$ & $\mathbf{0.300}$ & $\mathbf{0.300}$ & $\mathbf{0.393}$ & $2.186$ & $\mathbf{0.916}$ & $3.878$ & $\mathbf{0.359}$ & $0.401$ & $3.480$ \\
\midrule
\textit{Case 2, $t$-distribution} & & & & & & & & & \\
ART (s) & $\mathbf{68}$ & $1272$ & $\mathbf{30}$ & $481$ & $\mathbf{40}$ & $722$ & $46$ & $\mathbf{45}$ & $846$ \\
Avg. $\rho_g(\theta_K)$ & $\mathbf{-35.501}$ & $-35.495$ & $\mathbf{-34.994}$ & $-34.317$ & $\mathbf{-33.119}$ & $-28.612$ & $\mathbf{-34.986}$ & $-34.935$ & $-33.151$ \\
RMSE of $\rho_g(\theta_K)$ & $\mathbf{0.015}$ & $0.021$ & $\mathbf{0.059}$ & $0.825$ & $\mathbf{0.336}$ & $5.634$ & $\mathbf{0.071}$ & $0.131$ & $2.221$ \\
RMSE of $\theta_K$ & $\mathbf{0.458}$ & $0.542$ & $\mathbf{0.770}$ & $2.829$ & $\mathbf{1.357}$ & $5.296$ & $\mathbf{0.864}$ & $1.145$ & $4.573$ \\
\bottomrule
\end{tabular}}
{\textit{Note.} Computational results based on 40 independent replications. ART denotes the average running time. A1 and P denote Algorithm~\ref{alg} and the partition-based method, respectively. A1-B and A1-C denote Algorithm~\ref{alg} combined with the Bernstein and Chebyshev approximations, respectively. Bold indicates the best displayed value within each measure; ties are all bolded.}
\end{table}

\newsavebox{\boxRateToy}
\sbox{\boxRateToy}{%
    \begin{minipage}{\textwidth}
        \centering

        \begin{minipage}{0.23\textwidth}
            \centering
            \includegraphics[width=\linewidth]{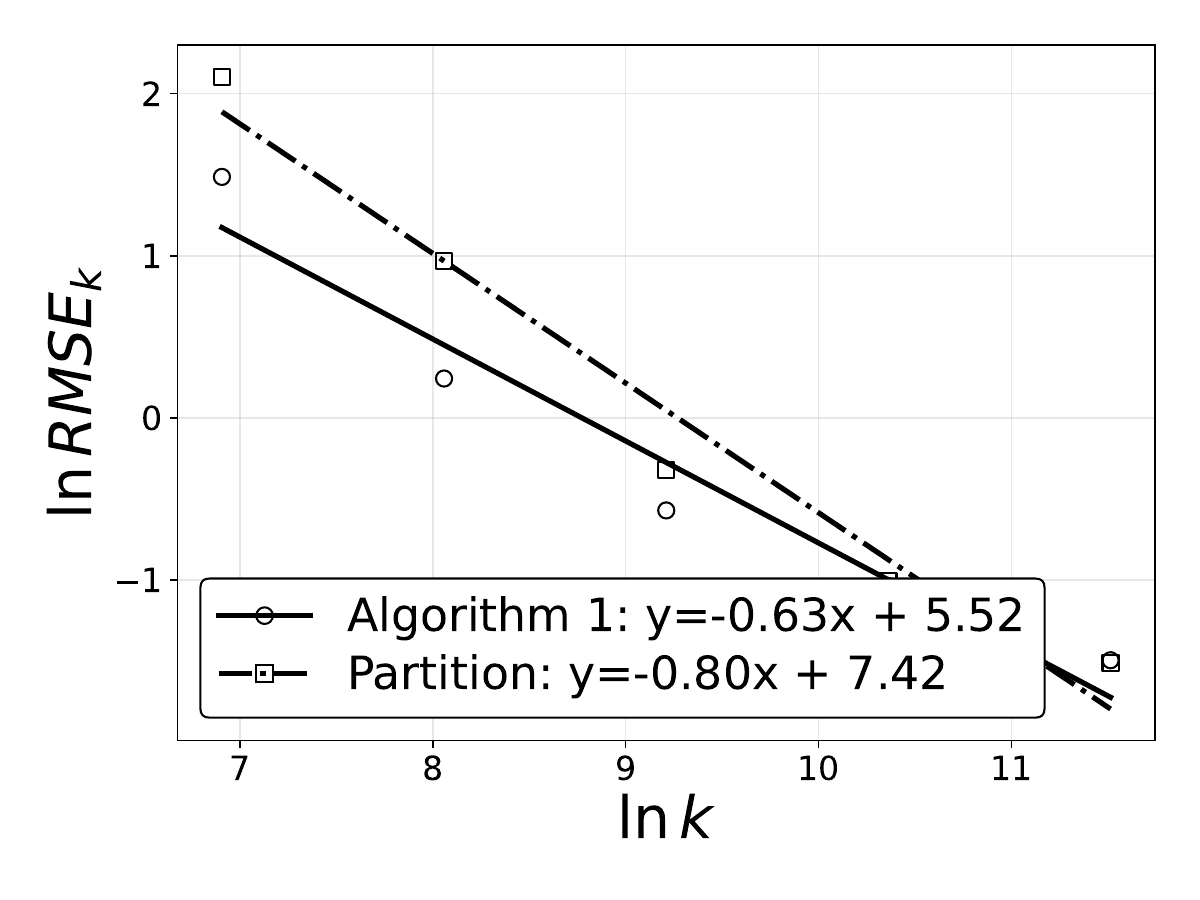}\par
            \vspace{3pt}\footnotesize (a) Case 1, \(\exp\), RVaR
        \end{minipage}\hfill
        \begin{minipage}{0.23\textwidth}
            \centering
            \includegraphics[width=\linewidth]{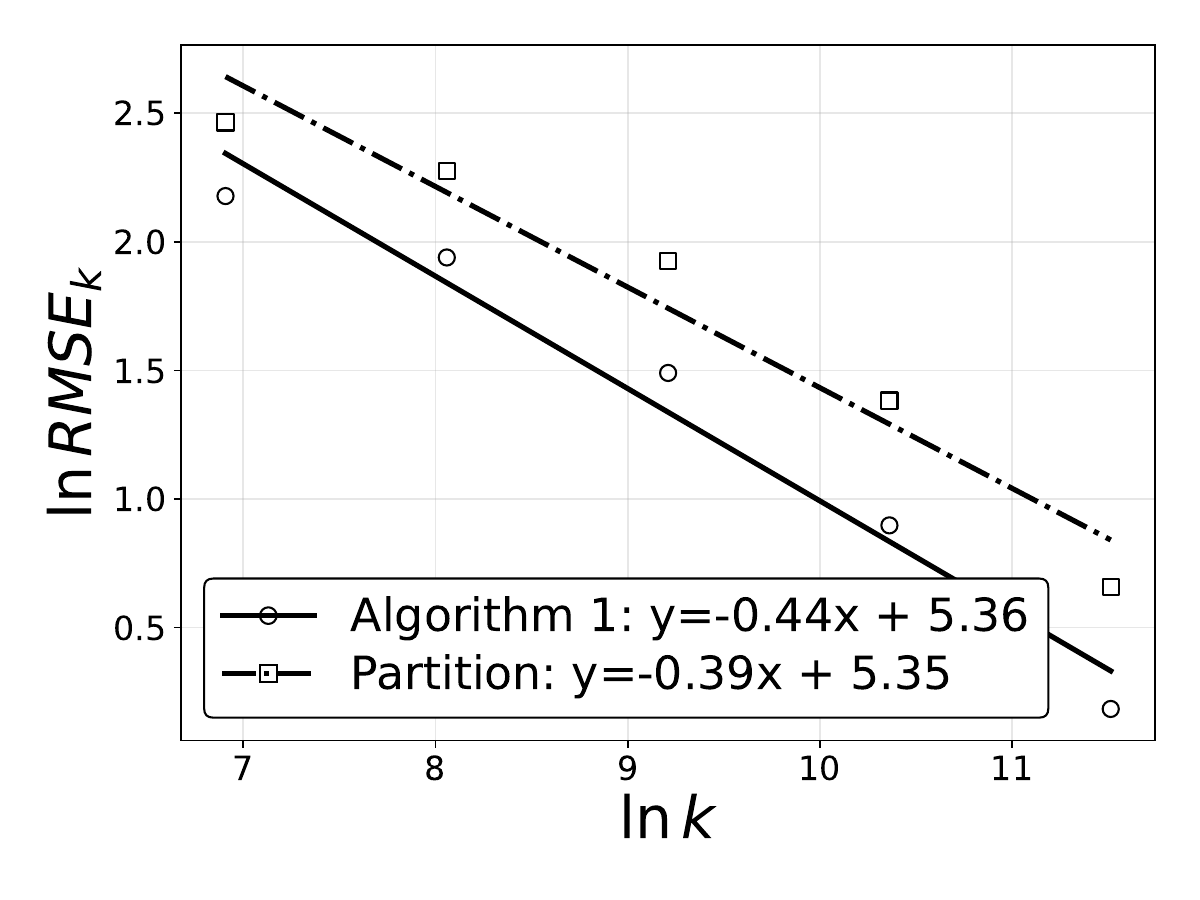}\par
            \vspace{3pt}\footnotesize (b) Case 1, \(t\), RVaR
        \end{minipage}\hfill
        \begin{minipage}{0.23\textwidth}
            \centering
            \includegraphics[width=\linewidth]{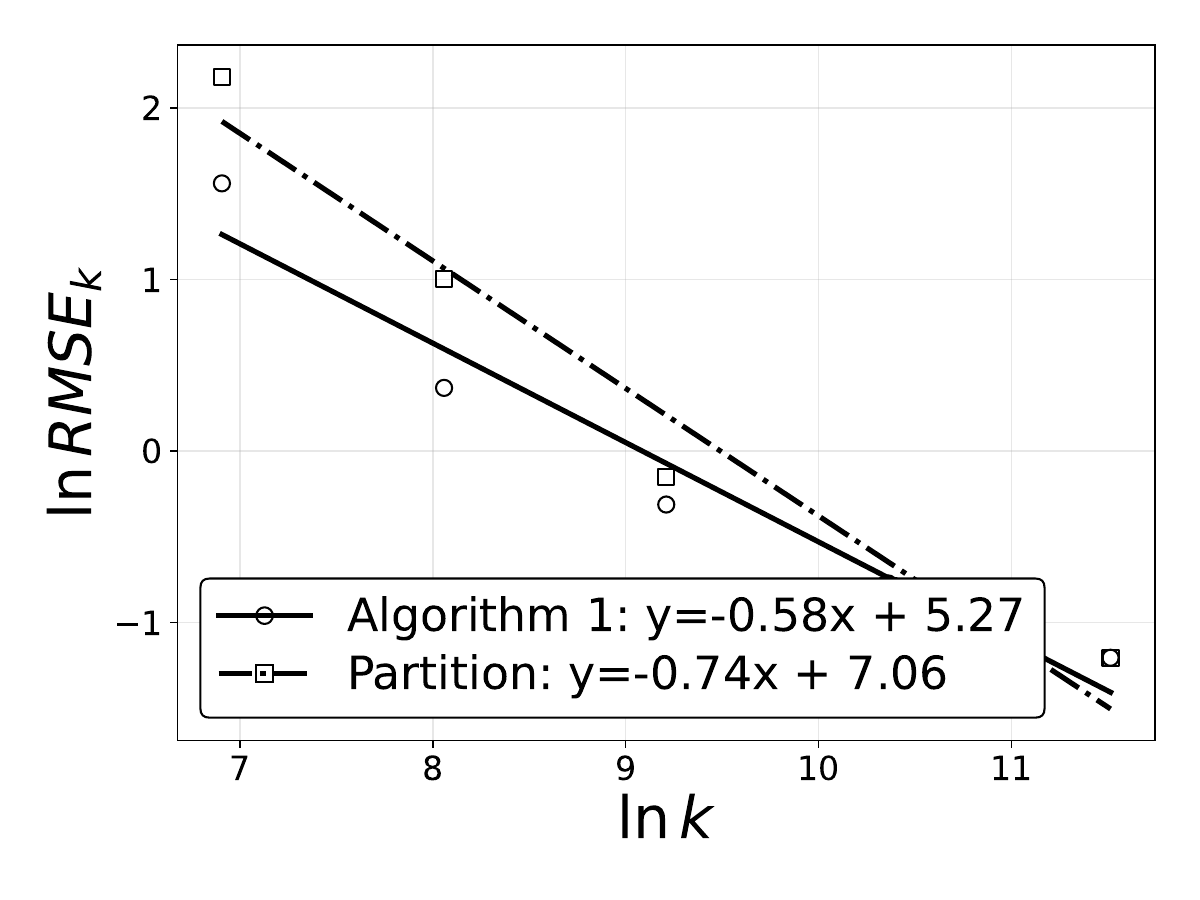}\par
            \vspace{3pt}\footnotesize (c) Case 2, \(\exp\), RVaR
        \end{minipage}\hfill
        \begin{minipage}{0.23\textwidth}
            \centering
            \includegraphics[width=\linewidth]{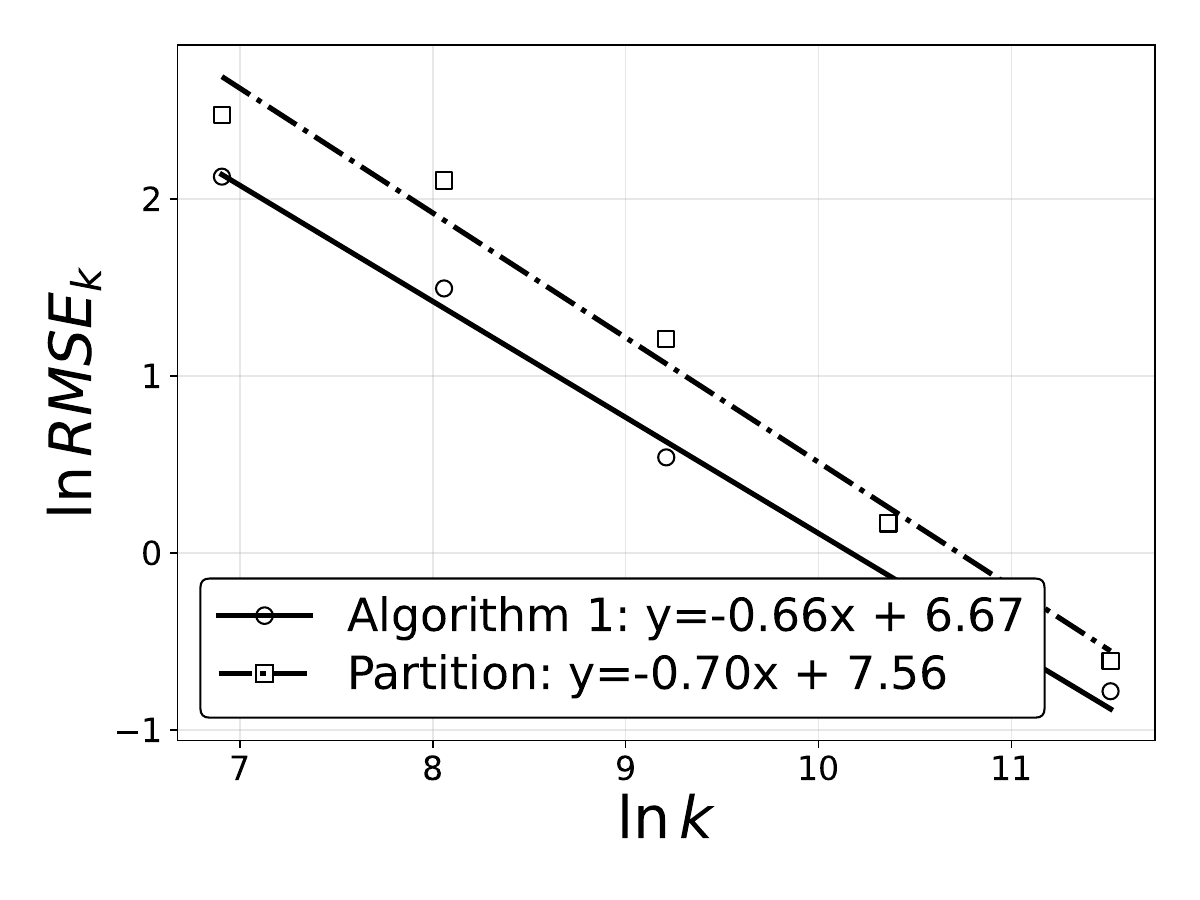}\par
            \vspace{3pt}\footnotesize (d) Case 2, \(t\), RVaR
        \end{minipage}

        \par\vspace{3ex}

        \begin{minipage}{0.23\textwidth}
            \centering
            \includegraphics[width=\linewidth]{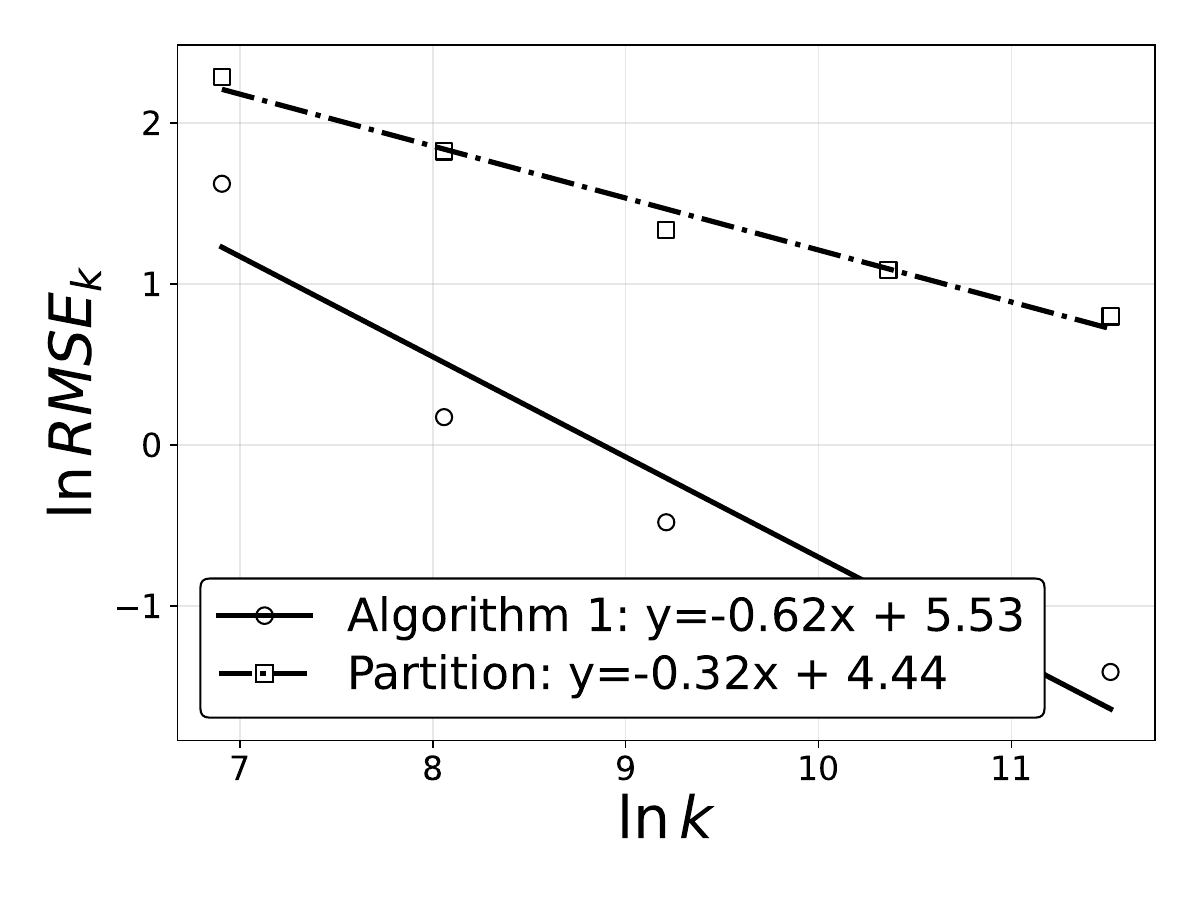}\par
            \vspace{3pt}\footnotesize (e) Case 1, $\exp$, M-G
        \end{minipage}\hfill
        \begin{minipage}{0.23\textwidth}
            \centering
            \includegraphics[width=\linewidth]{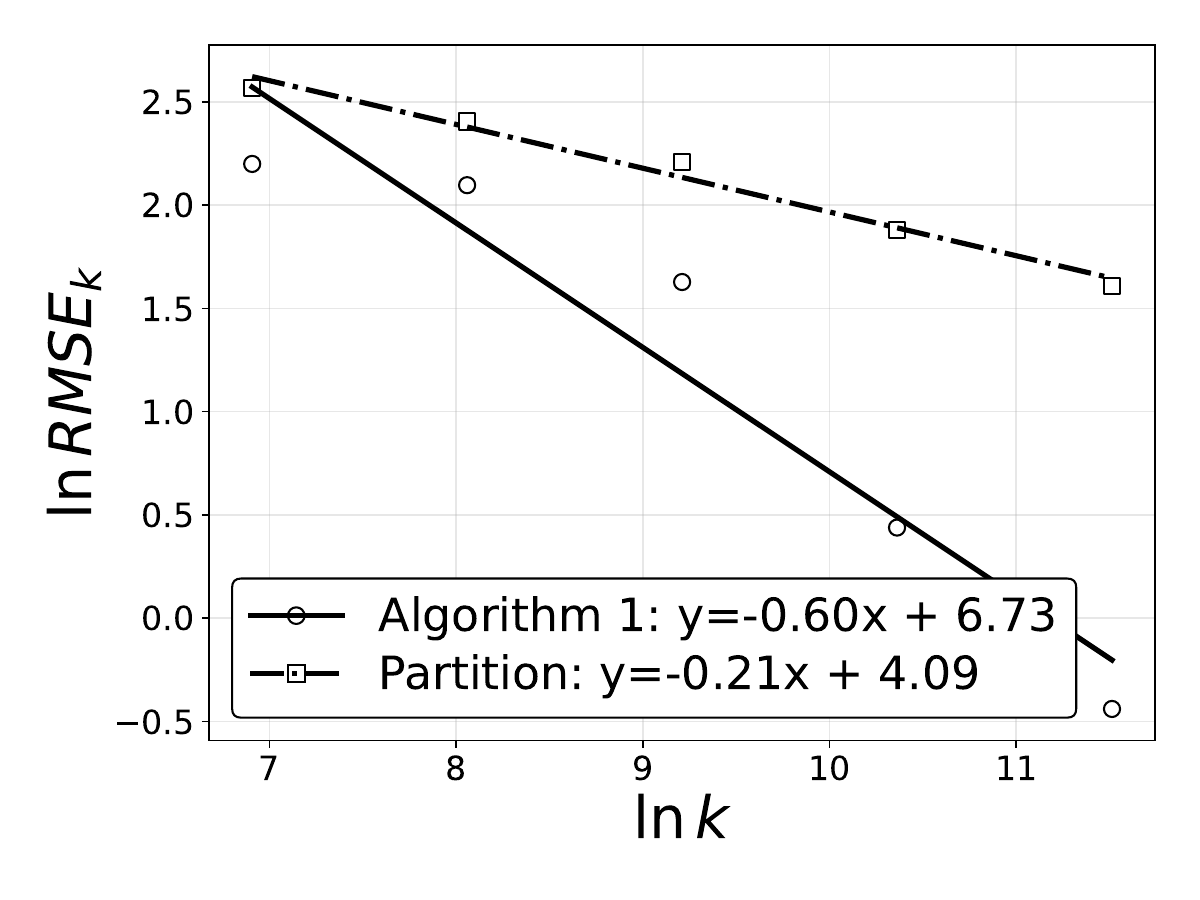}\par
            \vspace{3pt}\footnotesize (f) Case 1, $t$, M-G
        \end{minipage}\hfill
        \begin{minipage}{0.23\textwidth}
            \centering
            \includegraphics[width=\linewidth]{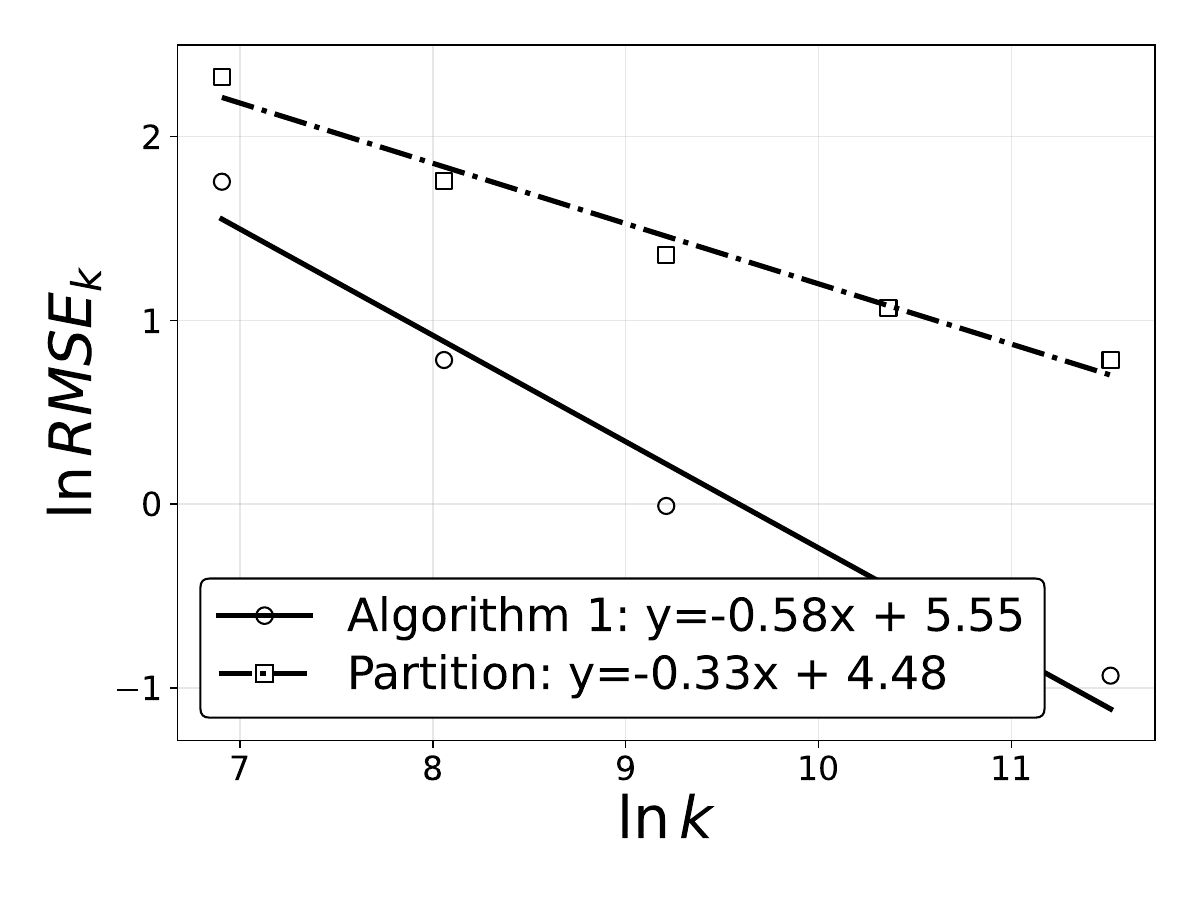}\par
            \vspace{3pt}\footnotesize (g) Case 2, $\exp$, M-G
        \end{minipage}\hfill
        \begin{minipage}{0.23\textwidth}
            \centering
            \includegraphics[width=\linewidth]{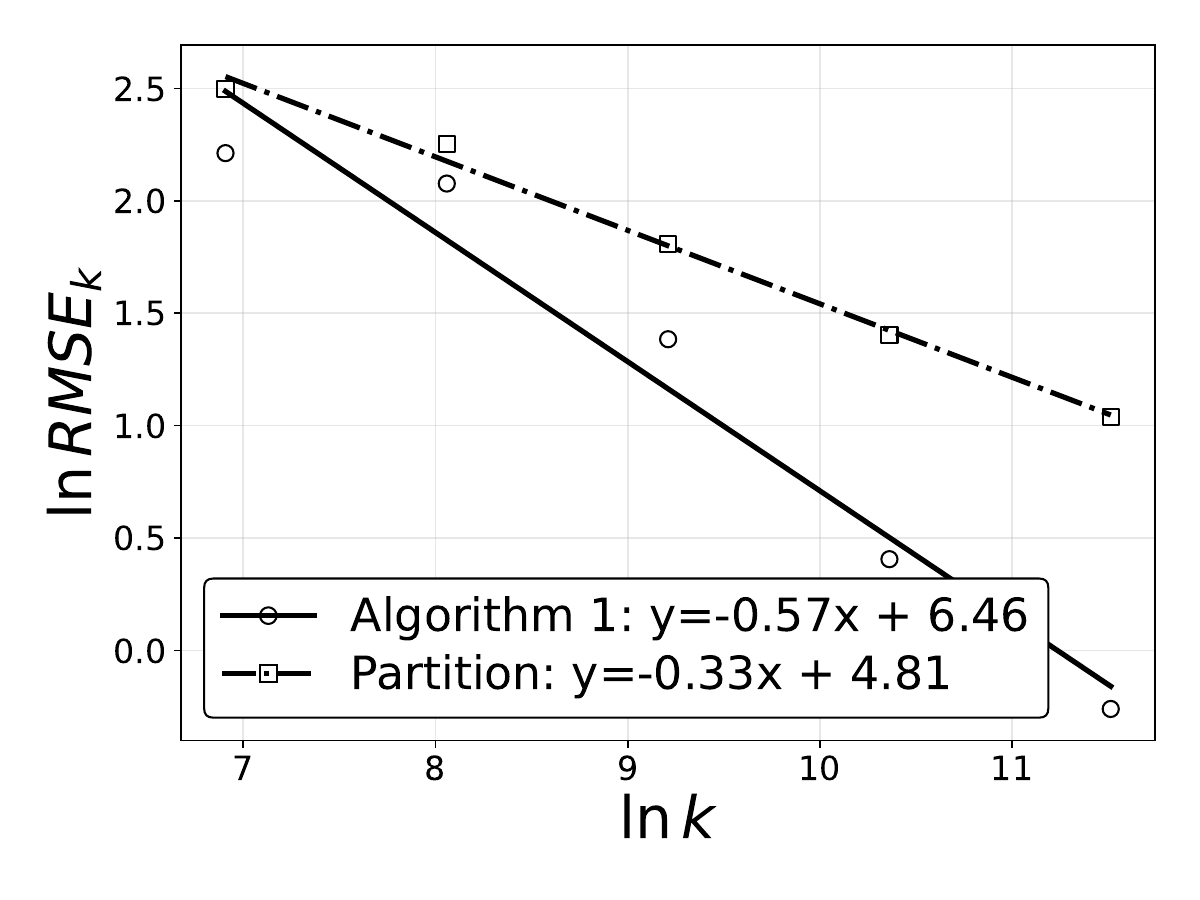}\par
            \vspace{3pt}\footnotesize (h) Case 2, $t$, M-G
        \end{minipage}

        \par\vspace{3ex}

        \begin{minipage}{0.23\textwidth}
            \centering
            \includegraphics[width=\linewidth]{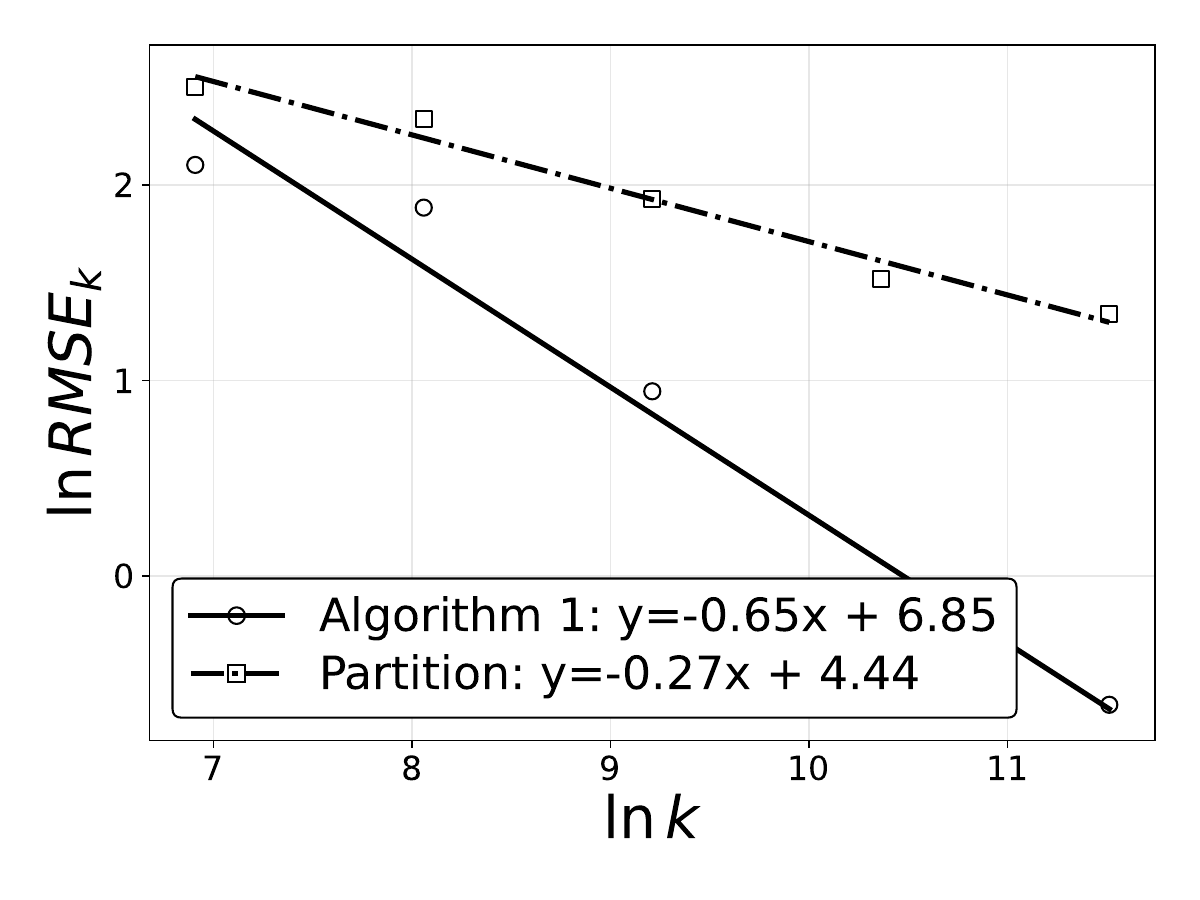}\par
            \vspace{3pt}\footnotesize (i) Case 1, $\exp$, G-S
        \end{minipage}\hfill
        \begin{minipage}{0.23\textwidth}
            \centering
            \includegraphics[width=\linewidth]{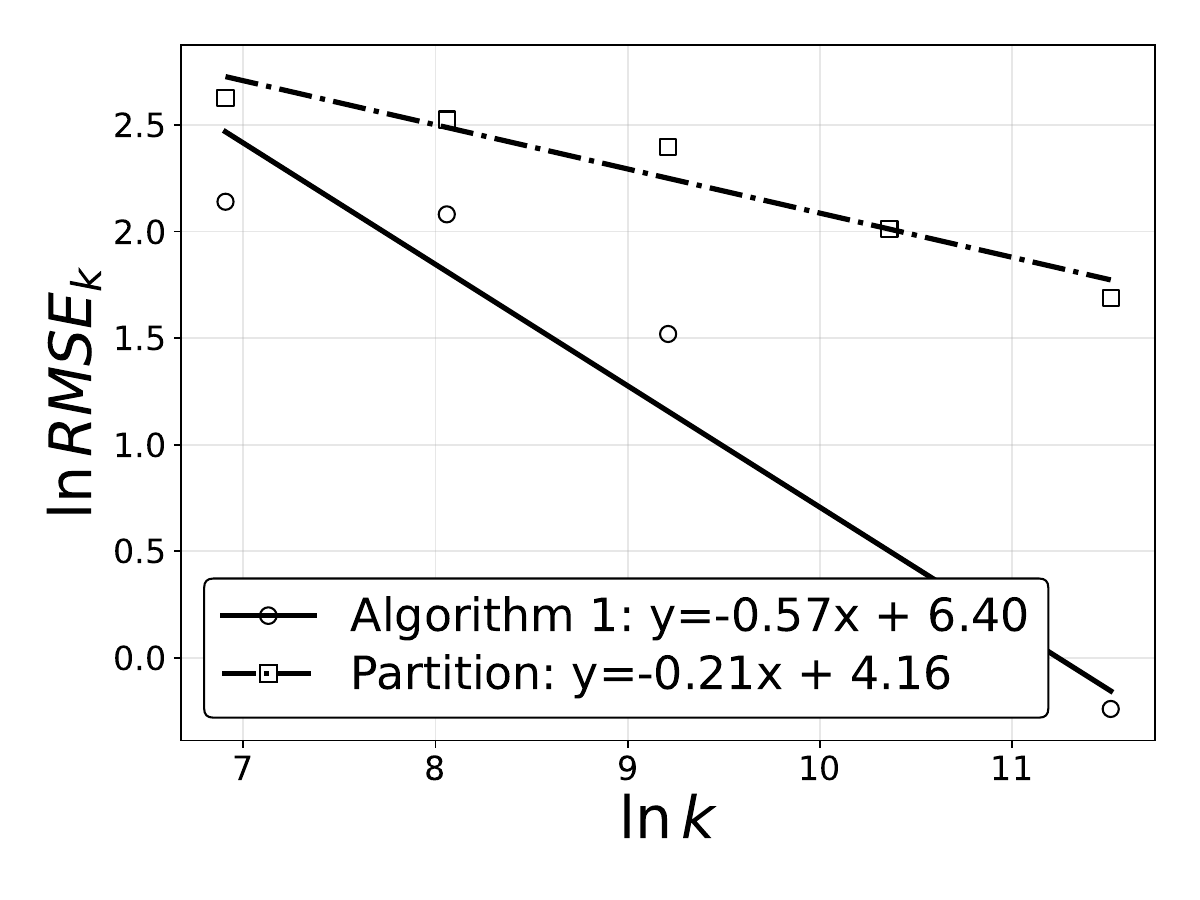}\par
            \vspace{3pt}\footnotesize (j) Case 1, $t$, G-S
        \end{minipage}\hfill
        \begin{minipage}{0.23\textwidth}
            \centering
            \includegraphics[width=\linewidth]{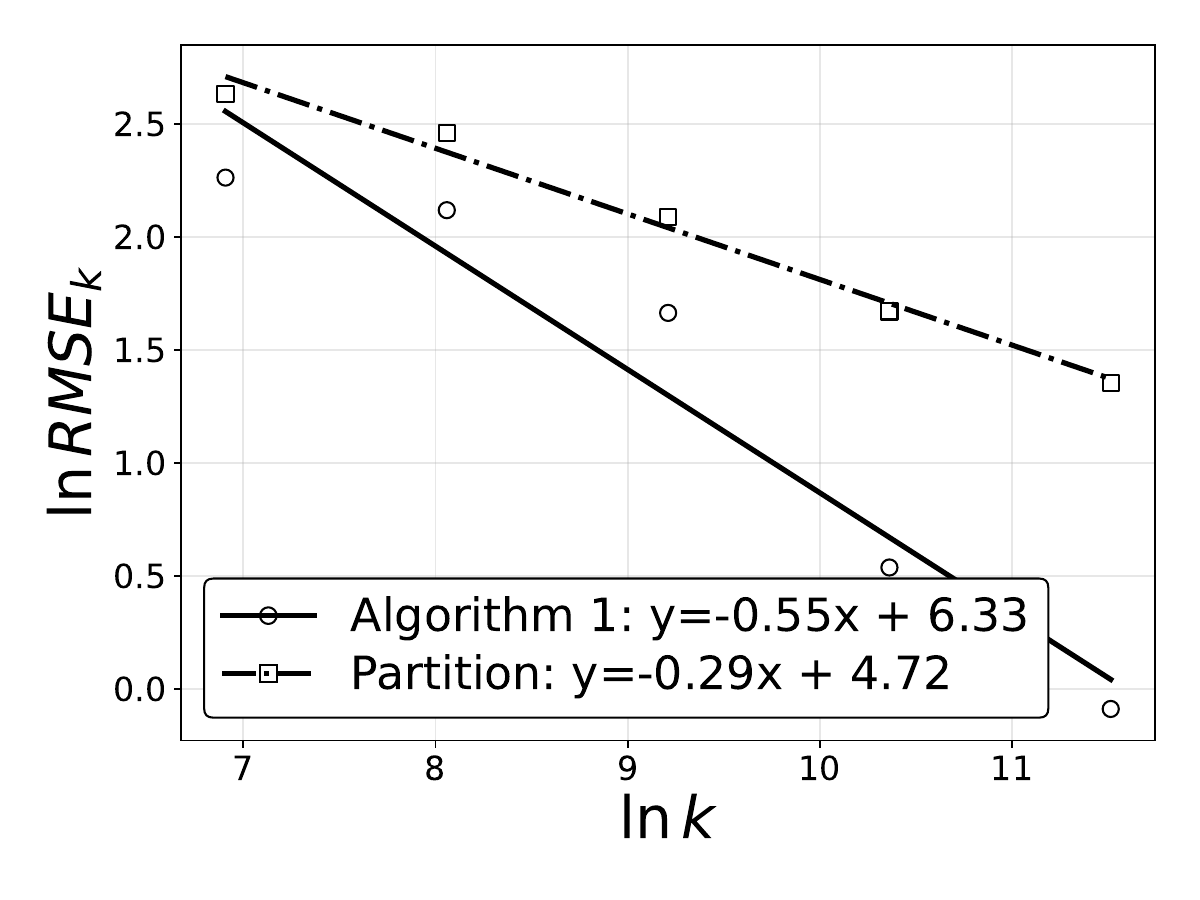}\par
            \vspace{3pt}\footnotesize (k) Case 2, $\exp$, G-S
        \end{minipage}\hfill
        \begin{minipage}{0.23\textwidth}
            \centering
            \includegraphics[width=\linewidth]{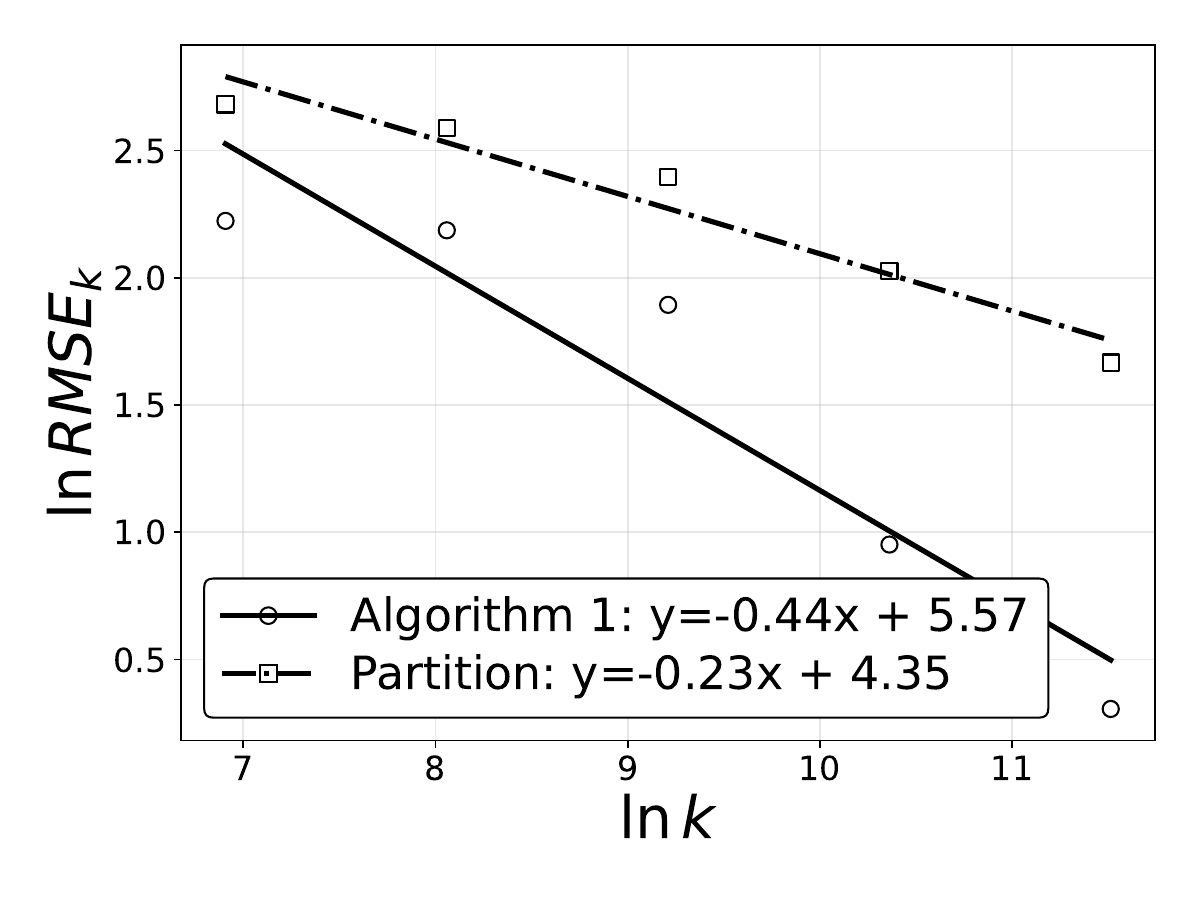}\par
            \vspace{3pt}\footnotesize (l) Case 2, $t$, G-S
        \end{minipage}

        \par\vspace{3ex}

        \begin{minipage}{0.23\textwidth}
            \centering
            \includegraphics[width=\linewidth]{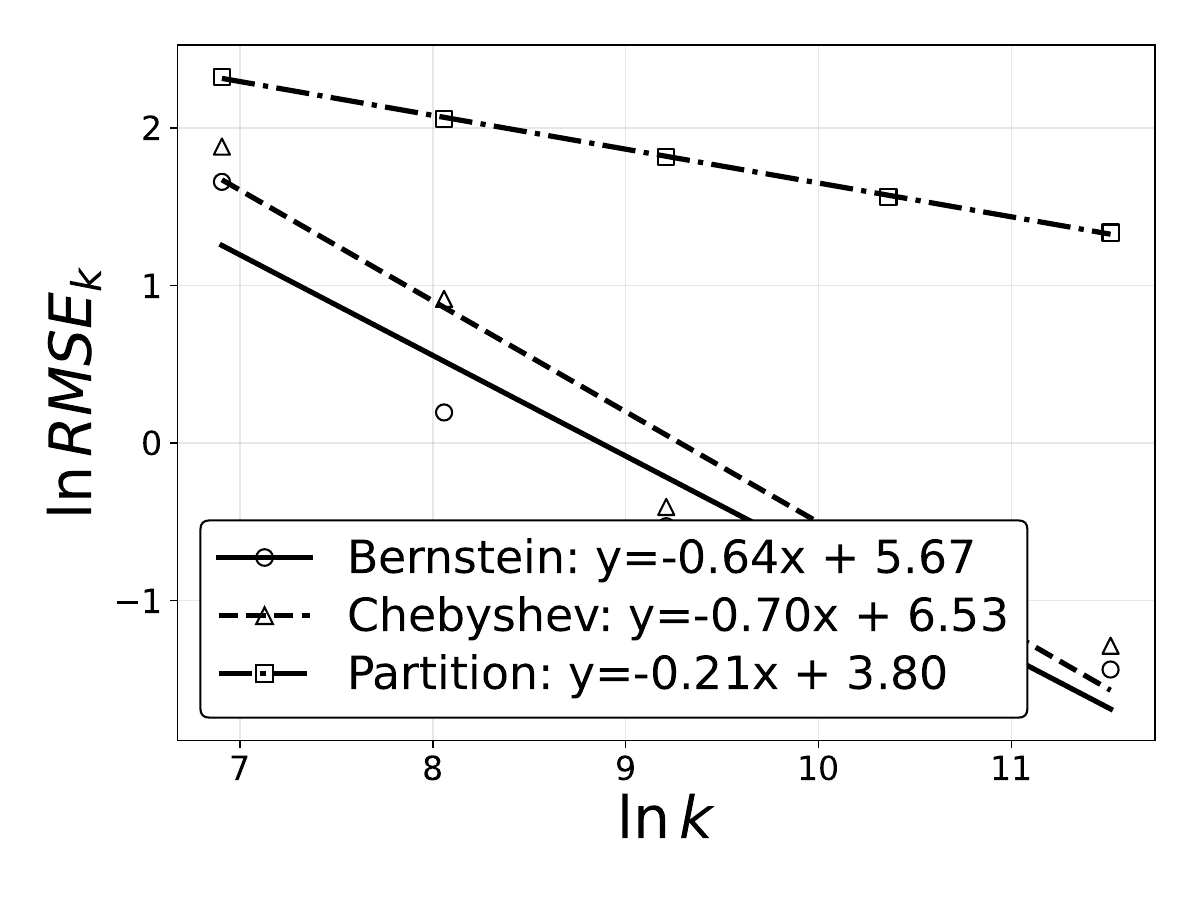}\par
            \vspace{3pt}\footnotesize (m) Case 1, $\exp$, Wang
        \end{minipage}\hfill
        \begin{minipage}{0.23\textwidth}
            \centering
            \includegraphics[width=\linewidth]{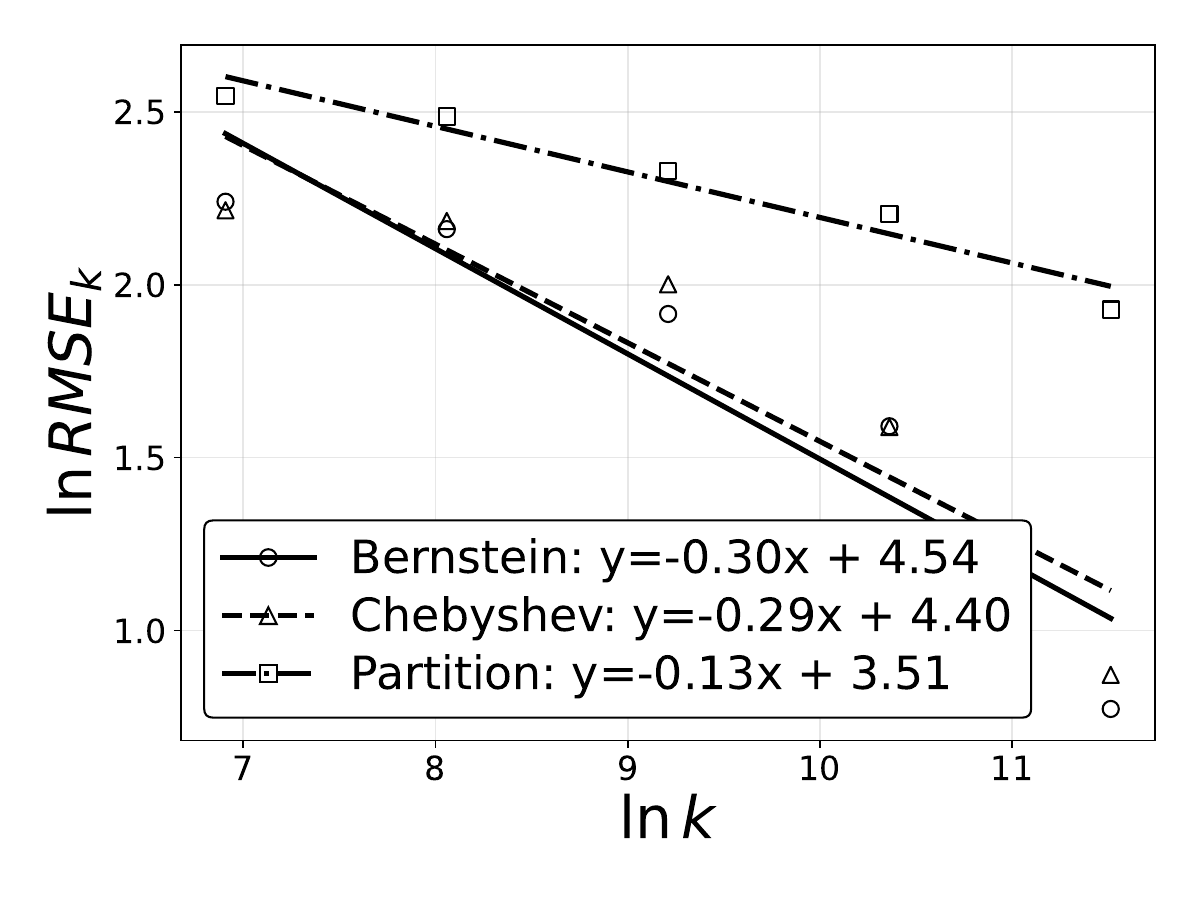}\par
            \vspace{3pt}\footnotesize (n) Case 1, $t$, Wang
        \end{minipage}\hfill
        \begin{minipage}{0.23\textwidth}
            \centering
            \includegraphics[width=\linewidth]{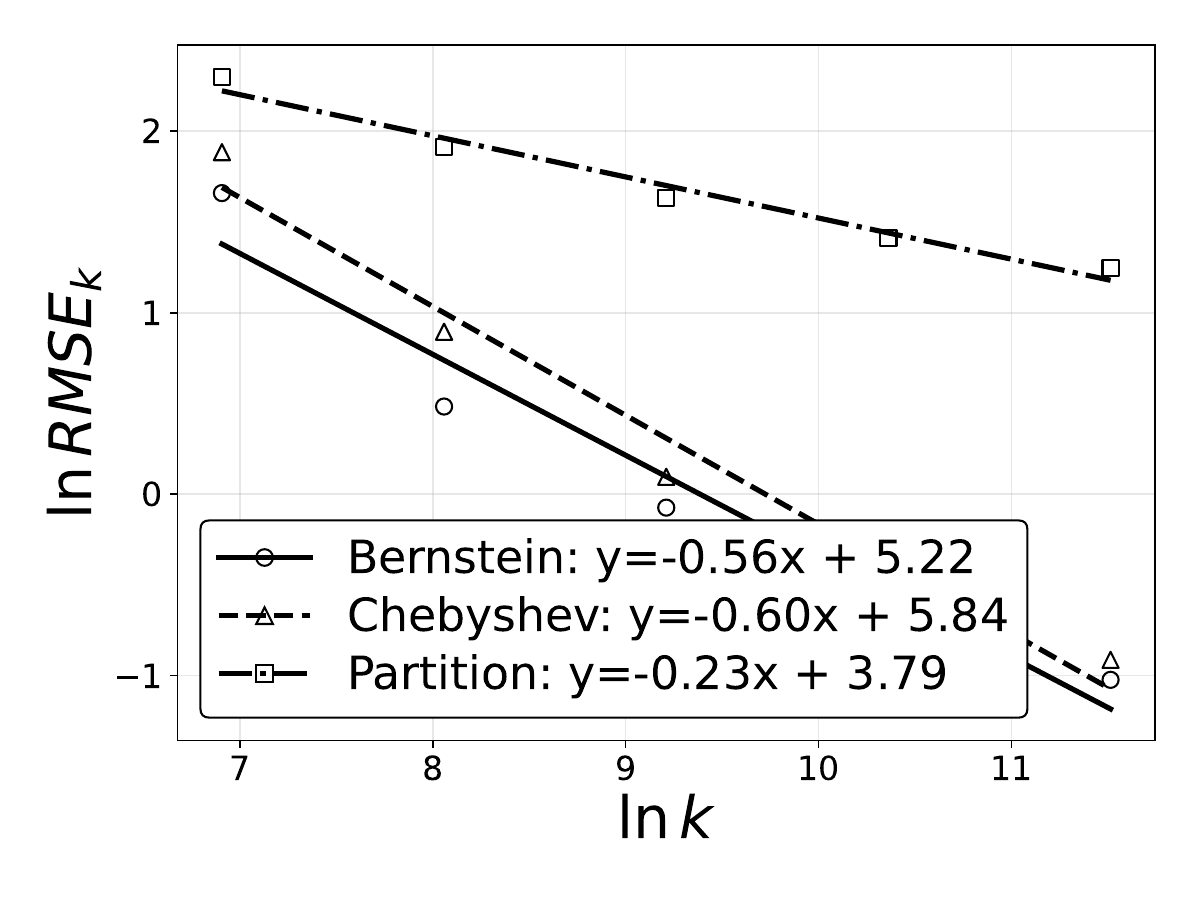}\par
            \vspace{3pt}\footnotesize (o) Case 2, $\exp$, Wang
        \end{minipage}\hfill
        \begin{minipage}{0.23\textwidth}
            \centering
            \includegraphics[width=\linewidth]{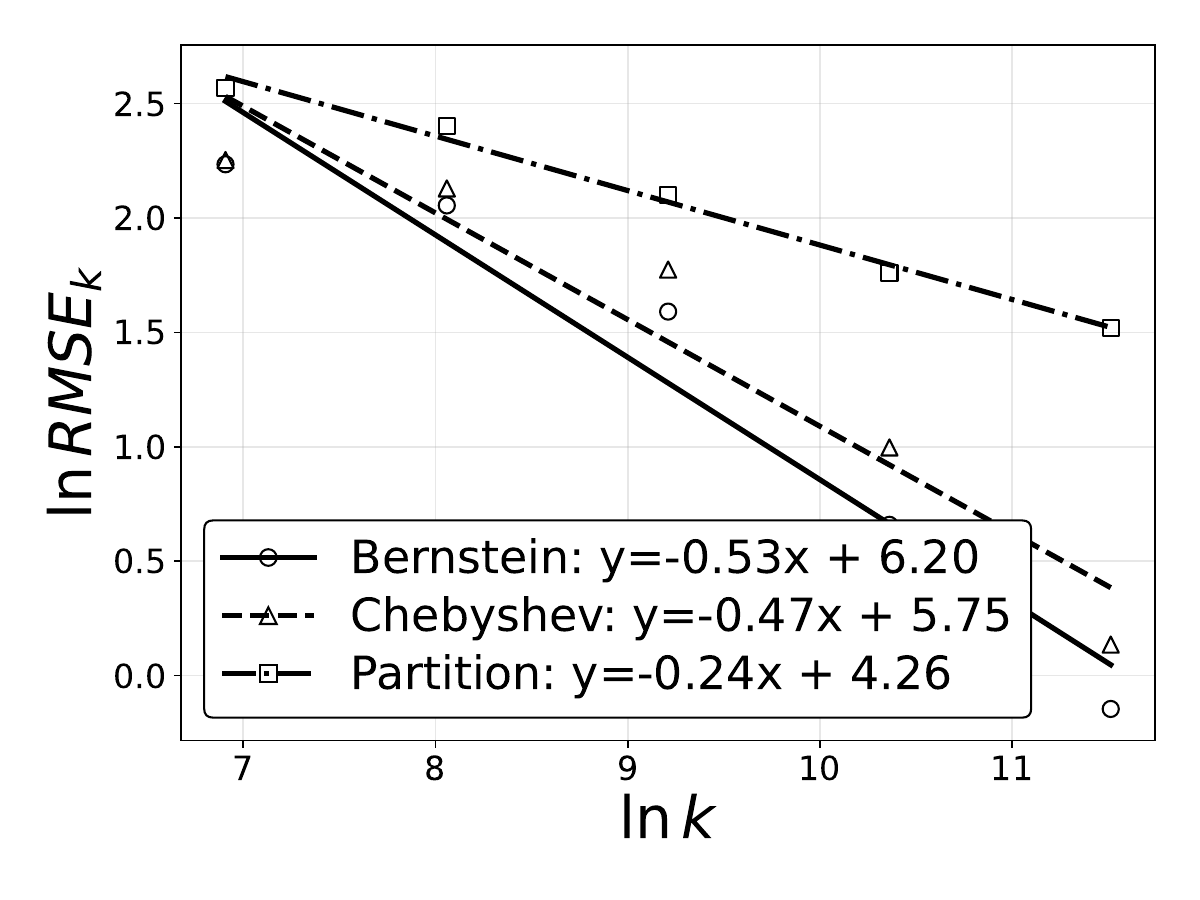}\par
            \vspace{3pt}\footnotesize (p) Case 2, $t$, Wang
        \end{minipage}
    \end{minipage}%
}

\begin{figure}
    \FIGURE
    {\usebox{\boxRateToy}}
    {Convergence Rates of the Empirical RMSEs of $\{\theta_k\}$ on Synthetic Test Cases.\label{fig:rate_Toy}}
    {Rows from top to bottom represent RVaR, mean-Gini (M-G), Gini shortfall (G-S), and the Wang transform (Wang), respectively.}
\end{figure}

\subsection{Barrier Option}\label{sec_option}
A financial barrier option grants the holder the right to purchase an underlying asset at a predetermined strike price $X$ at maturity $T$. We consider a discrete up-and-out European call option, which is nullified (knocked out) if the underlying asset price exceeds a specified barrier level $L$ at any monitoring time $t\in \mathcal{T}:= \{t_1, t_2, \dots, t_{M}\}$, where $t_j:=(j/M)T$.
We assume the asset price $S_t$ follows a geometric jump-diffusion process: $S_t = S_0 \exp[(r - \sigma^2/2)t + \sigma W_t + \sum_{i=1}^{N_t} J_i],$ where $S_0$ is the initial price, $r$ is the risk-free interest rate, $\sigma$ is the diffusion volatility, $W_t$ is a standard Brownian motion,
$\{N_t\}$ is a Poisson process with intensity $\lambda_J$, and  $J_i \sim \mathcal{N}(\mu_J, \sigma_J^2)$ denotes the $i$th jump size in the log-price.
The discounted terminal payoff $Z$ of the barrier option is given by
\begin{align*}
    Z = e^{-rT} \max\{S_T-X, 0\} \cdot \I{\max\nolimits_{1 \le j \le M} S_{t_j} < L}
    =\Lambda_{M}\I{0<\Lambda_{M}<b}
    \prod\nolimits_{j=1}^{M-1}\I{\Lambda_j< b},
\end{align*}
where $b:=e^{-rT}(L-X)$ and the discounted payoff $\Lambda_j:=e^{-rT}(S_{t_j}-X)=e^{-rT}\{S_0\exp[(r-\sigma^2/2)t_j+\sigma\sqrt{T/M}\sum_{\ell=1}^j\Upsilon_\ell+\sum_{i=1}^{N_{t_j}}J_i]-X\},$ with $\Upsilon_\ell=(W_{t_\ell}-W_{t_{\ell-1}})/\sqrt{T/M}$ being i.i.d. standard normal. Setting $L=\infty$ reduces the model to the vanilla European option $Z=e^{-rT} \max\{S_T-X, 0\}$.
We adopt parameter values from \cite{bakshi1997empirical}: $\mu_J=-0.05$, $\sigma_J=0.07$, $\lambda_J=0.59$, $r=0.05$, $\sigma=0.14$, $T=1$, $M=12$, $S_0=X=100$, and $L=130$.
Let \(\theta=S_0\in\Re\) be the parameter of interest and \(Y(\theta)=-Z\) the discounted terminal loss. At the fixed value \(\theta=100\), we estimate its mean-Gini risk with \(\lambda_{\operatorname{MG}}=1\) and the sensitivity \(\partial\rho_g(\theta)/\partial\theta\). Here \(S_0\) is an exogenous simulation input varied for sensitivity analysis.
This sensitivity quantifies the option's local risk exposure to movements in the initial asset price under the mean-Gini criterion.
It accounts for both changes in the potential terminal payoff and changes in the probability of being knocked out.
Since the barrier payoff has a point mass at zero, this example also provides a useful robustness test for the proposed estimation and sensitivity recursions.

For comparison with the ground truth, assuming the underlying distribution is known, the DRM and its sensitivity can be estimated by combining Equations~\eqref{drm} and \eqref{eq_drho} with the generalized likelihood ratio estimators \citep{glynn2021computing}:
\begin{align*}
    \bar{F}(y;\theta)&=\bE[\I{Z< -y}]\\
    & =\bE\Big[
    \I{\Lambda_{M}\le -y}
    \prod\nolimits_{j=1}^{M-1}\I{\Lambda_j< b}
    +\sum\nolimits_{j=1}^{M}
    \prod\nolimits_{\ell=1}^{j-1}\I{\Lambda_\ell< b}
    \I{\Lambda_j>b}
    \Big],\\
    \frac{\partial \bar{F}(y;\theta)}{\partial\theta}&=
    \bE\Big[\Big(\I{\Lambda_{M}\le -y}
    \prod\nolimits_{j=1}^{M-1}\I{\Lambda_j< b}\\
    &\qquad{}+\sum\nolimits_{j=1}^{M}
    \prod\nolimits_{\ell=1}^{j-1}\I{\Lambda_\ell< b}
    \I{\Lambda_j>b}\Big)\frac{\Upsilon_1}{\theta\sigma\sqrt{T/M}}
    \Big],
\end{align*}
for \(-b\le y<0\), with \(\bar F(0;\theta)=0\). This yields benchmark estimates $\rho_g(\theta)\approx-1.7996$ and $\partial\rho_g(\theta)/\partial\theta\approx -0.13185$ based on $10^8$ sample paths.

Because the DRM, VaR, and their gradients may differ in scale, we employ separate step-sizes in Procedures~\eqref{proc_q}--\eqref{proc_theta}: $\gamma_{k,q}=\gamma_{k,D_q}=10^3/k^{0.999}$, $\gamma_{k,\rho}=\gamma_{k,G}=1/(k^{0.999}+10^3)$, $\alpha_k\equiv0$, and $c_k=10/k^{1/6}$.
The initial values are set as $\theta_0=100$ and $\rho_0=q_0(\cdot)=G_0=D_{q,0}(\cdot)=0$. A projection operator is applied in Procedures~\eqref{proc_q}--\eqref{proc_drho} to constrain the estimates within the bounded set $\{x:\|x\|\le50\}$. In addition to the standard sensitivity estimators (\ref{proc_dq}) and
(\ref{proc_drho}), we implement common random number (CRN) versions of these recursions to evaluate their effectiveness in reducing the variance of the gradient estimators \citep{hu2025quantile}.

Figures~\ref{fig:barrier_convergence}(a)--\ref{fig:barrier_convergence}(c) show the empirical mean absolute errors (MAEs) of the estimators over $40$ independent replications, i.e., $(1/40)\sum_{r=1}^{40}|\rho_k(\omega_r)-\rho_g(\theta)|$, as functions of the number of iterations.
Figures~\ref{fig:barrier_convergence}(d)--\ref{fig:barrier_convergence}(f) illustrate the empirical convergence rates of the DRM and its sensitivity estimators. Table~\ref{tab:barrier_results} reports the average running times and descriptive statistics of the final estimates. Algorithm~\ref{alg} achieves smaller MAEs and RMSEs than the partition-based method for both the DRM and its sensitivity, with running times of 147--148 seconds compared with 839--842 seconds. For Algorithm~\ref{alg}, CRN reduces the sensitivity RMSE from \(0.013\) to \(0.008\).

\newsavebox{\boxBarrierConvergence}
\sbox{\boxBarrierConvergence}{%
    \begin{minipage}{\textwidth}
        \centering

        \begin{minipage}{0.32\textwidth}
            \centering
            \includegraphics[width=\linewidth]{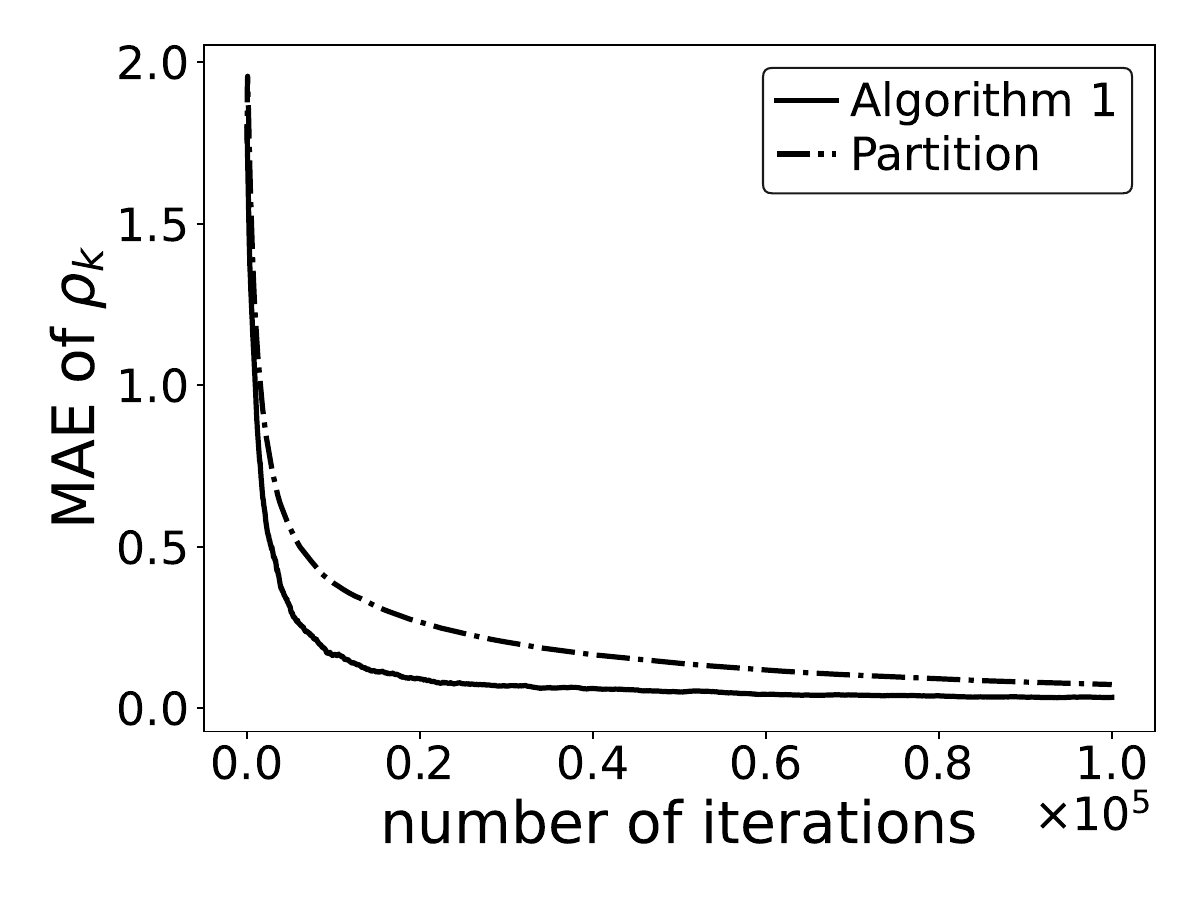}\par
            \vspace{3pt}\footnotesize (a) Mean-Gini Estimation
        \end{minipage}\hfill
        \begin{minipage}{0.32\textwidth}
            \centering
            \includegraphics[width=\linewidth]{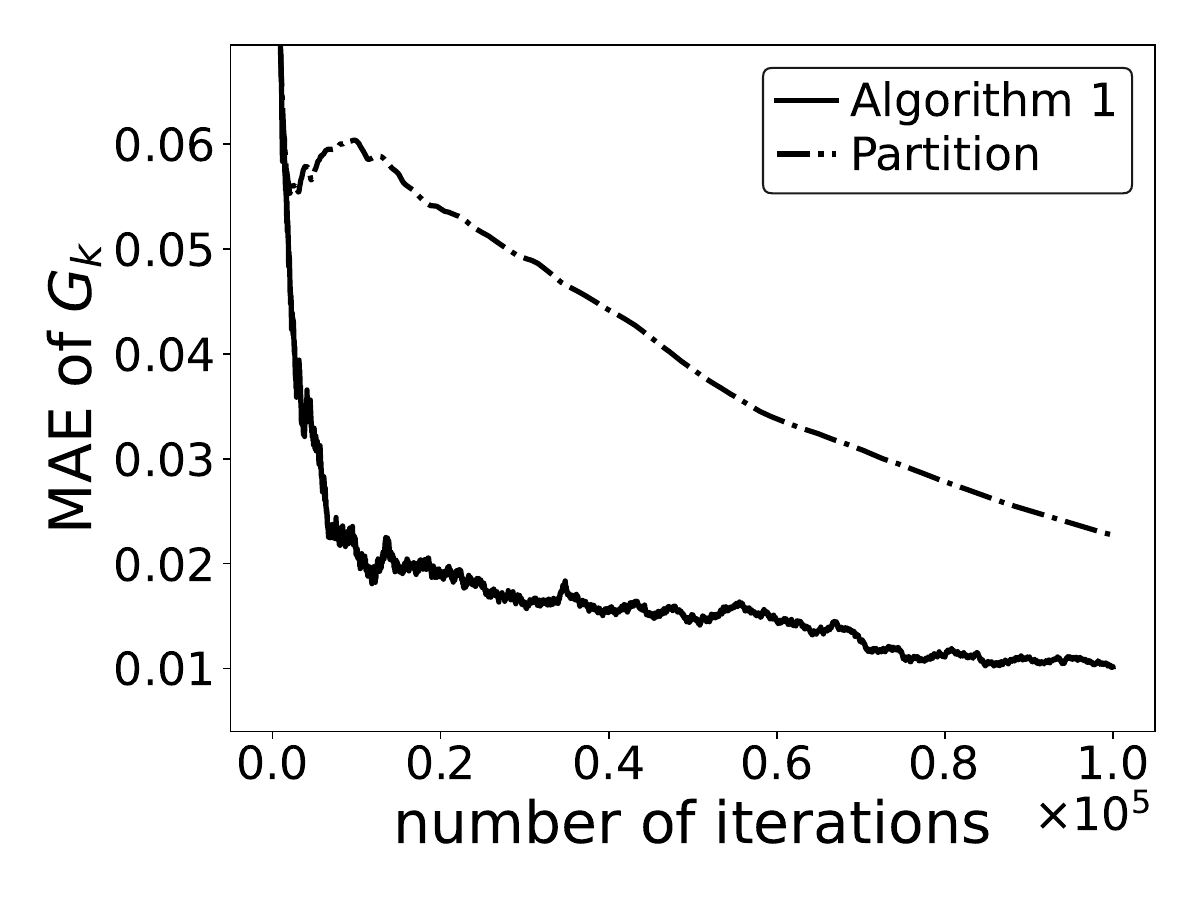}\par
            \vspace{3pt}\footnotesize (b) Sensitivity (No CRN)
        \end{minipage}\hfill
        \begin{minipage}{0.32\textwidth}
            \centering
            \includegraphics[width=\linewidth]{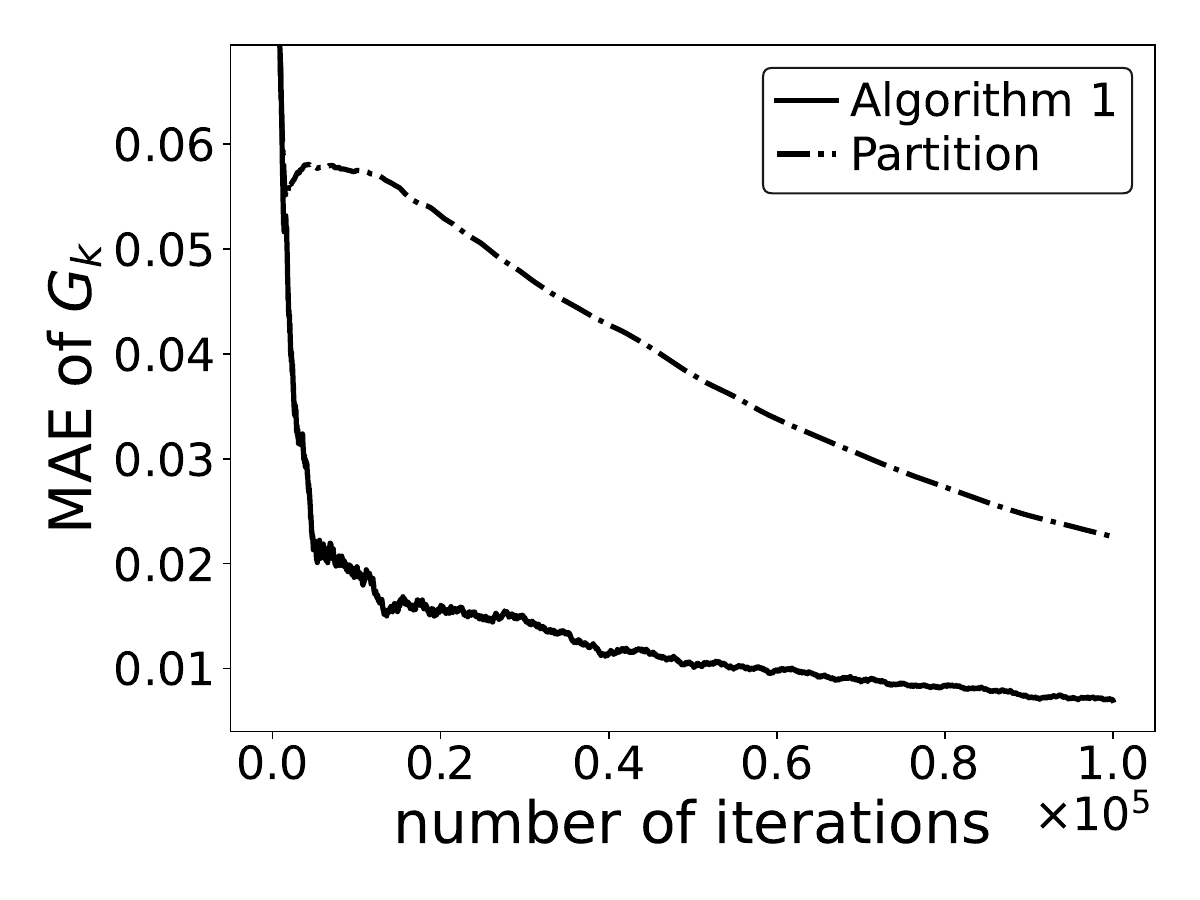}\par
            \vspace{3pt}\footnotesize (c) Sensitivity (CRN)
        \end{minipage}

        \par\vspace{3ex}

        \begin{minipage}{0.32\textwidth}
            \centering
            \includegraphics[width=\linewidth]{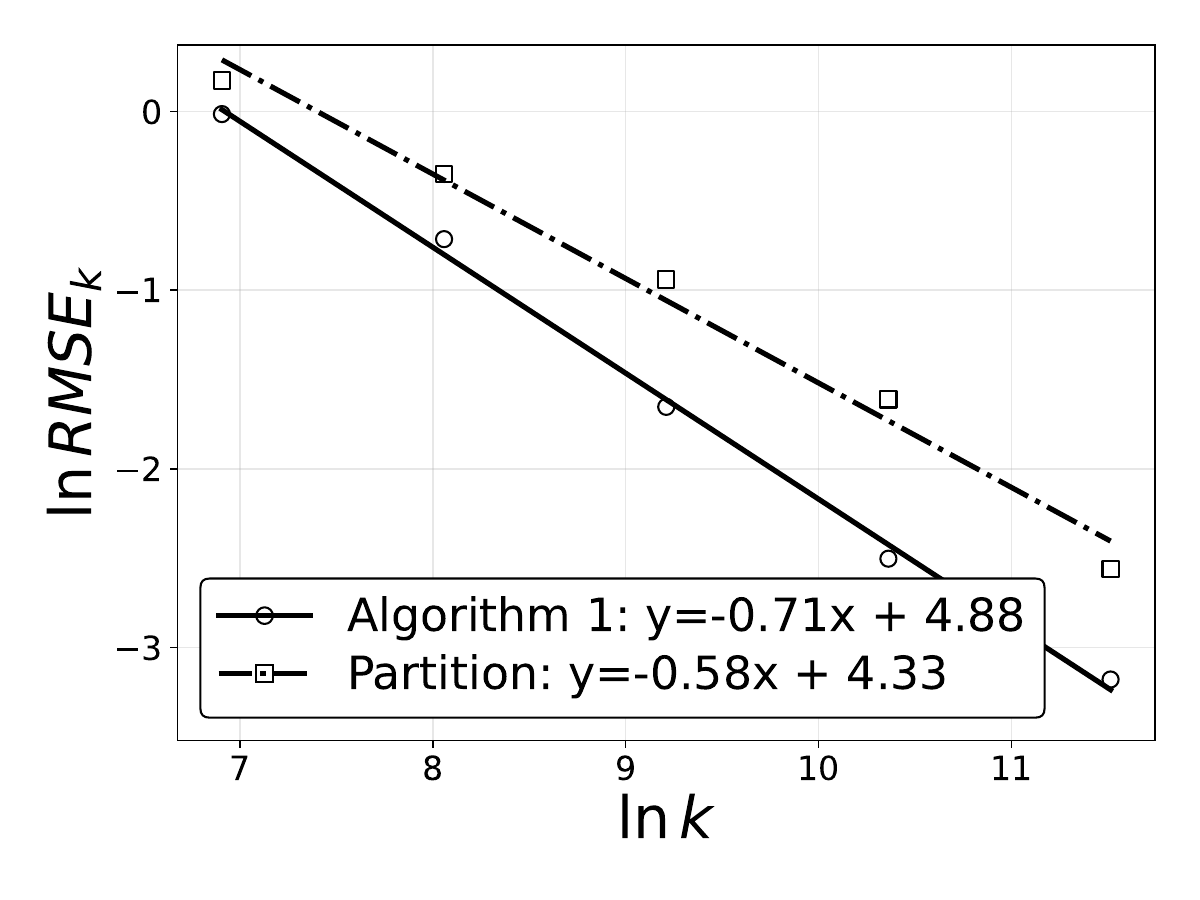}\par
            \vspace{3pt}\footnotesize (d) Mean-Gini Estimation
        \end{minipage}\hfill
        \begin{minipage}{0.32\textwidth}
            \centering
            \includegraphics[width=\linewidth]{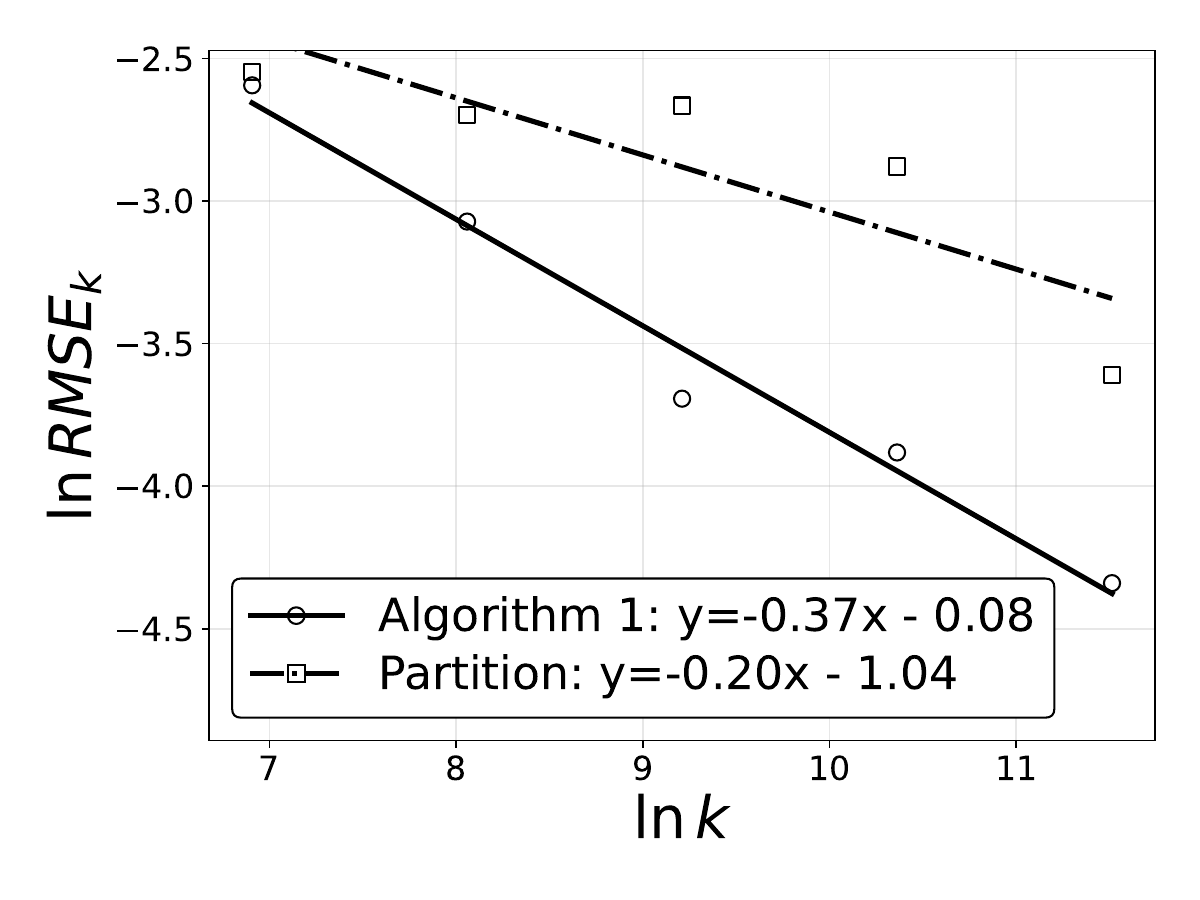}\par
            \vspace{3pt}\footnotesize (e) Sensitivity (No CRN)
        \end{minipage}\hfill
        \begin{minipage}{0.32\textwidth}
            \centering
            \includegraphics[width=\linewidth]{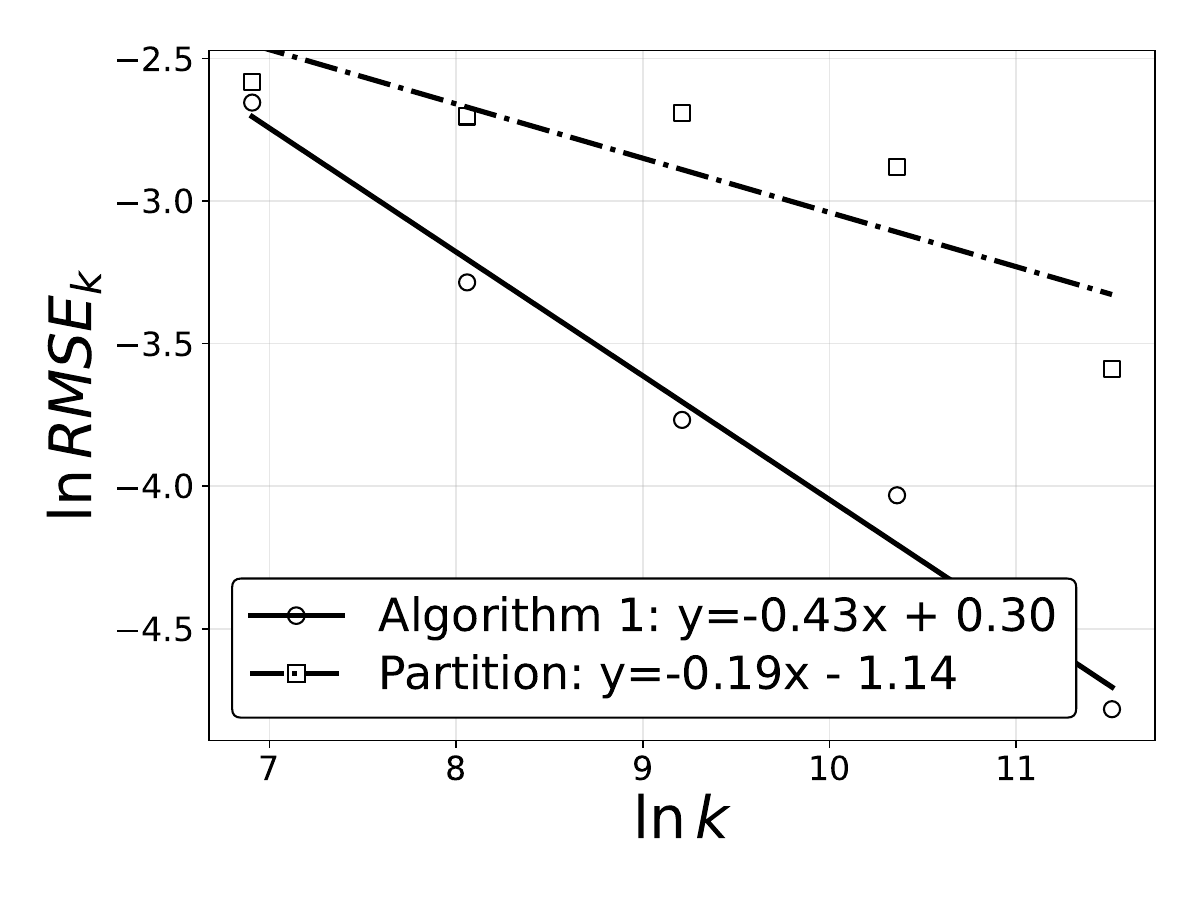}\par
            \vspace{3pt}\footnotesize (f) Sensitivity (CRN)
        \end{minipage}

    \end{minipage}%
}

\begin{figure}
    \FIGURE
    {\usebox{\boxBarrierConvergence}}
    {Performance of the Proposed Algorithm and the Partition-Based Method on Barrier Option.\label{fig:barrier_convergence}}
    {The top row displays the MAEs for $\rho_k$ and $G_k$ (without and with CRN). The bottom row shows the corresponding log-log plots for the RMSEs.}
\end{figure}

\begin{table}
\centering
\TABLE
{Comparison Between Algorithm~\ref{alg} and the Partition-Based Method on Barrier Option.\label{tab:barrier_results}}
{\setlength{\tabcolsep}{4pt}\footnotesize
\begin{tabular*}{0.88\linewidth}{@{}l @{\extracolsep{\fill}} *{4}{c}@{}}
\toprule
& \multicolumn{2}{c}{No CRN} & \multicolumn{2}{c}{CRN} \\
\cmidrule(lr){2-3} \cmidrule(lr){4-5}
& Algorithm 1 & Partition-based & Algorithm 1 & Partition-based \\
\midrule
ART (s) & $\mathbf{147}$ & $839$ & $\mathbf{148}$ & $842$ \\
Bias of $\rho_K$ & $\mathbf{0.011}$ & $0.070$ & $\mathbf{0.024}$ & $0.074$ \\
MAE of $\rho_K$ & $\mathbf{0.028}$ & $0.070$ & $\mathbf{0.034}$ & $0.074$ \\
RMSE of $\rho_K$ & $\mathbf{0.035}$ & $0.074$ & $\mathbf{0.042}$ & $0.077$ \\
\midrule
Bias of $G_K$ & $\mathbf{-0.002}$ & $0.005$ & $0.006$ & $\mathbf{0.004}$ \\
MAE of $G_K$ & $\mathbf{0.010}$ & $0.023$ & $\mathbf{0.007}$ & $0.023$ \\
RMSE of $G_K$ & $\mathbf{0.013}$ & $0.027$ & $\mathbf{0.008}$ & $0.028$ \\
\bottomrule
\end{tabular*}}
{\textit{Note.} Computational results based on 40 independent replications. ART denotes the average running time. Bold indicates the best displayed value within each CRN setting; ties are all bolded.}
\end{table}

\subsection{Portfolio Optimization}\label{sec_portfolio}
We consider a portfolio comprising three primary industries: consumer, healthcare, and manufacturing. Let $X\in \Re^3$ denote the investment losses of these assets, which are assumed to follow a multivariate location-scale $t_3$ distribution: $X=\mu+\Sigma^{1/2}\xi$, where $\xi\in\Re^3$ follows a standard multivariate $t$-distribution with 3 degrees of freedom {and identity scale matrix}. {We estimate \(\mu\) and the annualized sample covariance \(\widehat C\) from daily returns in the \href{https://mba.tuck.dartmouth.edu/pages/faculty/ken.french/data_library.html}{Kenneth R. French Data Library}. Since \(\operatorname{Cov}(\xi)=3I_3\), we set \(\Sigma=\widehat C/3\). In the order consumer, healthcare, and manufacturing, the rounded estimates are}:
\[\mu = \begin{pmatrix}
-0.114 \\
-0.121 \\
-0.110
\end{pmatrix},
\quad
{\widehat C} = \begin{pmatrix}
0.0267 & 0.0217 & 0.0254 \\
0.0217 & 0.0305 & 0.0224 \\
0.0254 & 0.0224 & 0.0313
\end{pmatrix}.\]

We consider the risk-averse portfolio optimization problem {\(\min_{\theta\in\Theta}\rho_g(\theta)\)}, where $\rho_g(\theta)$ denotes the Gini shortfall of the portfolio loss $Y(\theta)=\theta_1X_1+\theta_2X_2+(1-\theta_1-\theta_2)X_3$ with $\lambda_{\operatorname{GS}}=1$ and $\varphi=0.75$. The feasible set is the solid unit simplex  $\Theta=\{\theta\in\Re^2:\sum_{i=1}^2\theta_i\le 1,\theta_1,\theta_2\ge0\}$.
The analytical DRM function, which serves as the ground truth, is {\(\rho_g(\theta)=w(\theta)^\intercal\mu+\rho_g[Z]\sqrt{w(\theta)^\intercal\Sigma w(\theta)}\)}, where \(w(\theta)=(\theta_1,\theta_2,1-\theta_1-\theta_2)^\intercal\) and {\(Z\) is a univariate standard \(t_3\) random variable. The value of \(\rho_g[Z]\) is given in Table~\ref{tab:base_risks}.}

In the implementation of Algorithm~\ref{alg}, we select step- and perturbation-sizes $\alpha_k=5/k^{0.999}$, $\gamma_{k,q}=\gamma_{k,D_q}=10/k^{3/4}$, $\gamma_{k,\rho}=\gamma_{k,G}=1/k^{3/4}$, and $c_k=1/(k+10^3)^{1/8}$.  The algorithm is initialized with $\rho_0=q_0(\cdot)=0$, $G_0=D_{q,0}(\cdot)=(0,0)^\intercal$, $\theta_0=(0,0)^\intercal$. A projection operator is incorporated into Procedures~\eqref{proc_q} and \eqref{proc_rho} to constrain the estimates within the bounded set $\{x:\|x\|\le1\}$.
We compare Algorithm~\ref{alg} with the L-statistic benchmark, both with and without CRN, using \(40\) independent replications and a common budget of \(600{,}000\) simulator evaluations. 
Implementation details of the L-statistic benchmark are given in Section~F of the Appendix.

Figures~\ref{fig:portfolio_convergence}(a)--\ref{fig:portfolio_convergence}(b) plot the average objective values against the cumulative number of simulator evaluations. Algorithm~\ref{alg} approaches the optimum more rapidly than the L-statistic benchmark under both sampling schemes. Table~\ref{tab:portfolio_convergence} reports the final results. Without CRN, the decision RMSE is \(0.143\) for Algorithm~\ref{alg} and \(0.314\) for the L-statistic benchmark. With CRN, these values decrease to \(0.106\) and \(0.202\), respectively. Algorithm~\ref{alg} also achieves a smaller RMSE of the objective value in both comparisons.

Figures~\ref{fig:portfolio_convergence}(c)--\ref{fig:portfolio_convergence}(d) show the empirical decision RMSEs on logarithmic scales. The fitted slopes are approximately \(-0.23\) for Algorithm~\ref{alg} and \(-0.09\) for the L-statistic benchmark under both sampling schemes. CRN reduces the empirical error levels for both methods, with little change in the fitted slopes. These results illustrate the tradeoff between batch estimation and decision updates under a limited simulation budget. As the L-statistic batch size grows, each update consumes more simulator evaluations. Algorithm~\ref{alg} uses a fixed number of outputs per iteration and can therefore perform more decision updates within the same budget. 

\newsavebox{\boxPortfolioConvergence}
\sbox{\boxPortfolioConvergence}{%
     \begin{minipage}{\textwidth}
        \centering

        \begin{minipage}{0.24\textwidth}
            \centering
            \includegraphics[width=\linewidth]{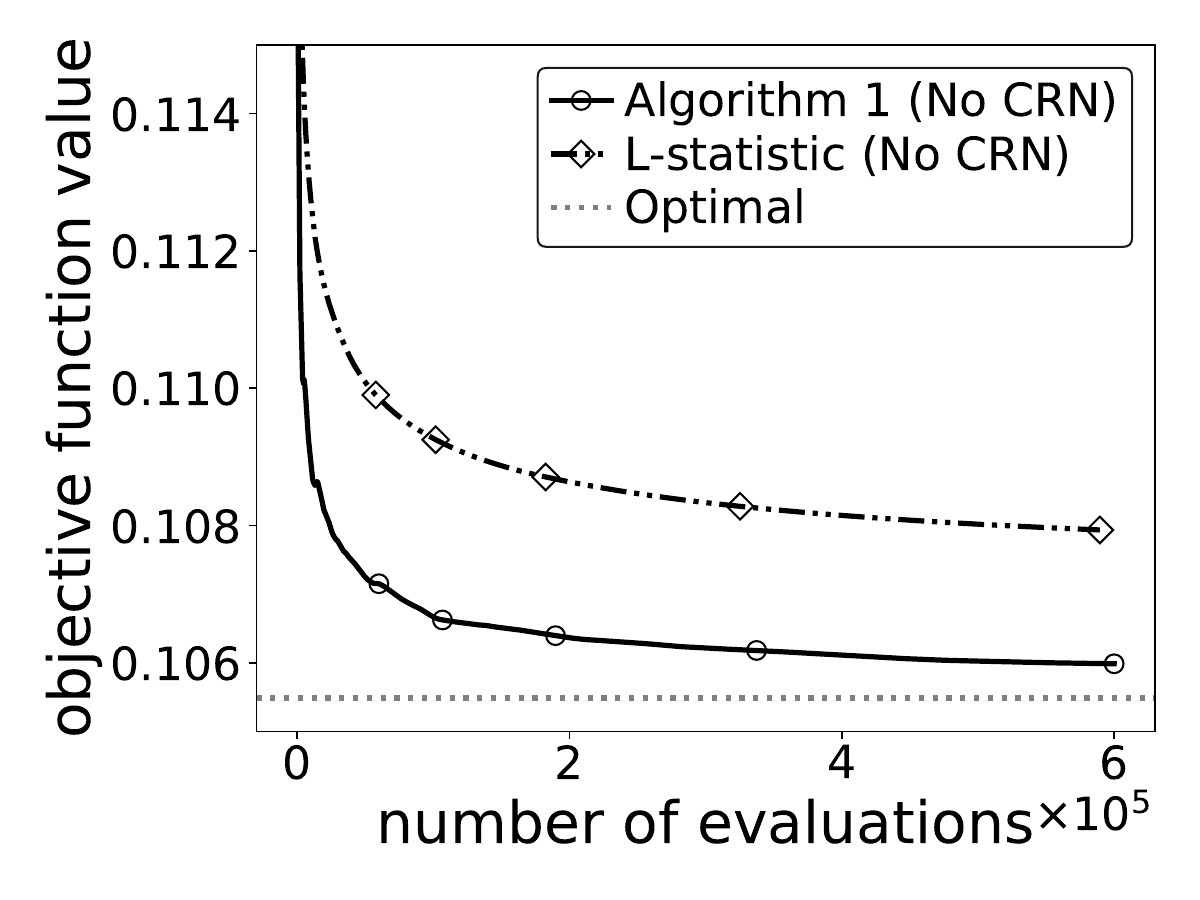}\par
            \vspace{3pt}\footnotesize (a) Objective (No CRN)
        \end{minipage}\hfill
        \begin{minipage}{0.24\textwidth}
            \centering
            \includegraphics[width=\linewidth]{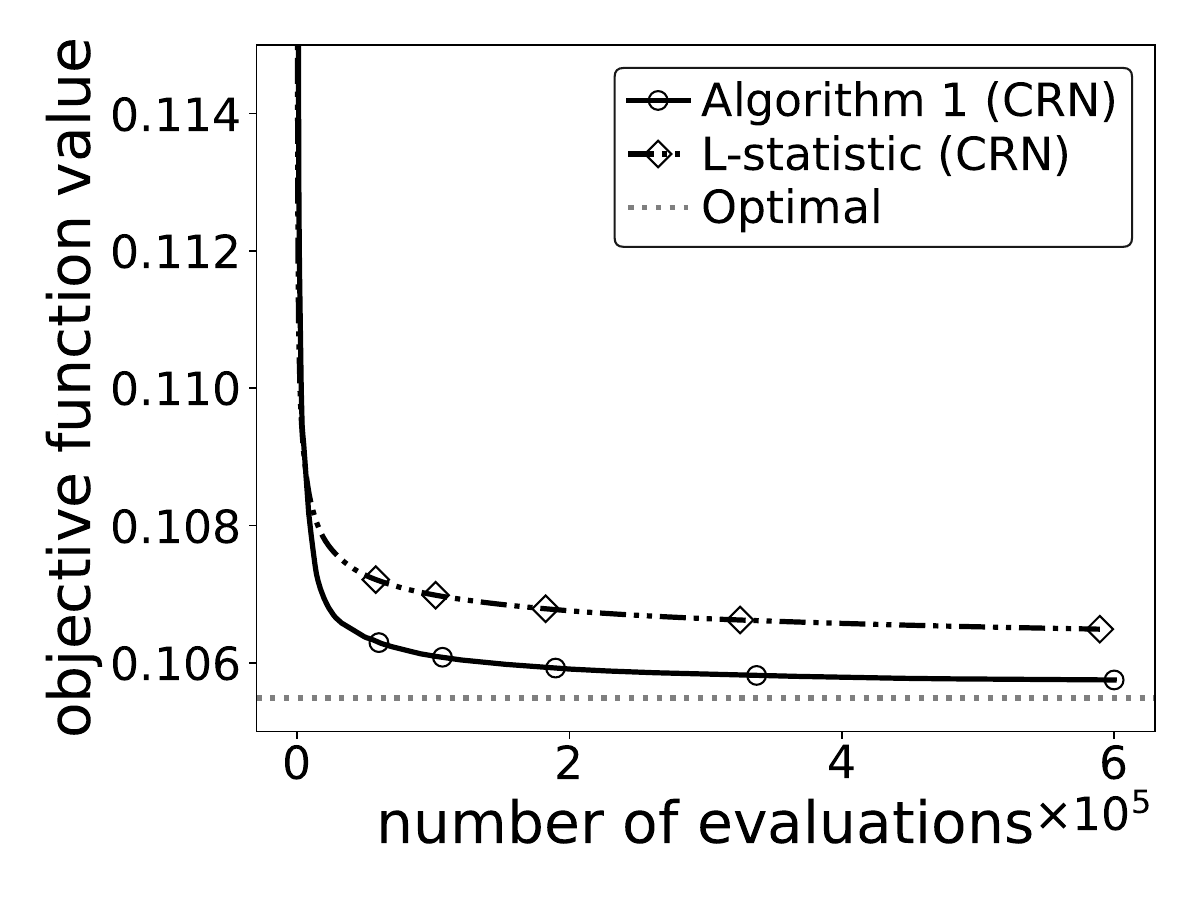}\par
            \vspace{3pt}\footnotesize (b) Objective (CRN)
        \end{minipage}\hfill
        \begin{minipage}{0.24\textwidth}
            \centering
            \includegraphics[width=\linewidth]{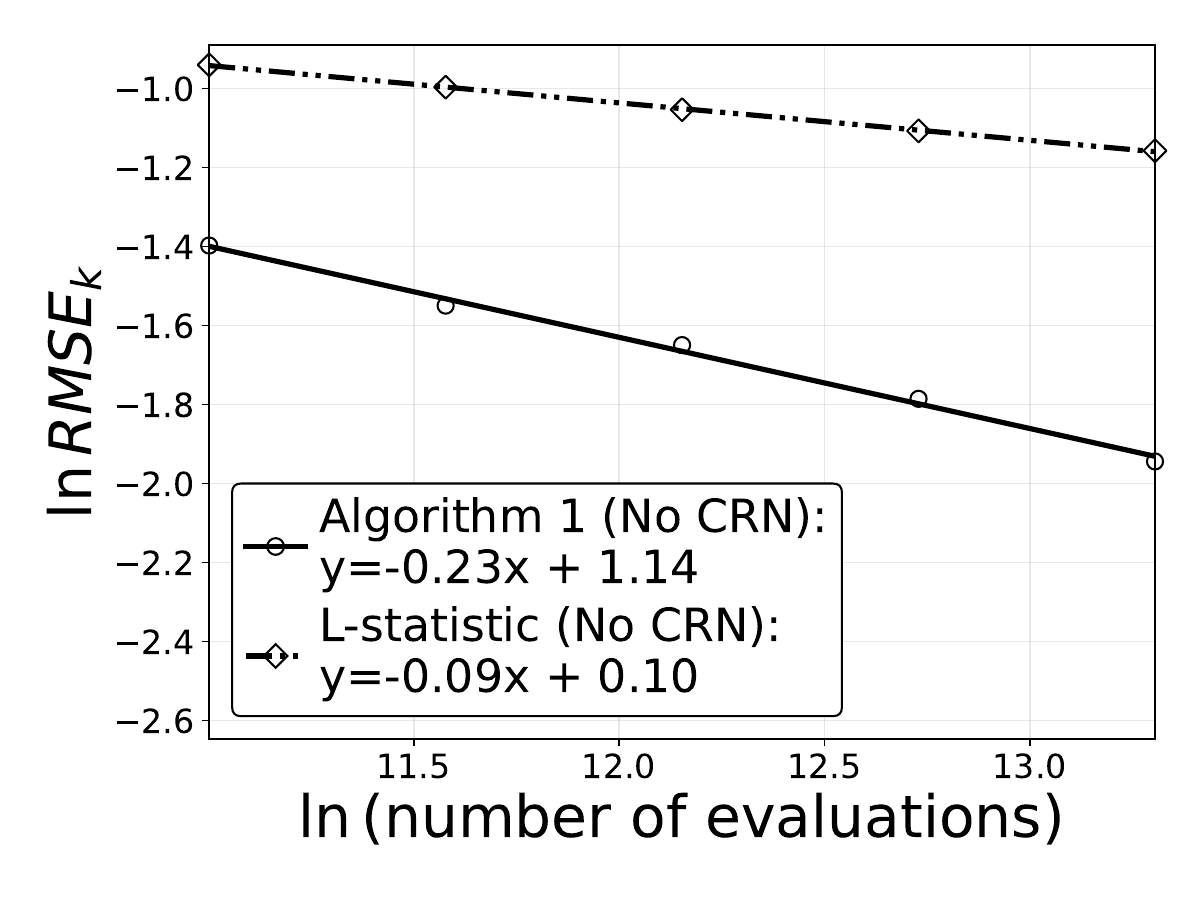}\par
            \vspace{3pt}\footnotesize (c) RMSE (No CRN)
        \end{minipage}\hfill
        \begin{minipage}{0.24\textwidth}
            \centering
            \includegraphics[width=\linewidth]{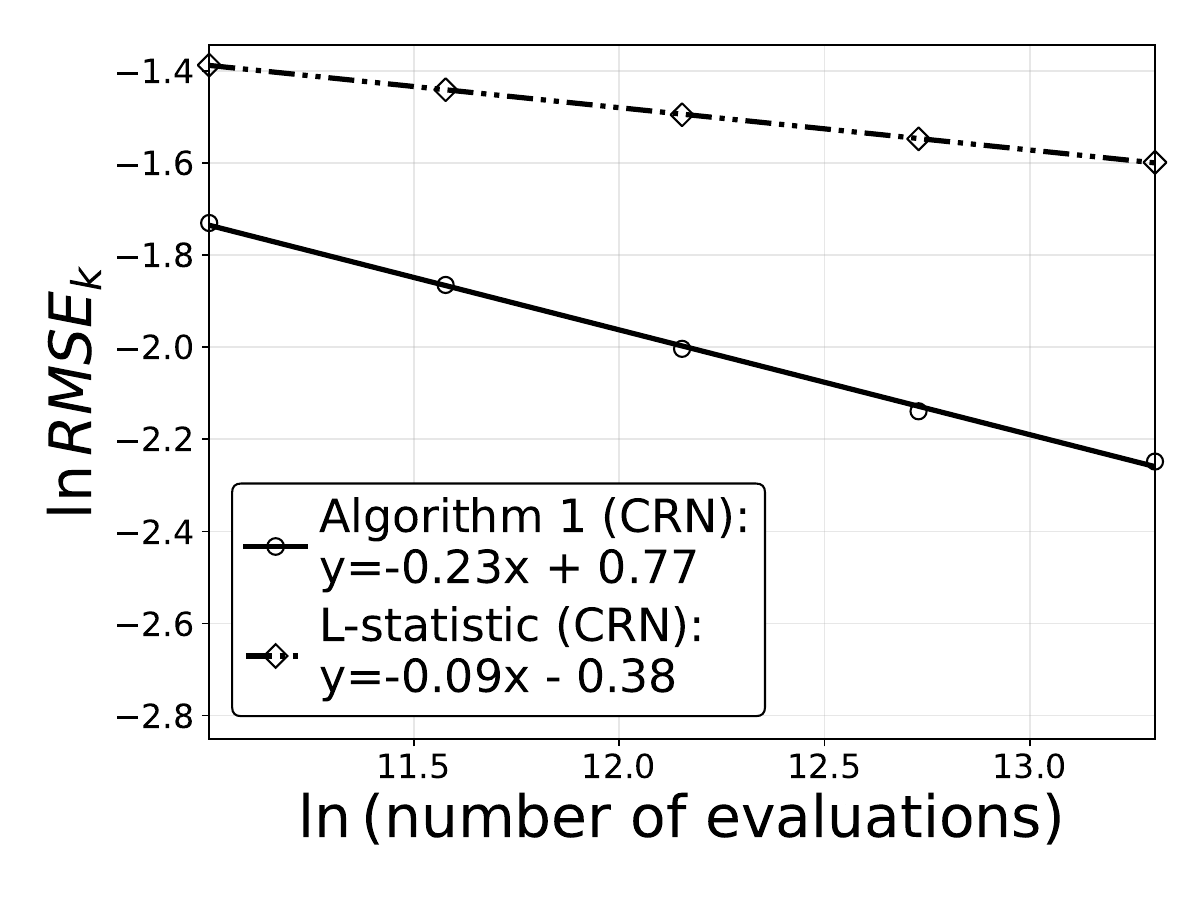}\par
            \vspace{3pt}\footnotesize (d) RMSE (CRN)
        \end{minipage}
    \end{minipage}%
}

\begin{figure}
     \FIGURE
    {\usebox{\boxPortfolioConvergence}}
    {Performance of Algorithm~\ref{alg} and the L-Statistic Benchmark on Portfolio Optimization.\label{fig:portfolio_convergence}}
    {Panels (a) and (b) show average objective values, and panels (c) and (d) show empirical decision RMSEs on logarithmic scales. The horizontal axis measures cumulative simulator evaluations.}
\end{figure}

\begin{table}
\centering
\TABLE
{Comparison Between Algorithm~\ref{alg} and the L-Statistic Method on Portfolio Optimization.\label{tab:portfolio_convergence}}
{\setlength{\tabcolsep}{4pt}\footnotesize
\begin{tabular*}{0.88\linewidth}{@{}l @{\extracolsep{\fill}} *{4}{c}@{}}
\toprule
& \multicolumn{2}{c}{No CRN} & \multicolumn{2}{c}{CRN} \\
\cmidrule(lr){2-3} \cmidrule(lr){4-5}
& Algorithm 1 & L-statistic & Algorithm 1 & L-statistic \\
\midrule
Avg. $\rho_g(\theta_K)$ & $\mathbf{0.106}$ & $0.108$ & $\mathbf{0.106}$ & $\mathbf{0.106}$ \\
RMSE of $\rho_g(\theta_K)$ & $\mathbf{0.001}$ & $0.003$ & $\mathbf{0.00039}$ & $0.002$ \\
RMSE of $\theta_K$ & $\mathbf{0.143}$ & $0.314$ & $\mathbf{0.106}$ & $0.202$ \\
\bottomrule
\end{tabular*}}
{\textit{Note.} Computational results based on 40 independent replications. 
All methods use a common budget of 600000 simulator outputs and the last completed round. Both methods are evaluated with and without common random numbers (CRN). Bold indicates the best displayed value within each CRN setting; ties are all bolded.}
\end{table}

\section{Conclusions}\label{sec_conclude}
In this paper, we have developed a unified approach for the estimation, sensitivity analysis, and optimization of DRMs in a black-box simulation setting. For piecewise polynomial distortion functions, we derive a reformulation of the DRM and its gradient that enables their estimation using only simulation outputs, without requiring explicit model knowledge or numerical partitioning. The resulting coupled recursions use a fixed number of new outputs per iteration, avoiding increasingly large simulation batches and repeated sorting. We establish strong consistency and convergence rates for the proposed method and extend the framework to non-polynomial distortion functions through polynomial approximation. Numerical experiments on synthetic and financial examples demonstrate the effectiveness of the approach in terms of both accuracy and computational efficiency.

The black-box nature of the proposed approach makes it applicable to complex stochastic systems that can be repeatedly evaluated through simulation.
Future research may explore variance reduction techniques, adaptive sampling strategies, and integration with online decision-making methods to further improve the computational efficiency and practical applicability of the proposed framework.

\bibliographystyle{informs2014}
\bibliography{references}

\end{document}